\documentclass[12pt,a4paper]{amsart}
\usepackage[utf8]{inputenc}%
\usepackage[T1]{fontenc}%
\usepackage[english]{babel}%
\usepackage[margin=2cm]{geometry}%
\usepackage{dsfont}%
\usepackage[pdftex]{graphicx}% Already loaded
\usepackage{amsmath,amssymb,amsthm}%
\usepackage{color}%
\usepackage{float}%
\usepackage{ifthen}%
\usepackage{hyperref}%
\usepackage{tikz}
\usepackage{pgfplots}
\usepackage{subcaption}
\usetikzlibrary{decorations.markings}
\pgfplotsset{compat=1.18}

\DeclareSymbolFont{Symbols}{OMS}{cmsy}{m}{n}
\DeclareMathSymbol{\Setminus}{\mathbin}{Symbols}{"6E}
\providecommand{\tightlist}{%
  \setlength{\itemsep}{2pt}\setlength{\parskip}{0pt}}% item-Zwischenraum
\newcounter{prob}% Numerate problems
\newcounter{exmpl}% Numerate examples
\DeclareMathOperator\Real{Re}

\def\grad{\boldsymbol{\nabla}}%
\def\q{\quad}%
\def\qq{\qquad}%
\def\qqq{\qq\q}%

\newcommand{\IN}{\mathbb{N}}% natural
\newcommand{\IP}{\mathbb{P}}% polynomials
\newcommand{\IR}{\mathbb{R}}% real
\newcommand{\cA}{\mathcal{A}}%
\newcommand{\cB}{\mathcal{B}}%
\newcommand{\cE}{\mathcal{E}}%
\newcommand{\cO}{\mathcal{O}}%
\newcommand{\cS}{\mathcal{S}}%
\def\bbT{\mathbb{T}}%
\def\bsF{\boldsymbol{F}}%
\def\bsf{\boldsymbol{f}}%
\def\bsG{\boldsymbol{G}}%
\def\bsH{\boldsymbol{H}}%
\def\bsM{\boldsymbol{M}}%
\def\bsm{\boldsymbol{m}}%
\def\bsn{\boldsymbol{n}}%
\def\bsP{\boldsymbol{P}}%
\def\bsp{\boldsymbol{p}}%
\def\bsq{\boldsymbol{q}}%
\def\bsV{\boldsymbol{V}}%
\def\bsv{\boldsymbol{v}}%
\def\bsw{\boldsymbol{w}}%
\def\bsx{\boldsymbol{x}}%
\def\bsy{\boldsymbol{y}}%
\def\bs0{\boldsymbol{0}}%
\def\bsPhi{\boldsymbol{\Phi}}%
\def\bsphi{\boldsymbol{\varphi}}%
\def\ez{\mathrm{e}}%
\def\de{\,\mathrm{d}}%
\def\fdg{\mid}%
\def\Id{\boldsymbol{\mathit{Id}}}%
\def\ptl{\partial}%
\def\pkt{\cdot}%
\def\Tri{\mathcal{K}}%
\def\tol{\mathrm{TOL}}
\def\tensor{\,\raise2pt\hbox{${}_{\otimes}$}}% Tensor
\newcommand{\LTE}{\text{LTE}}%
\newcommand{\ukappa}{\underline{\kappa}}

\def\Bmat#1{\begin{bmatrix}#1\end{bmatrix}}% Matrix in Displaystyle ...
\def\bmat#1{\left[\begin{smallmatrix}#1\end{smallmatrix}\right]}% Text
\def\dual#1#2{\langle#1\,,\,#2\rangle}
\def\norm#1#2{\|#1\|_{#2}}%
\renewcommand{\bar}{\overline}
\renewcommand{\hat}{\widehat}
\renewcommand{\tilde}{\widetilde}

\renewcommand{\leq}{\le}

\theoremstyle{plain}% default
\newtheorem{theorem}{Theorem}
\newtheorem{Proposition}[theorem]{Proposition}
\newtheorem{Lemma}[theorem]{Lemma}
\newtheorem{Corollary}[theorem]{Corollary}
\newtheorem{Algorithm}[theorem]{Algorithm}

\newtheorem{Remark}[theorem]{Remark}

\numberwithin{equation}{section}
\numberwithin{table}{section}
\numberwithin{figure}{section}

\def\OO{{\mathcal O}}

\newcommand{\unit}[1]{\,\mathrm{#1}}

\definecolor{pastelmagenta}{rgb}{0.96, 0.6, 0.76}% Art rot
\definecolor{carrotorange}{rgb}{0.93, 0.57, 0.13}% Art orange
\definecolor{kellygreen}{rgb}{0.3, 0.73, 0.09}% Art gruen
\definecolor{ashgrey}{rgb}{0.7, 0.75, 0.71}% Art grau
\title[Adaptive mesh refinement for the Landau--Lifshitz--Gilbert equation]{Adaptive mesh refinement for the\\[4pt]
Landau--Lifshitz--Gilbert equation}
\author{Willy D\"orfler, Michael Feischl, Stefan Karch} 

\begin{document}

\begin{abstract}
We propose a new adaptive algorithm for the approximation of the Landau--Lifshitz--Gilbert equation via a higher-order tangent plane scheme. We show that the adaptive approximation satisfies an energy inequality and demonstrate numerically that the adaptive algorithm outperforms uniform approaches. 
\end{abstract}

%\begin{keyword}
%\end{keyword}

\maketitle

\textbf{Keywords}\quad \def\sep{;\ }
adaptive method \sep BDF methods \sep energy technique \sep finite elements \sep Landau--Lifshitz--Gilbert equation \sep stability

%%%%%%%%%%%%%%%%%%%%%%%%%%%%%%%%%%%%%%%%%%%%%%%%%%%%%%%%%%%%%%%%%%%%%%%
\section{Introduction} \label{SEC:Intro}

The Landau--Lifshitz--Gilbert equation (LLG) serves as an important practical tool and as a valid model of micromagnetic phenomena occurring in, e.g., magnetic sensors, recording heads, and magneto-resistive storage devices~\cite{Gilbert:1955,LandauLifshitz:1935,Prohl:2001}.
Other developments, such as \emph{magneto-resistive memory} and \emph{skyrmion based race-track memory}, focus on new ways to store data magnetically, see, e.g.,~\cite{fert}.  
All of these applications have in common that magnetic materials exhibit a strong locality (in time and space) of magnetic effects~\cite{Schrefl:1999}:
The so-called domain walls reveal a change of direction of the magnetization on small length scales.
Similarly, temporal behavior is characterized by rapid changes (switching processes) followed by long relatively stable periods.
Those effects render nonadaptive simulations tremendously inefficient (see~\cite{Schrefl:1999} for early experiments and~\cite{adaptiveLLG1} for current applications) and classical methods based on uniform meshes often fail to reach meaningful accuracy even on the largest computers.

Even uniform time-stepping for the LLG equation is challenging.
Although weak convergence of the approximations has been known since at least 2008 (see, e.g.,~\cite{Alouges:2008,BartelsProhl:2006}), strong a~priori convergence of uniform time-stepping schemes that obey energy bounds has first been proved recently in~\cite{FeischlTran:2017_fembem} and then extended to higher-order in~\cite{AFKL:2021}.
The latter two works are built on the tangent plane idea first introduced in~\cite{Alouges:2008} to remove the nonlinear solver required in~\cite{BartelsProhl:2006}.
This is achieved by solving for the time derivative of the magnetization $\ptl_t\bsm$ instead of the magnetization itself.
The constraint $|\bsm|=1$ is translated into orthogonality $\ptl_t\bsm\cdot \bsm = 0$ almost everywhere.
In discretization, this can only hold nodewise~\cite{Alouges:2008} or in an average sense~\cite{AFKL:2021}.
As an attractive side effect, the discrete approximations still obey the energy decay, see also~\cite{BKW:2024}, where a similar approach is used for the harmonic map heat flow.
This is in contrast to other higher-order methods for LLG in~\cite{PRS:2018,KVBAT:2014,DPPRS:2020} (without strong error analysis) and~\cite{Gao:2014,An:2016} for which no discrete energy bounds were proved.
In \cite{ChengShen:2023} a projection algorithm is proposed that enforces energy decay, however, no error analysis is provided.

While the higher-order convergence depends on certain regularity assumptions on the exact solution (see~\cite{FeischlTran:2017_existence} for existence proofs of smooth solutions under a smallness condition on the initial data), it was proved in~\cite{AFP:2025} that even without regularity beyond $H^1$ one can still get weak convergence (a consequence of the energy bound) and thus is at least as good as traditional first-order schemes.

\smallskip

The aim of this work is to exploit the higher-order time-stepping schemes developed in~\cite{AFKL:2021} to obtain a heuristic estimate of the temporal approximation error.
This will allow us to construct an adaptive time-stepping scheme for LLG.
Combined with gradient recovery estimators for the spatial error, we derive a fully adaptive integrator for LLG.
While a convergence proof seems out of reach for the moment, we demonstrate in several experiments the effectiveness of the method.

In Section~\ref{SSEC:LLG} we introduce the initial value problem in its strong form.
In Section~\ref{SEC:disretisationLLG} we present the spatial discretization based on the standard higher-order conforming finite element method and a Lagrangian setting to cope with the nonlinear constraint.
For time-stepping, we use a collocation method based on the BDF($k$) method with a predictor based on extrapolation.
In the uniform setting, this method has been shown to converge with optimal rates in \cite{AFKL:2021}.
In Section~\ref{SEC:adaptivityTime} we formulate the BDF method with variable time steps.
Furthermore, we present approximations to higher-order derivatives with finite differences to estimate the truncation error of a given order.
This allows us to suggest new time step sizes and likewise orders.
In Theorem~\ref{THM:stability} we show that under regularity requirements, restricted time step changes, and a certain size of the damping term, the time adaptive method satisfies an energy bound.
In Sections~\ref{SEC:adaptivitySpace} and \ref{SEC:adaptivity} we formulate our adaptive method and conclude in Section~\ref{SEC:numerical-experiments} with some numerical experiments.

%%%%%%%%%%%%%%%%%%%%%%%%%%%%%%%%%%%
\subsection{Notation} \label{SSEC:notation}

Let $\Omega \subset\IR^d$, $d=2,3$, be a bounded domain with a polygonal boundary and $T>0$ a prescribed final time.
By $L^2(\Omega)$ we denote the usual Lebesgue spaces with scalar product $(v,w)_{\Omega} =\int_{\Omega}vw$ and norm $\norm{v}{L^2(\Omega)} =(v,v)_{\Omega}^{1/2}$.
Vector-valued functions are in $L^2(\Omega)^d$ if they are component wise in $L^2(\Omega)$ and their scalar product is defined as $(\bsv,\bsw)_{\Omega} =\int_{\Omega}\bsv\cdot\bsw$.
By $H^1(\Omega)$, we denote the space of functions in $L^2(\Omega)$ with weak derivatives in $L^2(\Omega)$.
The norm in $H^1(\Omega)$ is given by $(\norm{v}{L^2(\Omega)}^2+\norm{\grad v}{L^2(\Omega)}^2)^{1/2}$.

In the following, we work with vector-valued functions $\bsm =[m_1,m_2,m_3]:(0,T) \times\Omega\to\IR^3$
and use the notation $\bsm(t) =\bsm(t,\pkt):\Omega\to\IR^3$.
For $\grad\bsm$ define $\norm{\grad\bsm}{L^2(\Omega)}^2 =\int_{\Omega}\grad\bsm:\grad\bsm$, where '$:$' denotes the Frobenius product of two matrices.
For vectors or tensors $|\pkt|$ denotes the Euclidean norm or the Frobenius norm.

%%%%%%%%%%%%%%%%%%%%%%%%%%%%%%%%%%%
\subsection{The Landau--Lifshitz--Gilbert equation} \label{SSEC:LLG}

The strong form of the LLG equations we will use here is as follows: Find $\bsm:[0,T]\times\bar\Omega$ such that
\begin{alignat}{2}
   \label{EQ:LLG-general-pde}
   \alpha\ptl_t\bsm + \bsm\times\ptl_t\bsm
      &= \bsP(\bsm) \bsH_{\mathrm{eff}}(\bsm)\qq
      && \text{on }(0,T)\times\Omega\,,\\
   \label{EQ:LLG-general-bc-hom}
   \ptl_{\bsn}\bsm &= \bs0
      && \text{on }(0,T)\times\ptl\Omega\,,\\
   \label{EQ:LLG-general-iv}
   \bsm(0) &= \bsm^0 && \text{in }\Omega\,,
\end{alignat}
where $\bsn$ is the exterior normal vector field on $\ptl\Omega$ and $\ptl_{\bsn}\bsm =\bsn\cdot\grad\bsm$ the normal derivative. $\bsm^0: \Omega\to\IR^3$ is a given vector field of unit length, i.e.\ $|\bsm^0| =1$, and with $\ptl_{\bsn}\bsm^0 =\bs0$. $\alpha>0$ is a given parameter.
$\bsP$ is the projection operator defined by $\bsP(\bsm) =\Id-\frac{\bsm\tensor\bsm}{|\bsm|^2}$ for $\bsm \not=\bs0$.
In this work, we consider an effective field of the form
\begin{align}\label{EQ:LLG-Heff}
   \bsH_{\text{eff}}(\bsm)
    =  C_{\mathrm{e}} \Delta\bsm
     + \bsH_{\mathrm{ext}}
\end{align}
for some $C_{\mathrm{e}}>0$ and $\bsH_{\mathrm{ext}}$ a given time- and space-depending vector field.
Formally, from \eqref{EQ:LLG-general-pde}, by multiplying by $\bsm$, one can derive $\ptl_t|\bsm(t)|^2 =2\bsm\cdot\ptl_t\bsm =0$, hence $|\bsm(t)| =|\bsm^0| =1$.
Multiplying with $\ptl_t\bsm$ and integrating will lead to the stability bound
\begin{align}\label{EQ:stability}
   \cS(t) &= C_{\mathrm{e}}\norm{\grad\bsm(t)}{L^2(\Omega)}^{2}
      + \alpha\int_0^t
         \norm{\ptl_t\bsm(s)}{L^2(\Omega)}^{2}\de s\\
   \notag
   &\le C_{\mathrm{e}}\norm{\grad\bsm(0)}{L^2(\Omega)}^{2}
      + \frac1{\alpha}\int_0^t\norm{\bsH_{\mathrm{ext}}(s)}{L^2(\Omega)}^{2}
      \de s
\end{align}
for all $t\in[0,T]$.
The inequality \eqref{EQ:stability} is usually called \emph{energy bound}.
In the same way, but with additional integration by parts in time, one can derive conservation of the physical energy
\begin{align*}%\label{EQ:energy-conservation}
   \cE(t) + \alpha\int_0^t\norm{\ptl_t\bsm(s)}{L^2(\Omega)}^{2}\de s
      = \cE(0) - \int_0^t(\ptl_t\bsH_{\mathrm{ext}},\bsm)_\Omega,
\end{align*}
with $\cE$ defined by
\begin{align}\label{EQ:def_energy}
   \cE(t) = \frac12C_{\mathrm{e}}\norm{\grad\bsm(t)}{L^2(\Omega)}^{2}
      - (\bsH_{\mathrm{ext}},\bsm)_\Omega.
\end{align}
Note that the Landau--Lifshitz--Gilbert equation is historically
\begin{align*}
   \ptl_t\bsm - \alpha\bsm\times\ptl_t\bsm
      &= - \bsm\times \bsH_{\mathrm{eff}}(\bsm),
\end{align*}
but this is equivalent to \eqref{EQ:LLG-general-pde} and follows after multiplication with $\bsm\times$, using vector identities as well as $|\bsm| = 1$ and $\bsm\cdot\ptl_t\bsm = 0$.
The form \eqref{EQ:LLG-general-pde} has been used for the results in \cite{AFKL:2021}.

%%%%%%%%%%%%%%%%%%%%%%%%%%%%%%%%%%%%%%%%%%%%%%%%%%%%%%%%%%%%%%%%%%%%%%%
\section{Discretization of the Landau--Lifshitz--Gilbert equation} \label{SEC:disretisationLLG}

%%%%%%%%%%%%%%%%%%%%%%%%%%%%%%%%%%%
\subsection{Weak formulation} \label{SSEC:weakLLG}

The derivation of a weak form can be simplified if one uses the test functions $\bsphi \in C^\infty(\Omega)^3$ with $\bsm(t)\cdot\bsphi =0$ at time $t>0$, since it will allow us later to remove $\bsP(\bsm)$ from \eqref{EQ:LLG-general-pde}.
Thus, we define the solution-dependent constraint space by
\begin{align*}
   \bbT(\bsm(t)) = \big\{\bsphi \in H^1(\Omega)^3
      \fdg \bsm(t)\cdot\bsphi = 0\big\}\,.
\end{align*}
For ease of notation, we will omit the argument $t$ in the following and arrive formally at the following weak equation for $\bsm$
\begin{align*}
   &\alpha(\partial_t\bsm,\bsphi)_\Omega + (\bsm \times \partial_t\bsm,\bsphi)_\Omega
      + C_{\mathrm{e}}(\grad\bsm,\grad\bsphi)_\Omega
    = (\bsH_{\mathrm{ext}},\bsphi)_\Omega\qquad\text{for all }\bsphi\in\bbT(\bsm)\,.
\end{align*}
The existence of weak solutions $\bsm\in H^1(\Omega)^3$ is well known (see, e.g.,~\cite{AlougesSoyeur:1992}) and follows with a standard Galerkin--Faedo approximation and compactness arguments. The uniqueness of weak solutions is open in general and there are nonuniqueness results for the related harmonic map heat flow~\cite{harmonicmap,AlougesSoyeur:1992}. Note that in case one of the weak solutions is also a strong solution, we have uniqueness~\cite{weakstrong}. In terms of regularity, it is suspected that even smooth initial conditions can lead to the blow-up of $\nabla\bsm$ in $L^\infty(\Omega)^{3,3}$ in finite time, however, no proof is known. For LLG on smooth bounded domains with Neumann boundary conditions, smooth initial data close to constants lead to arbitrarily smooth solutions~\cite{FeischlTran:2017_existence}.

%%%%%%%%%%%%%%%%%%%%%%%%%%%%%%%%%%%
\subsection{Spatial discretization} \label{SSEC:spatialdiscretization}

Since $\Omega$ is assumed to be polygonally bounded, we let $\Tri_h$ be a regular triangulation such that $\bar\Omega =\bigcup\{K\in\Tri_h\}$.
On $\Tri_h$ we consider the finite element space $W_h \subset H^1(\Omega)$ of constant degree $p\ge 1$.
For vector fields, we then use $\bsV_h =W_h^3\subset H^1(\Omega)^3$.
It looks quite natural to correspondingly define the discrete constraint space by
\begin{align*}
   \bbT_h^\mathrm{s}(\bsm_h)
    = \big\{\bsphi_h\in\bsV_h \fdg \bsm_h\cdot\bsphi_h = 0
      \text{ at every node}\big\}
\end{align*}
for some $\bsm_h \in H^1(\Omega)^3$.
However, this natural choice leads to problems with higher-order convergence~\cite[Rem.~2.2]{AFKL:2021} and~\cite{AFKL:2021} proposes using a weak constraint in the form of
\begin{align*}
   \bbT_h(\bsm_h) = \bbT_h^\mathrm{w}(\bsm_h)
    = \big\{\bsphi_h\in\bsV_h \fdg
      (\bsm_h\cdot\bsphi_h,\psi_h)_{\Omega} = 0
      \text{ for all }\psi_h\in W_h\big\}\,.
\end{align*}
Then it is straightforward to define the semi-discrete solution $\bsm_h:[0,T]\to\bsV_h$ by
\begin{align*}
   \alpha(\partial_t\bsm_h,\bsphi_h)_\Omega
      + (\bsm_h\times\partial_t\bsm_h,\bsphi_h)_\Omega
      + C_{\mathrm{e}}(\grad\bsm_h,\grad\bsphi_h)_\Omega
   &= (\bsH_{\mathrm{ext}},\bsphi_h)_\Omega
      \q\text{for all }\bsphi_h\in\bbT_h(\bsm_h)\,.
\end{align*}

%%%%%%%%%%%%%%%%%%%%%%%%%%%%%%%%%%%
\subsection{Temporal discretization} \label{SSEC:temporaldiscretization}

We complete the discretization by defining a time-stepping method.
For this, let $0 =t_0<\dots<t_{N_T} =T$ be a decomposition of the interval $[0,T]$.
The finite element space will change with time $t_n$, $n=0,\dots,N_T$, due to regular refinement (and coarsening) of the preceding one and is denoted by $\bsV_h^n$.
We want to set up an equation for the approximation $\bsm^{n} \approx\bsm(t_n)$ at time $t_n$.
In addition, we introduce the unknown $\bsv^{n} \approx\ptl_t\bsm(t_n)$ and assume that we can impose a relation between $\bsv_h^n$ and $\bsm_h^{n-j}$ for a set of indices $j\ge 0$.
For example, if we assume a linear time dependence of $\bsm_h$, we can state $\bsm_h^n =\bsm_h^{n-1}+\tau_{n}\bsv_h^n$ for $\tau_n :=t_n-t_{n-1}$.
Thus, the resulting equation
\begin{align*}%\label{EQ:k-step-BDF-LLG}
   &\alpha(\bsv_h^n,\bsphi_h)_\Omega
      + (\bsm_h^n\times \bsv_h^n,\bsphi_h)_\Omega
      + C_{\mathrm{e}}(\grad\bsm_h^n,\grad\bsphi_h)_\Omega
   = (\bsH_{\mathrm{ext}}^n,\bsphi_h)_\Omega
   \qquad\text{for all }\bsphi_h\in\bbT_h(\bsm_h^n),
\end{align*}
where $\bsH_{\mathrm{ext}}^n = \bsH_{\mathrm{ext}}(t_n)$,
can be seen as a nonlinear equation for $\bsv_h^n$ from which $\bsm_h^n$ could be easily obtained.
However, we must note that the test space also depends on $\bsm_h^n$, and hence on $\bsv_h^n$. This approach is called \emph{tangent plane scheme}~\cite{Alouges:2008,HPPRSS:2019}.

The following linearization has proved to be useful: construct an approximation (or \emph{predictor}) $\hat{\bsm}_h^n$ to $\bsm_h^n$ and then define $\bsm_h^n$ by
\begin{align*}%\label{EQ:k-step-BDF-LLG}
   &\alpha(\bsv_h^n,\bsphi_h)_\Omega + (\hat{\bsm}_h^n\times \bsv_h^n,\bsphi_h)_\Omega
      + C_{\mathrm{e}}(\grad\bsm_h^n,\grad\bsphi_h)_\Omega
   = (\bsH_{\mathrm{ext}}^n,\bsphi_h)_\Omega
   \quad\text{for all }\bsphi_h\in\bbT_h(\hat{\bsm}_h^n)
\end{align*}
which is now a linear problem for $\bsv_h^n$.
There is still the problem of dealing with the constraint; for example, one might construct a local basis \cite[Sect.~2.2]{AFKL:2021}.
Here we will add the constraint as a separate weak equation.
Thus we introduce the Lagrange variable $\lambda_h^n \in W_h$ and state the saddle point problem
\begin{alignat}{2}
   \label{EQ:LLG-weak-lambda}
   \alpha(\bsv_h^n,\bsphi_h)_\Omega
      + (\hat{\bsm}_h^n\times \bsv_h^n,\bsphi_h)_\Omega
      + C_{\mathrm{e}}(\grad\bsm_h^n,\grad\bsphi_h)_\Omega&\ &&\\
   \notag
    {}+ (\hat{\bsm}_h^n\cdot\bsphi_h,\lambda_h^n)_\Omega
   &= (\bsH_{\mathrm{ext}}^n,\bsphi_h)_\Omega
   &&\qq\text{for all }\bsphi_h\in\bsV_h^n,\\
   \notag
   (\hat{\bsm}_h^n\cdot\bsv_h^n,\psi_h)_\Omega &= 0 
   &&\qq\text{for all }\psi_h\in W_h^n\,.
\end{alignat}
Note that the analysis of \cite{AFKL:2021} is based on this formulation.
It is straightforward to show that the discrete formulation is well posed and leads to a unique discrete solution $\bsv_h^n$.
All numerical results in this paper are based on this formulation.

%%%%%%%%%%%%%%%%%%%%%%%%%%%%%%%%%%%
\subsection{BDF time-stepping} \label{SSEC:BDF}

For time-stepping, we employ the BDF (\emph{backward differencing formula}) scheme of order $k \in\IN$ proposed in~\cite{AFKL:2021}.
A BDF scheme provides a general relation
\begin{align}\label{EQ:m-dotm-relation}
   \xi_k \bsm_h^n
   = \tilde\bsPhi_k(\bsm_h^{n-1},\dots,\bsm_h^{n-k})
     + \tau_n \bsv_h^n
\end{align}
for $n\ge k$, some number $\xi_k >0$ and an affine linear mapping $\tilde\bsPhi_k: (\IR^3)^k \to\IR^3$.
The case $k=1$ has been given as an example in the previous Section ~\ref{SSEC:temporaldiscretization}, more details will be provided in Section~\ref{SEC:adaptivityTime}.
We insert this relation into \eqref{EQ:LLG-weak-lambda} and get the following linear equation for $\bsv_h^n$
\begin{alignat}{2}
   \notag
   \alpha(\bsv_h^n,\bsphi_h)_\Omega
      + (\hat{\bsm}_h^n\times \bsv_h^n,\bsphi_h)_\Omega
      + \beta_k\tau_n(\grad\bsv_h^n,\grad\bsphi_h)_\Omega&\ &&\\
   \label{EQ:LLG-weak-lambda-bdf}
    {}+ (\hat{\bsm}_h^n\cdot\bsphi_h,\lambda_h^n)_\Omega
   &= \dual{\bsf^n}{\bsphi_h}_\Omega
   &&\qq\text{for all }\bsphi_h\in\bsV_h^n,\\
   \notag
   (\hat{\bsm}_h^n\cdot\bsv_h^n,\psi_h)_\Omega &= 0 
   &&\qq\text{for all }\psi_h\in W_h^n\,,
\end{alignat}
with $\beta_k =C_{\mathrm{e}}/\xi_k$ and the new right-hand side
\begin{align}\label{EQ:new-rhs}
   \dual{\bsf^n}{\bsphi_h}_\Omega
   := (\bsH_{\mathrm{ext}}^n,\bsphi_h)_\Omega
     - \beta_k(\grad\tilde\bsPhi_k(\bsm_h^{n-1},\dots,\bsm_h^{n-k}),\grad\bsphi_h)_\Omega\,.
\end{align}
$\bsm_h^n$ is then defined by \eqref{EQ:m-dotm-relation}. Note that a discrete solution $\bsv_h^n$ exists since the form is coercive on $\bbT_h(\hat{\bsm}_h^n)$.

As an example for a predictor $\hat\bsm_h^n$ one can take an extrapolation of order $k-1$, that is, find the unique polynomial $\bsq \in\IP_{k-1}$ with $\bsq(t_{n-1-j}) =\bsm_h^{n-1-j}$ for $j=0,\dots,k-1$, and set $\hat\bsm_h^n =\bsq(t_n)/|\bsq(t_n)|$.
This problem has been analyzed in \cite{AFKL:2021} and it is possible to prove, for uniform step size $\tau$, unconditional convergence in time up to order two \cite[Thm.~3.1]{AFKL:2021} and convergence up to order five in time with some restrictions on $\alpha$ \cite[Thm.~3.4]{AFKL:2021}.

%%%%%%%%%%%%%%%%%%%%%%%%%%%%%%%%%%%
\subsection{The uniform discretization} \label{SSEC:uniform-results}

For the moment, we consider a fixed (almost) uniform mesh of mesh size $h$ with a finite element method of polynomial order $p$ and in time BDF($k$) for $k=1,\dots,5$.
We will assume regularity of the weak solutions as in \cite[(3.2)]{AFKL:2021}, i.e.,
%\begin{Assumption} For $k,p\in\IN$ assume
\begin{align}\label{EQ:regularity_m}
\begin{split}
   \bsm \in &\ C^{k+1}([0,T],L^{\infty}(\Omega)^3)
      \cap C^{1}([0,T],W^{p+1,\infty}(\Omega)^3)\,,\\
   &\Delta\bsm+\bsH_{\mathrm{ext}}
      \in C^{0}([0,T],W^{p+1,\infty}(\Omega)^3)\,.
\end{split}
\end{align}
%\end{Assumption}
Then it has been proved in \cite{AFKL:2021} that for sufficiently small $\tau$ and $h$ the discrete solutions exist and it holds
\begin{align}\label{EQ:error-bound}
   \norm{\grad(\bsm(t_n)-\bsm_h^n)}{L^2(\Omega)}
      \le C(\tau^k+h^p)
\end{align}
for $t_n =n\tau \le T$.
In case $k\in\{1,2\}$ one needs the additional condition $\tau \le C\sqrt{h}$ \cite[Thm.~3.1]{AFKL:2021}, while in case $k\in\{3,4,5\}$ one requires $p\ge 2$, $\tau \le Ch$, and $\alpha\ge\alpha_k>0$ for some constant $C$ and some threshold values $\alpha_k$ \cite[Thm.~3.4]{AFKL:2021}.
As a consequence of this result it is further concluded that normalization is approximated as (cf.~\cite[Rem.~3.1]{AFKL:2021})
\begin{align*}%\label{EQ:norm-bound}
   \norm{1-|\bsm_h^n|\,}{L^2(\Omega)}
      \le C(\tau^k+h^p).
\end{align*}

%%%%%%%%%%%%%%%%%%%%%%%%%%%%%%%%%%%%%%%%%%%%%%%%%%%%%%%%%%%%%%%%%%%%%%%
\section{Adaptivity in Time}\label{SEC:adaptivityTime}

%%%%%%%%%%%%%%%%%%%%%%%%%%%%%%%%%%%
\subsection{The variable step BDF method}\label{SSEC:variable-BDF}\ %

To describe the BDF method of order $k$ for $k\in\IN$ we find it convenient to work on the grid
\begin{align*}
   [0,t_1,\dots,t_k] = [0,\tau_1,\tau_1+\tau_2,\tau_1+\dots+\tau_k]
\end{align*}
for positive step sizes $\tau_j$, $j=1,\dots,k$.
This can then be translated into a mesh of the form $t_{n-k},\dots,t_n$ (for $n\ge k$) by an index shift.
We will use abbreviations
\begin{align*}
   \tau_{12} = \tau_1+\tau_2,\ \tau_{23} = \tau_2+\tau_3,
   \ \tau_{122} = \tau_1+2\tau_2,\ \dots
\end{align*}
and let $\tau =\max\{\tau_1,\dots,\tau_k\}$.
The definition of a variable step BDF method to solve the ordinary differential equation (ODE) $\bsy' =\bsF(t,\bsy)$, $\bsy: [0,T]\to\IR^N$, $\bsF: [0,T]\times\IR^N\to\IR^N$, is as follows:
Given $\{\bsy^j\}_{j=0}^{k-1}$, determine $\bsy^k$ (approximating $\bsy(t_k)$) such that there is a polynomial $\bsq\in\IP_k^m$ with
\begin{align*}
   \bsq(t_j) = \bsy^j,\ j=0,\dots,k,\qq \bsq'(t_k)
      = \bsF(t_k,\bsy^k)\,.
\end{align*}
This leads to a nonlinear equation of the form
\begin{align}\label{EQ:BDF-general}
    \bsPhi_k(\bsy)
    = \xi_k \bsy - \tilde\bsPhi_k(\bsy^{k-1},\dots,\bsy^0)
       - \tau_k \bsF(t_k,\bsy) = \bs0
\end{align}
to determine $\bsy^k$. 
For $k =1,2,3$ we obtain the explicit formulae \cite{CAC:2014},~\cite[App.~A]{DGLL:2021}
\begin{align}
   \label{EQ:def_Phi1}
   \bsPhi_1(\bsy) &= \bsy-\bsy^0-\tau_1\bsF(t_1,\bsy)\,,\\
   \label{EQ:def_Phi2}
   \bsPhi_2(\bsy)
   &= \frac{\tau_{122}}{\tau_{12}}\bsy - \frac{\tau_{12}}{\tau_1}\bsy^1
        + \frac{\tau_2^2}{\tau_1\tau_{12}}\bsy^0 
        - \tau_2\bsF(t_2,\bsy)\,,\\
  \label{EQ:def_Phi3}
  \bsPhi_3(\bsy)
   &= \frac{\tau_{23}\tau_{1233} + \tau_{123}\tau_3}{\tau_{23}\tau_{123}} \bsy
      - \frac{\tau_{123}\tau_{23}}{\tau_{12}\tau_{2}} \bsy^2
      + \frac{\tau_{123}\tau_{3}^2}{\tau_{1}\tau_{2}\tau_{23}} \bsy^1
      - \frac{\tau_{23}\tau_{3}^2}{\tau_{1}\tau_{12}\tau_{123}} \bsy^0
      - \tau_3\bsF(t_3,\bsy)\,.
\end{align}
A formula for $\bsPhi_4$ is for example provided in \cite[App.~A]{DGLL:2021}.
If we compare this with \eqref{EQ:m-dotm-relation}, we can get explicit values for $\xi_k$ and expressions for $\tilde\bsPhi_k(\bsy^{k-1},\dots,\bsy^0)$.
These formulas are in agreement with those for the constant $\tau$ \cite[Ch.~V]{HairerWanner:2010}.
Recall that for $k=1$ we get the implicit Euler method.

We will, in the provided notation and for later use, formulate consistent approximations for second- to fourth-order derivatives
\begin{align}
   \label{EQ:def_dh2y}
   \ptl_{\tau}^{2}\bsy(t_2)
   &= 2\Big(\frac1{\tau_{12}\tau_2}\bsy_2
      - \frac1{\tau_1\tau_2}\bsy_1
      + \frac1{\tau_1\tau_{12}}\bsy_0\Big)\,,\\
   \label{EQ:def_dh3y}
   \ptl_{\tau}^{3}\bsy(t_3)
   &= 6\Big(\frac1{\tau_{123}\tau_{23}\tau_3}\bsy_3
      - \frac1{\tau_{12}\tau_2\tau_3}\bsy_2
      + \frac1{\tau_1\tau_2\tau_{23}}\bsy_1
      - \frac1{\tau_1\tau_{12}\tau_{123}}\bsy_0\Big)\,,\\
  \label{E:def_dh4y}
  \ptl_{\tau}^{4}y(t_4)
  &= 24\Big(\frac1{\tau_{1234}\tau_{234}\tau_{34}\tau_4}\bsy_4
     - \frac1{\tau_{123}\tau_{23}\tau_3\tau_4}\bsy_3
     + \frac1{\tau_{12}\tau_2\tau_3\tau_{34}}\bsy_2\\
  \notag
  &\qq\q{} - \frac1{\tau_1\tau_2\tau_{23}\tau_{234}}\bsy_1
     + \frac1{\tau_1\tau_{12}\tau_{123}\tau_{1234}}\bsy_0\Big).
\end{align}

%%%%%%%%%%%%%%%%%%%%%%%%%%%%%%%%%%%
\subsection{Time step selection}\label{SSEC:stepselection}

To choose the next time step, we will always start with the first-order method, i.e.\ the implicit Euler method ($\bsPhi_1$).
Since the local truncation error $\LTE$ at time $t$ is known to be $1/2\,|\bsy''(t)|\tau^2$, we suggest for the next time step
\begin{align}\label{EQ:def_tau-1}
   \tau_n^{[1]} := \Big(\frac{2\,\tol}{|\bsy''(t_n)|}\Big)^{1/2}\,,
\end{align}
with the idea to get $\tol \ge \LTE \approx \frac12 |\bsy''(t_n)|\tau_n^2$.
In order to evaluate this, we use
\begin{align}\label{EQ:d2y}
   \bsy''(t_n) = \ptl_t\bsF(t_n,\bsy(t_n))
      + \ptl_{\bsy}\bsF(t_n,\bsy(t_n))\bsF(t_n,\bsy(t_n))\,,
\end{align}
but alternatively, we may use the finite difference approximation $\ptl_{\tau}^{2}\bsy(t_n)$ \eqref{EQ:def_dh2y} for $n\ge 2$.

The next order method is given by $\bsPhi_2$ and the truncation error in the uniform case is known to be $1/3\,|\bsy'''(t)|\tau^3$.
In the case of variable step sizes we use the Taylor expansion of $\bsPhi_2(\bsy(t))$ in $t_n$ and get the local truncation error $1/6\,\tau_{n-1,n}\tau_n^2\,|\bsy'''(t_n)|$ which inspires the time step prediction
\begin{align}\label{EQ:def_tau-2}
   \tau_n^{[2]} := \Big(\frac{6\,\tol}{\tau_{n-1,n}|\bsy'''(t_n)|}\Big)^{1/2}\,.
\end{align}
At this point, we will rely on the finite difference approximation $\ptl_{\tau}^{3}\bsy(t_n)$ \eqref{EQ:def_dh3y} for $n\ge 3$.
For the third-order method, we utilize the time step prediction
\begin{align}\label{EQ:def_tau-3}
   \tau_n^{[3]} := \Big(\frac{24\,\tol}{\tau_{n-2,n-1,n}\tau_{n-1,n}| \bsy^{(4)}(t_n)|}\Big)^{1/2}\,.
\end{align}
Note that in this work we do not need to compute $k+1$-order finite difference approximations for BDF($k$), since we can take advantage of the fact that we have already solved for the time derivative and use the approximations
$\norm{\ptl_t^{k+1}\bsm(t_n)}{L^2(\Omega)} \approx\norm{\ptl_\tau^{k+1}\bsm_h(t_{n})}{L^2(\Omega)} \approx\norm{\ptl_\tau^{k}\bsv_h(t_{n})}{L^2(\Omega)}$.

Each time step choice for BDF($k$) will be subject to a restriction compared to the previous step in the form
\begin{align}\label{eq:tn}
  \tau_{n+1} = 
  \begin{cases}
     \max\big\{\tau_{\min},\frac14\tau_n,\min\{\tau_n^{[k]},
       2\,\tau_n,\tau_{\max}\}\big\}\qquad\quad &\text{for }k\le 2\,,\\[4pt]
    \max\big\{\tau_{\min},\frac12\tau_n,\min\{\tau_n^{[k]},
       \sqrt{2}\,\tau_n,\tau_{\max}\}\big\}\qquad &\text{for }k>2
   \end{cases}
\end{align}
that is slightly above the bound in Theorem \ref{THM:stability}.
Positive numbers $\tau_{\min}$ and $\tau_{\max}$ are prescribed lower and upper bounds for the size of the time steps, here taken to be $10^{-10}$ and $T/10$.
Note that derivatives $|\bsy^{(k)}(t_n)|$ close to or equal to zero are handled by the maximum cutoff value $\tau_{\max}$ in~\eqref{eq:tn}.
In case we allow order selection, the order $k_{n+1} \in\{\max\{1,k_n-1\}, k_n,\min\{k_n+1,k_{\max}\}\}$ is chosen that provides the maximal time step, where $k_{\max}$ is a prescribed maximal order.
There exist many variable step size (also variable order) BDF algorithms, for example \cite{ShampineReichelt:1997}.
Our method is a variant of the one proposed in \cite{CAC:2014}.

%%%%%%%%%%%%%%%%%%%%%%%%%%%%%%%%%%%
\subsection{Application to LLG}\label{SSEC:time steps-LLG}%

We iteratively solve \eqref{EQ:LLG-weak-lambda-bdf}, starting from a normalized initial condition $\bsm_h^0$ that approximates $\bsm^0$.
The BDF(2) method with variable step size is then defined in Sections~\ref{SEC:disretisationLLG}, \ref{SSEC:variable-BDF}, and \ref{SSEC:stepselection} in case $k=2$, taking $L^2(\Omega)$-norms for the difference quotients in time.
Then we use an extrapolation of order $k-1$ to define the prediction $\hat{\bsm}_h^n$ \cite[(2.1)]{AFKL:2021}.
An important statement about the uniform BDF(2) method for LLG is the energy bound \cite[Prop.~3.1]{AFKL:2021}.
We can derive a corresponding result for BDF(2) with variable step size under the condition that the coarsening factor is limited and $\alpha$ has a positive lower bound.

There are results that prove stability for larger step size relations for $\tau_n/\tau_{n-1}$ than $1+\sqrt{2}$, e.g., \cite{LiaoZhang:2021} \cite{LTZ:2020} for diffusion equations.
However, these results do not directly apply to our equations. Nevertheless, as optimal error estimates are concerned, the recommended step size relation can still be the classical result of \cite{Grigorieff:1983}, see also \cite[p.~1223]{LiaoZhang:2021}.

Besides the bounds $\ukappa\le\kappa_n\le\bar\kappa$ on the step ratios $\kappa_n=\tau_n/\tau_{n-1}$, Theorem~\ref{THM:stability} below requires that the step ratios vary smoothly, see \eqref{EQ:assumpt-smooth-kappa}. 
Restrictions on the variation of the step sizes are known in the stability analysis of variable step multistep methods.
Gear and Tu \cite{GearTu:1974} and, for variable order, Gear and Watanabe \cite[Thm.~2.5]{GearWatanabe:1974} prove stability in the sense of \cite[Def.~2]{GearWatanabe:1974} for small step changes $\tau_{n+1}/\tau_n = 1+\OO(\tau_n)$.
This implies $|\kappa_{n+1}-\kappa_n| = \OO(\tau_n)$ and hence \eqref{EQ:assumpt-smooth-kappa}.
However, our assumption also allows step ratios away from $1$.
In practice, for the time step selection of Section~\ref{SSEC:stepselection} based on estimates of derivatives, smoothly varying step ratios are expected for smooth solutions \cite[Sect.~5]{GearTu:1974}.
Moreover, $|\kappa_{n+1}-\kappa_n|/\tau_n$ can be tracked throughout the computation.

\begin{theorem}[Energy bound for orders $k =1,2$]
\label{THM:stability}
Consider the discretization \eqref{EQ:LLG-weak-lambda-bdf} of the LLG equation \eqref{EQ:LLG-general-pde} for $k \in\{1,2\}$ with finite elements of polynomial degree $p \ge 1$.
We assume $\tau =\mathcal{O}(h)$ and, for $k=2$, that the time steps satisfy $0 < \ukappa \le \kappa_n := \tau_n/\tau_{n-1} \le \bar \kappa$ with fixed $\ukappa \le \bar \kappa < \kappa_0 = (2 + \sqrt{13})/3 = 1.868\dots$ and $\alpha \ge\alpha_2 >0$ for some $\alpha_2$ that depends on $\bar \kappa$.
Further, for $k=2$ we assume that the step-ratios vary smoothly in the sense that
\begin{align}\label{EQ:assumpt-smooth-kappa}
    | \kappa_{n+1} - \kappa_n | \le c_\kappa \tau_n
\end{align}
for all $n$ and some fixed constant $c_\kappa > 0$.
Then, the numerical solution satisfies the following discrete energy bound: For $n >k$ it holds
\begin{align*}
   \norm{\grad\bsm_h^n}{L^2(\Omega)}^2
   + \frac{\alpha}{2}\sum_{j=k}^{n}\tau_j\norm{\bsv_h^j}{L^2(\Omega)}^2
   \le C\Big(\sum_{i=0}^{k-1}\norm{\grad\bsm_h^i}{L^2(\Omega)}^2
   + \sum_{j=k}^{n}\tau_j
   \norm{\bsH_{\mathrm{ext}}(t_j)}{L^2(\Omega)}^2\Big)
\end{align*}
under the regularity requirements \eqref{EQ:regularity_m} and for some positive constants $C$ that is independent of the spatial mesh size $h$ and the time steps $\tau$, but depends exponentially on $T$ for $k=2$.
\end{theorem}

\begin{proof}
We first note that the BDF($k$) scheme provides coefficients $\delta_j$, $j=0,\dots,k$, such that
\begin{align}
    \label{EQ:BDFvh}
   \tau_n\bsv_h^n = \sum_{j=0}^{k}\delta_j\bsm_h^{n-j}\,,
\end{align}
see Section~\ref{SSEC:variable-BDF}.
If we let
\begin{align*}
   \delta(s) := \sum_{i=0}^k \delta_i s^i\,,
\end{align*}
and if we can show that for some $\eta>0$ we have $\Real(\delta(s)/(1-\eta s)) >0$ for all complex $s$ with $|s|<1$, then the BDF method is G-stable~\cite[Lem.~A.2]{AFKL:2021}.
Such an $\eta$ is called \emph{multiplier}.
Then we test the weak equation for $\bsm_h^n-\eta\bsm_h^{n-1}$ with $\bsv_h^n$ and will obtain a lower bound for $(\grad(\bsm_h^n-\eta\bsm_h^{n-1}),\grad\bsv_h^n)_{\Omega}$ due to the G-stability, see \cite{AFKL:2021} for the uniform case.

The case $k=1$ is clear.

\textbf{Step~1, Multiplier:}
In the case $k =2$, we let $\tau_n =\kappa\tau$ and $\tau_{n-1} =\tau$ for some $\kappa >0$. From \eqref{EQ:def_Phi2} we get
\begin{align*}
    \delta_0 = \frac{1+2\kappa}{1+\kappa}\,,\quad
    \delta_1 = -(1+\kappa)\,,\quad
    \delta_2 = \frac{\kappa^2}{1+\kappa}\,.
\end{align*}
We immediately see that $\delta$ has two zeros $s_1 =1$ and $s_2 =(1+2\kappa)/\kappa^2 >1$ as long as $\kappa <1+\sqrt{2}$.
Thus, we have
\begin{align*}
   \frac{\delta(s)}{1-\eta s}
      = \frac{\delta_2(s-1)(s-s_2)}{1-\eta s}
      = \frac{\delta_2(s-1)(s-s_2)(1-\eta\overline{s})}{|1-\eta s|^2}\,.
\end{align*}
Clearly, the choice $\eta =1/s_2$ would suffice to ensure $\Real(\delta(s)/(1-\eta s)) >0$ for $|s|<1$.
However, we obtain a sharper bound by resolving the inequality
\begin{align*}
   \Real\big((s-1)(s-s_2)(1-\eta\overline{s})\big) > 0
      \qq\text{for }|s|<1\,.
\end{align*}
Considering the edge case $|s|=1$ leads to the following quadratic inequality for the real part $s_r$ of $s$ (the cubic terms cancel)
\begin{align*}
   2s_r^2 - (1+\eta)(1+s_2)s_r + (1+\eta)s_2 + \eta -1 > 0
   \qq\text{for all } -1< s_r< 1\,.
\end{align*}
For given $0 <\kappa <1+\sqrt{2}$, this can be fulfilled if
\begin{align*}
   \eta>\frac{3-s_2}{s_2+1}
   = g(\kappa) := \frac{3\kappa^2-2\kappa-1}{\kappa^2+2\kappa +1}\,.
\end{align*}
For $\kappa \le 1$ we can take $\eta =0$ and for $\kappa \in(1,1+\sqrt{2})$ we find $\eta <1$.
The bound $\kappa <1+\sqrt{2}$ dates back to \cite{Grigorieff:1983}, but was proved here with the multiplier technique from \cite{NevanlinnaOdeh:1981}.
Moreover, we have $1 - \frac{1}{2}g(\kappa) ( 1 + \kappa) > 0$ for $\kappa \in (0, (2 + \sqrt{13})/3) =(0, \kappa_0)$.
Thus for any $\kappa<\kappa_0$ there exists $\omega >0$ with $1-\frac12 g(\kappa)(1+\kappa)>2\omega$.
Since $g$ is strictly increasing,~cf.~Lemma~\ref{LEM:gmono}, and $\bar\kappa<\kappa_0$, the interval $\big(g(\bar\kappa),g(\kappa_0)\big)$ is nonempty, and we fix a single multiplier
\begin{align*}
    \eta = \eta_\star\in\big(g(\bar\kappa),\,g(\kappa_0)\big)\qq
       \textrm{with } \eta_\star \ge 0,
\end{align*}
valid for all step-ratios $\kappa=\kappa_n\in[\ukappa,\bar\kappa]$.
For convenience, for $k=2$ we define $\bsG_n := \bsG(\kappa_n, \eta_\star)$, where $\bsG$ is provided in Lemma~\ref{LEM:Gconstruction} in the appendix, and denote by $\gamma_-, \gamma_+$ the eigenvalues of $\bsG_n$. 
By Proposition~\ref{PROP:Lipschitz} with $K = [\ukappa, \bar \kappa]$, we have $0 < \gamma_-^{\min} \le \gamma_- \le \gamma_+ \le \gamma_+^{\max}$.
For $k=1$ we take $\eta = 0$ and $\bsG_n = \frac12$. 
Then $\bsG$ does not depend on $n$ and the related sum in Step~2 vanishes.

%%%%%%%%%%%%%%%%%%%%%%%%%%%%%%%%%%%%%%%%%%%%%%%%%%%%%%%%%
\textbf{Step~2, General energy bound:} Let us first take, for simplicity, $\bsH_{\mathrm{ext}}^n =\bs0$.
We proceed recalling the weak equations for $\bsv_h^n$ and $\bsv_h^{n-1}$, for $n\ge k$,
\begin{alignat*}{2}%\label{}
   \alpha(\bsv_h^n,\bsphi_h)_\Omega
      + (\hat{\bsm}_h^n\times \bsv_h^n,\bsphi_h)_\Omega
      + C_{\mathrm{e}}(\grad\bsm_h^n,\grad\bsphi_h)_\Omega
   &= 0\quad
   && \text{for all }\bsphi_h\in\bbT_h(\hat{\bsm}_h^n)\,,\\
   \alpha(\bsv_h^{n-1},\bsphi_h)_\Omega
      + (\hat{\bsm}_h^{n-1}\times \bsv_h^{n-1},\bsphi_h)_\Omega
      + C_{\mathrm{e}}(\grad\bsm_h^{n-1},\grad\bsphi_h)_\Omega
   &= 0\quad
   && \text{for all }\bsphi_h\in\bbT_h(\hat{\bsm}_h^{n-1})\,.
\end{alignat*}
In order to derive an equation for $\bsv_h^n-\eta\bsv_h^{n-1}$ we use the test function $\bsv_h^n \in\bbT_h(\hat{\bsm}_h^n)$ in the first equation and $\hat{\bsP}_h^{n-1}\bsv_h^n$ with $\hat{\bsP}_h^{n-1} :=\bsP_h(\hat{\bsm}_h^{n-1}) \in\bbT_h(\hat{\bsm}_h^{n-1})$, where $\bsP_h$ denotes the $L^2$-orthogonal projection onto $\bbT_h$, in the second equation.
The latter term will be written in the form $\hat{\bsP}_h^{n-1}\bsv_h^n =\bsv_h^n-(\hat{\bsP}_h^n-\hat{\bsP}_h^{n-1})\bsv_h^n =\bsv_h^n-\bsp_h^n$.
Incorporation of these test functions results in
\begin{align*}%\label{EQ:}
   &\alpha(\bsv_h^n,\bsv_h^n)_\Omega
      + C_{\mathrm{e}}(\grad\bsm_h^n,\grad\bsv_h^n)_\Omega
   = 0\,,\\
   &\alpha(\bsv_h^{n-1},\bsv_h^n)_\Omega
      + C_{\mathrm{e}}(\grad\bsm_h^{n-1},\grad\bsv_h^n)_\Omega\\
   &\qqq= \alpha(\bsv_h^{n-1},\bsp_h^n)_\Omega
      - (\hat{\bsm}_h^{n-1}\times\bsv_h^{n-1},\bsv_h^n-\bsp_h^n)_\Omega
      + C_{\mathrm{e}}(\grad\bsm_h^{n-1},\grad\bsp_h^n)_\Omega\\
   &\qqq= \alpha(\bsv_h^{n-1},\bsp_h^n)_\Omega
      - (\hat{\bsm}_h^{n-1}\times \bsv_h^{n-1},
         \bsv_h^n-\sigma\bsv_h^{n-1}-\bsp_h^n)_\Omega
      + C_{\mathrm{e}}(\grad\bsm_h^{n-1},\grad\bsp_h^n)_\Omega\,,
\end{align*}
where we can choose either $\sigma =0$ or $\sigma =1$.
From these two equations, we derive
\begin{align}\label{EQ:vn-minus-eta_vn-1}
\begin{split}
   &\alpha(\bsv_h^n-\eta\bsv_h^{n-1},\bsv_h^n)_\Omega
      + C_{\mathrm{e}}(\grad(\bsm_h^n-\eta\bsm_h^{n-1}),
         \grad\bsv_h^n)_\Omega\\
   &\quad= -\eta\alpha(\bsv_h^{n-1},\bsp_h^n)_\Omega
    + \eta(\hat{\bsm}_h^{n-1}\times \bsv_h^{n-1},
       \bsv_h^n-\sigma\bsv_h^{n-1}-\bsp_h^n)_\Omega
    - \eta C_{\mathrm{e}}(\grad\bsm_h^{n-1},\grad\bsp_h^n)_\Omega\,.
\end{split}
\end{align}
By construction of $\eta$ in the first part of the proof and applying Lemma~\ref{LEM:Gconstruction}, we have for the left-hand side of \eqref{EQ:vn-minus-eta_vn-1} the lower bound \cite[Lem.~8.1(proof)]{AFKL:2021}
\begin{align*}%\label{EQ:}
   &\alpha(\bsv_h^n-\eta\bsv_h^{n-1},\bsv_h^n)_\Omega
      + C_{\mathrm{e}}(\grad(\bsm_h^n-\eta\bsm_h^{n-1}),
         \grad\bsv_h^n)_\Omega\\
   &\qq\ge \big(1-\frac{\eta}2\big)\alpha\norm{\bsv_h^n}{L^2(\Omega)}^2
      - \frac{\eta}2\,\alpha\norm{\bsv_h^{n-1}}{L^2(\Omega)}^2
      + \frac{C_{\mathrm{e}}}{\tau_n}
          \big(\norm{\grad\bsM_h^n}{\bsG_n,\Omega}^2
        - \norm{\grad\bsM_h^{n-1}}{\bsG_n,\Omega}^2\big)\,,
\end{align*}
where $\bsM_h^n =[\bsm_h^{n-k+1}, \dots, \bsm_h^n]$ and 
\begin{align*}
    \norm{\grad \bsM_h^n}{\bsG_n, \Omega}
    := \sum\limits_{i, \ell = 1}^k g_{i\ell}
       (\grad \bsm_h^{n-k+i}, \grad \bsm_h^{n-k+\ell})_{\Omega}.
\end{align*}
for $\bsG_n \in \IR^{k,k}$.
For details, see \cite[Appendix]{AFKL:2021}.
In the following, we will give an upper bound for the right-hand side of \eqref{EQ:vn-minus-eta_vn-1}.

To simplify the following terms, we will introduce the abbreviations
\begin{align*}
   \mu_n = \norm{\hat{\bsm}_h^n-\hat{\bsm}_h^{n-1}}{\infty,\Omega},\quad
   \nu_n = \min\Big\{\frac{\norm{\bsv_h^n-\bsv_h^{n-1}}{L^2(\Omega)}}%
      {\norm{\bsv_h^n}{L^2(\Omega)}},1\Big\},\quad
   {\Lambda_n
      = \frac{\norm{\grad\bsp_h^n}{L^2(\Omega)}}{\norm{\bsp_h^n}{L^2(\Omega)}}\,}.
\end{align*}
As a further preparation, we provide bounds for $\bsp_h^n$.
By~\cite[p.~1022~and~(8.5)]{AFKL:2021}, we obtain $\hat \bsm_h^n \in W^{1, \infty}(\Omega)$. Thus, using estimates from~\cite[Lem.~5.2]{AFKL:2021} yields
\begin{align*}%\label{EQ:}
   \norm{\bsp_h^n}{L^2(\Omega)}
   &\le c\norm{\hat{\bsm}_h^n-\hat{\bsm}_h^{n-1}}{\infty,\Omega}
      \norm{\bsv_h^n}{L^2(\Omega)}
    \le c\mu_n\norm{\bsv_h^n}{L^2(\Omega)}
\end{align*}
for some constant $c>0$ and
\begin{align*}
    \norm{\grad\bsp_h^n}{L^2(\Omega)}
       \leq c\mu_n\Lambda_n\norm{\bsv_h^n}{L^2(\Omega)}.
\end{align*}
Then we get the following estimates for the three terms on the right-hand side in \eqref{EQ:vn-minus-eta_vn-1}:
For the first term
\begin{align*}
   \eta\alpha(\bsv_h^{n-1},\bsp_h^n)_\Omega
   &\le c\mu_n\eta\alpha
      \norm{\bsv_h^{n-1}}{L^2(\Omega)}\norm{\bsv_h^n}{L^2(\Omega)}
    \le \frac{c}{2} \mu_n\eta\alpha\,
      \big(\norm{\bsv_h^{n-1}}{L^2(\Omega)}^2
      + \norm{\bsv_h^n}{L^2(\Omega)}^2\big)\,.
\end{align*}
For the second term, we can choose
\begin{align*}%\label{EQ:}
   &\eta(\hat{\bsm}_h^{n-1}\times
      \bsv_h^{n-1},\bsv_h^n-\bsv_h^{n-1}-\bsp_h^n)_\Omega\\
   &\qq\qq\le \eta\norm{\bsv_h^{n-1}}{L^2(\Omega)}
      \big(\norm{\bsv_h^n-\bsv_h^{n-1}}{L^2(\Omega)}
      + c\mu_n\norm{\bsv_h^n}{L^2(\Omega)}\big)\\
   &\qq\qq\le \frac12\Big(\frac{\norm{\bsv_h^n-\bsv_h^{n-1}}{L^2(\Omega)}}%
       {\norm{\bsv_h^n}{L^2(\Omega)}}+c\mu_n\Big)\eta\,
       \big(\norm{\bsv_h^{n-1}}{L^2(\Omega)}^2
       +\norm{\bsv_h^n}{L^2(\Omega)}^2\big)
\end{align*}
for $\sigma =1$, or
\begin{align*}
   \eta(\hat{\bsm}_h^{n-1}\times
      \bsv_h^{n-1},\bsv_h^n -\bsp_h^n)_\Omega
      \le \frac12(1+c\mu_n)\eta\,
      \big(\norm{\bsv_h^{n-1}}{L^2(\Omega)}^2
      + \norm{\bsv_h^n}{L^2(\Omega)}^2\big)
\end{align*}
for $\sigma =0$.
Thus, we can take
\begin{align*}%\label{EQ:}
   \eta(\hat{\bsm}_h^{n-1}\times
      \bsv_h^{n-1},\bsv_h^n-\bsv_h^{n-1}-\bsp_h^n)_\Omega
      \le \frac12(\nu_n+c\mu_n)\eta\,
       \big(\norm{\bsv_h^{n-1}}{L^2(\Omega)}^2
       +\norm{\bsv_h^n}{L^2(\Omega)}^2\big)\,.
\end{align*}
For the third term
\begin{align*}%\label{EQ:}
 \eta(\grad\bsm_h^{n-1},\grad\bsp_h^n)_\Omega
    &\le \eta\norm{\grad\bsm_h^{n-1}}{L^2(\Omega)}
      \norm{\grad\bsp_h^n}{L^2(\Omega)}
      \le c\eta\mu_n\Lambda_n
      \norm{\grad\bsm_h^{n-1}}{L^2(\Omega)}\norm{\bsv_h^n}{L^2(\Omega)}.
\end{align*}
If we now also take $\bsH_{\mathrm{ext}}^n \not=\bs0$ into account, the right-hand side will have the additional term 
\begin{align*}
   (1+c\mu_n)H_{\mathrm{ext}}^{n,n-1}\norm{\bsv_h^n}{L^2(\Omega)}
\end{align*}
with $H_{\mathrm{ext}}^{n,n-1} :=\norm{\bsH_{\mathrm{ext}}^n}{L^2(\Omega)}
+\norm{\bsH_{\mathrm{ext}}^{n-1}}{L^2(\Omega)}$.
We collect the formulas and get
\begin{align}\label{EQ:presum}
\begin{split}
   &\Big(\big(1-(\frac12+\frac{c}{2}\mu_n)\eta\big)\alpha
      - \frac12(\nu_n+c\mu_n)\eta\Big)
   \norm{\bsv_h^n}{L^2(\Omega)}^2\\
   &\qq{}-\Big((\frac12+\frac{c}{2}\mu_n)\eta\alpha
      +\frac12(\nu_n+c\mu_n)\eta\Big)\,
   \norm{\bsv_h^{n-1}}{L^2(\Omega)}^2\\
   &\qq\qq{}+ \frac{C_{\mathrm{e}}}{\tau_n}
      \big(\norm{\grad\bsM_h^n}{\bsG_n,\Omega}^2
      - \norm{\grad\bsM_h^{n-1}}{\bsG_n,\Omega}^2\big)\\
   &\q\le C_{\mathrm{e}} c\eta\mu_n\Lambda_n
      \norm{\grad\bsm_h^{n-1}}{L^2(\Omega)}\norm{\bsv_h^n}{L^2(\Omega)}
      + (1+c\mu_n)H_{\mathrm{ext}}^{n,n-1}\norm{\bsv_h^n}{L^2(\Omega)}\,.
      \end{split}
\end{align}
By \cite[p.~1028]{AFKL:2021} we have $\mu_n = O(h)$ using $\tau_n = O(h)$.
Clearly, $\nu_n = O(1)$ at least.
Finally, $\mu_n\Lambda_n = O(1)$ since by inverse estimate $\Lambda_n =O(1/h)$.
We fix $h_0$ so that $\mu_n\Lambda_n$ has a definite (but not necessarily small) bound for $h\le h_0$.

For notational simplicity, we let $\eta_{n,\alpha} =(\frac12+\frac{c}{2}\mu_{n})\eta\alpha + \frac12(\nu_{n}+c\mu_{n})\eta$.
Multiplying by $\tau_n$ and summing up the estimate~\eqref{EQ:presum} over $n$ yields
\begin{align*}
    &(\alpha - \eta_{n,\alpha})\tau_n \norm{\bsv_h^n}{L^2(\Omega)}^2
      + \sum_{j=k+2}^{n} \big( (\alpha-\eta_{j-1,\alpha}) \tau_{j-1}
      - \eta_{j,\alpha} \tau_{j}\big) \norm{\bsv_h^{j-1}}{L^2(\Omega)}^2
      + C_{\mathrm{e}} \norm{\grad\bsM_h^n}{\bsG_n,\Omega}^2\\
    &\q\le
      \sum_{j = k+1}^{n} \Big\{ C_{\mathrm{e}} c
         \eta\mu_{j}\Lambda_j
         \norm{\grad\bsm_h^{j-1}}{L^2(\Omega)}
         + (1+c\mu_{j})H_{\mathrm{ext}}^{j,j-1}\Big\}
         \tau_j\norm{\bsv_h^{j}}{L^2(\Omega)}
         \\[4pt]
    &\qq{}
      + C_{\mathrm{e}} \norm{\grad\bsM_h^{k}}{\bsG_{k+1},\Omega}^2
      + \eta_{k+1,\alpha} \tau_{k+1} \norm{\bsv_h^{k}}{L^2(\Omega)}^2
         + C_{\mathrm{e}}\sum_{j=k+1}^{n-1}\Big(\norm{\grad\bsM_h^j}
         {\bsG_{j+1},\Omega}^2-\norm{\grad\bsM_h^j}{\bsG_j,\Omega}^2\Big).
\end{align*}
From our definition of $\kappa_0$ above, we have $1-\frac12 g(\kappa)(1+\kappa)>2\omega$ for $\kappa \le\bar \kappa < \kappa_0$ and some $\omega>0$.
Then we choose $h$ so small and $\alpha$ so large that
\begin{align*}
   (\alpha-\eta_{j-1,\alpha}) \tau_{j-1}
      - \eta_{j,\alpha} \tau_{j}
   \ge \big( \alpha-\eta_{j-1,\alpha}
      - \eta_{j,\alpha} \kappa_{0}\big)\tau_{j-1}
   > \omega\tau_{j-1}
\end{align*}
for all $j \ge k+2$.
Inserting this lower bound, as well as $\eta\le 1$, yields
\begin{align*}
    &\omega\sum_{j=k+1}^{n} \tau_{j} \norm{\bsv_h^{j}}{L^2(\Omega)}^2
    + C_{\mathrm{e}} \norm{\grad\bsM_h^n}{\bsG_n,\Omega}^2\\
    &\quad
    \le \sum_{j = k+1}^{n} \Big\{{ C_{\mathrm{e}}}
       c\mu_{j}\Lambda_j
       \norm{\grad\bsm_h^{j-1}}{L^2(\Omega)}
       + (1+c\mu_{j})
         H_{\mathrm{ext}}^{j,j-1}\Big\}\tau_j\norm{\bsv_h^{j}}{L^2(\Omega)}
         \\[4pt]
    &\qq\qq{}
      + C_{\mathrm{e}} \norm{\grad\bsM_h^{k}}{\bsG_{k+1},\Omega}^2
      + \eta_{k+1,\alpha} \tau_{k+1} \norm{\bsv_h^{k}}{L^2(\Omega)}^2
         + C_{\mathrm{e}}\sum_{j=k+1}^{n-1}\Big(\norm{\grad\bsM_h^j}
         {\bsG_{j+1},\Omega}^2-\norm{\grad\bsM_h^j}{\bsG_j,\Omega}^2\Big).
\end{align*}
On the right-hand side we proceed with Young's inequality for the first sum
\begin{align*}
   &\sum_{j = k+1}^{n} \Big\{{C_{\mathrm{e}}}
       c\mu_{j}\Lambda_j
       \norm{\grad\bsm_h^{j-1}}{L^2(\Omega)}
       + (1+c\mu_{j})H_{\mathrm{ext}}^{j,j-1}\Big\}
         \tau_j \norm{\bsv_h^{j}}{L^2(\Omega)}\\
   &\qq\le
      \frac{\omega}2\sum_{j=k+1}^{n}\tau_j\norm{\bsv_h^{j}}{L^2(\Omega)}^2
      + \frac1{\omega}
      \sum_{j=k+1}^{n}\tau_j \big(C_{\mathrm{e}}c\mu_{j}\Lambda_j\big)^2
      \norm{\grad\bsm_h^{j-1}}{L^2(\Omega)}^2\\
   &\qq\qq{}+ \frac2{\omega}
      \sum_{j=k+1}^{n}\tau_j(1+c\mu_{j})^2
      \big(H_{\mathrm{ext}}^{j,j-1}\big)^2
\end{align*}
and apply Corollary~\ref{COR:commutator} and the step-ratio smoothness assumption \eqref{EQ:assumpt-smooth-kappa} for the second sum to obtain
\begin{align*}
    \sum_{j=k+1}^{n-1}\Big(\norm{\grad\bsM_h^j}{\bsG_{j+1},\Omega}^2
       &- \norm{\grad\bsM_h^j}{\bsG_j,\Omega}^2\Big)
     \le L c_\kappa \sum_{j=k+1}^{n-1} \tau_j
       \Big(\norm{\grad\bsm_h^{j-1}}{L^2(\Omega)}^2
       + \norm{\grad\bsm_h^j}{L^2(\Omega)}^2\Big)\\
    &= L c_\kappa \sum_{j=k+1}^{n-1} \tau_j
       \norm{\grad\bsm_h^{j-1}}{L^2(\Omega)}^2
       + L c_\kappa \sum_{j=k+2}^{n} \tau_{j-1}
       \norm{\grad\bsm_h^{j-1}}{L^2(\Omega)}^2\\
    &\le L c_\kappa \sum_{j=k+1}^{n} \tau_j \big(1 + \frac{1}{\kappa_j}\big)
       \norm{\grad\bsm_h^{j-1}}{L^2(\Omega)}^2,
\end{align*}
where $L$ is the constant of Proposition~\ref{PROP:Lipschitz} for $K=[\ukappa, \bar \kappa]$.
Absorbing the $\bsv$-sum of the right-hand side on the left-hand side, we obtain
\begin{align*}
    &\frac{\omega}2\sum_{j=k+1}^{n}
     \tau_{j} \norm{\bsv_h^{j}}{L^2(\Omega)}^2
     + C_{\mathrm{e}} \norm{\grad\bsM_h^n}{\bsG_n,\Omega}^2\\
    &\qq\le C_{\mathrm{e}}\norm{\grad\bsM_h^k}{\bsG_{k+1},\Omega}^2
      + \eta_{k+1,\alpha} \tau_{k+1} \norm{\bsv_h^{k}}{L^2(\Omega)}^2
      + \sum_{j=k+1}^{n}\tau_j \beta_j
        \norm{\grad\bsm_h^{j-1}}{L^2(\Omega)}^2\\
    &\qq\qq{}+ \frac2{\omega}
      \sum_{j=k+1}^{n}\tau_j(1+c\mu_{j})^2
      \big(H_{\mathrm{ext}}^{j,j-1}\big)^2,
\end{align*}
where
\begin{align*}
    \beta_j := \frac1{\omega} \big({C_{\mathrm{e}}}c\mu_{j}\Lambda_j\big)^2
       + C_{\mathrm{e}} L c_\kappa \big(1 + \frac{1}{\kappa_j}\big).
\end{align*}
Using the definition of $\bsM_h^n$ yields
\begin{align*}
    &\frac{\omega}2\sum_{j=k+1}^{n} \tau_{j} \norm{\bsv_h^{j}}{L^2(\Omega)}^2
       + C_{\mathrm{e}}\gamma_-^{\min} \norm{\grad\bsm_h^n}{L^2(\Omega)}^2\\
    &\quad\le 
       C_{\mathrm{e}}\gamma_+^{\max} \sum\limits_{i=1}^{k}
       \norm{\grad\bsm_h^{i}}{L^2(\Omega)}^2
       + \eta_{k+1,\alpha} \tau_{k+1} \norm{\bsv_h^{k}}{L^2(\Omega)}^2
       + \sum_{j=k+1}^{n}\tau_j \beta_j
       \norm{\grad\bsm_h^{j-1}}{L^2(\Omega)}^2\\
    &\qq\qq{}+ \frac2{\omega}
      \sum_{j=k+1}^{n}\tau_j(1+c\mu_{j})^2
      \big(H_{\mathrm{ext}}^{j,j-1}\big)^2\,.
\end{align*}
We now treat the problem in the form
\begin{align*}
    C_{\mathrm{e}}\gamma_-^{\min} \norm{\grad\bsm_h^n}{L^2(\Omega)}^2
    \le A + 
      \sum_{j=k+1}^{n}\tau_j \beta_j
      \norm{\grad\bsm_h^{j-1}}{L^2(\Omega)}^2\,.
\end{align*}
If we let $\xi = \max_j\{\big(\beta_j/( C_{\mathrm{e}} \gamma_-^{\min})\}$, we obtain by a discrete Gronwall's inequality
\begin{align*}
    C_{\mathrm{e}}\gamma_-^{\min} \norm{\grad\bsm_h^n}{L^2(\Omega)}^2
       \le A \mathrm{e}^{\xi T}\,.
\end{align*}
This provides us with a stability bound.
We can further derive
\begin{align*}
    \frac{\omega}2\sum_{j=k+1}^{n} \tau_{j} \norm{\bsv_h^{j}}{L^2(\Omega)}^2
    &\le A + \sum_{j=k+1}^{n}\tau_j \beta_j
      \norm{\grad\bsm_h^{j-1}}{L^2(\Omega)}^2\\
    &\le A \big(1 + \xi T\mathrm{e}^{\xi T}\big)
     \le A\big(1 + \xi T\big)\mathrm{e}^{\xi T}\,.
\end{align*}
Lastly, testing the weak equation for $n=k$ with $\bsphi_h = \bsv_h^k$ and \eqref{EQ:BDFvh} together with Young’s inequalities, we obtain
\begin{align*}
    \tau_{k+1} \frac{1}{2\kappa_{k+1}} \alpha \norm{\bsv_h^k}{L^2(\Omega)}^2
     + C_{\mathrm{e}} \frac{\delta_0}2 \norm{\grad\bsm_h^k}{L^2(\Omega)}^2
    &\le 
     \frac{C_{\mathrm{e}}}{2\delta_0}
     \Big(\sum\limits_{j=0}^{k-1} |\delta_{k-j}|
     \norm{\grad \bsm_h^{j}}{L^2(\Omega)} \Big)^2
     + \frac{\tau_k}{2 \alpha} \norm{\bsH_{\mathrm{ext}}^k}{L^2(\Omega)}^2\\
    &\le
     C \sum\limits_{j=0}^{k-1} \norm{\grad \bsm_h^{j}}{L^2(\Omega)}^2
     + \frac{\tau_k}{2 \alpha} \norm{\bsH_{\mathrm{ext}}^k}{L^2(\Omega)}^2\,,
\end{align*}
which concludes the proof.
\end{proof}

\begin{Remark}\ %
\begin{enumerate}
 \item
 A smooth change in the $\kappa_n$ is expected for smooth solutions \cite[p.~1042]{GearTu:1974}.
 \item
 The quantities $\mu_n$, $\nu_n$, $\Lambda_n$, whose a priori control needed the regularity requirements \eqref{EQ:regularity_m}, are in fact a posteriori terms.
 \item
 The lower bound on the damping parameter $\alpha$ is to be expected, since according to \cite[p.~1004]{AFKL:2021} any stability bound for A($\alpha$)-stable methods with $\alpha < \pi$ breaks down as eigenvalues of the linearisation approach the imaginary axis, see also the discussion in \cite[Sect.~2.5]{HEGP:2015}.
 However, we did not observe any difficulties to take values for $\alpha$ about $0.1$.
 \item
 The proof of Theorem~\ref{THM:stability} carries over to $k\ge3$, whenever analogues of Lemma~\ref{LEM:Gconstruction}, Proposition~\ref{PROP:Lipschitz} and Corollary~\ref{COR:commutator} are available for the BDF($k$) coefficients. 
 In particular, this is the case for BDF($k$), $k \in \{3, 4, 5\}$, with uniform time-steps and thus a step $n$-independent matrix $\bsG \in \IR^{k \times k}$.
 For variable time-steps and $k \ge 3$, an explicit construction as in the appendix analog to \eqref{EQ:Gidentity} would be required.
\end{enumerate}
\end{Remark}
 
%%%%%%%%%%%%%%%%%%%%%%%%%%%%%%%%%%%%%%%%%%%%%%%%%%%%%%%%%%%%%%%%%%%%%%%
\section{Adaptivity in Space}\label{SEC:adaptivitySpace}

%%%%%%%%%%%%%%%%%%%%%%%%%%%%%%%%%%%
\subsection{The recovery error indicator}\label{SSEC:recovery}

Since $\norm{\grad\bsm(t)}{L^2(\Omega)}$ is a quantity that appears in the energy bound \eqref{EQ:stability}, it is natural to also measure the error in this norm.
A fairly general technique for this is gradient recovery, initially introduced in~\cite{zzest}, see also~\cite{Ainsworth:2000} for an overview.
For the finite element space $W_h$, this defines a mapping $\grad_hv_h$ that projects the gradient of $v_h \in W_h$ into $W_h^3$.
Precise conditions can be obtained from~\cite{Ainsworth:2000}.
Having defined such $\grad_h$ (componentwise), we then define a local error indicator by
\begin{align}\label{EQ:local-error}
   \rho_{K,n} = \norm{\grad\bsm_h^n-\grad_h\bsm_h^n}{L^2(K)}\,
\end{align}
%\commmf{Replace $\bsG_h$ by $\nabla_h$ and $\eta$ by $\rho$ to avoid double meaning}
While the sum of the local error indicators provide only a rigorous lower bound for the true error, it is usually a very tight estimate of the local error and can even be asymptotically exact in case of superconvergence~\cite{Ainsworth:2000}, \cite{guo2025recovery,bartels2002each}.
In fact, the estimate can be compared with the best approximation in $H^1$ if the solution is smooth and the mesh is sufficiently fine \cite{BankYserentant:2015}.

Here we define $\grad_h$ as the $L^2(\Omega)$-projection from $\grad W_h$ to $W_h^3$, i.e.
\begin{align*}%label{EQ:L2-recovery}
   (\grad_hv_h,\bsphi_h)_{\Omega} = (\grad v_h,\bsphi_h)_{\Omega}
      \qq\text{for all }\bsphi_h\in W_h^3\,.
\end{align*}
For vector fields, this is applied component-wise.

%%%%%%%%%%%%%%%%%%%%%%%%%%%%%%%%%%%
\subsection{Refinement and coarsening}\label{SSEC:refinement}

Having calculated $[\rho_K]_{K\in\Tri_h}$ as in \eqref{EQ:local-error}, we extract a minimal subset $\cA_h \subset\Tri_h$ such that 
\begin{align}\label{EQ:refinement-condition}
   \sum_{K\in\cA_h}\rho_K^2
   \ge(1-\theta_{\mathrm{r}})\sum_{K\in\Tri_h}\rho_K^2
\end{align}
for some $\theta_{\mathrm{r}} \in(0,1)$.
A new mesh will be established by two bisections of elements in $\cA_h$.
To choose a set of elements to be coarsened, we extract a maximal subset $\cB_h \subset\Tri_h$ such that 
\begin{align}\label{EQ:coarsening-condition}
   \sum_{K\in\cB_h}\rho_K^2
   \le(1-\theta_{\mathrm{c}})\sum_{K\in\Tri_h}\rho_K^2\,.
\end{align}
Specific values for $\theta_{\mathrm{c}}$ and $\theta_{\mathrm{r}}$ will be provided for the respective examples below.
For local refinement and coarsening techniques for triangulations, see, for example, \cite{Baensch:1991}.

%%%%%%%%%%%%%%%%%%%%%%%%%%%%%%%%%%%%%%%%%%%%%%%%%%%%%%%%%%%%%%%%%%%%%%%
\section{Adaptivity in Space and Time}\label{SEC:adaptivity}

Our algorithm works in the following way: In each time step, it first fixes a suitable time step size.
Unlike \cite{AFKL:2021}, we use the extrapolation in time of order $k$ (instead of $k-1$) as a predictor.
This defines the spatial problem \eqref{EQ:LLG-weak-lambda-bdf} that is solved with adaptive mesh refinement.

\begin{Algorithm} \label{ALG:Algo1}\ \smallskip

\underline{Input}
\begin{itemize}\tightlist
 \item[] A coarse initial mesh $\Tri_h^{\text{ini}}$, final time $T>0$ and tolerances $\tol_{\mathrm{s}}$ (in space) and $\tol_{\mathrm{t}}$ (in time), refinement/coarsening controls $\theta_{\mathrm{r}},\theta_{\mathrm{c}} >0$ and the spatial error indicator $\rho$ \eqref{EQ:local-error}.
\end{itemize}

\underline{Precomputation}
\begin{itemize}\tightlist
 \item We iteratively refine $\Tri_h^{\text{ini}}$ to $\Tri_h^0$ using mesh refinement on elements in $\cA_h$ chosen by \eqref{EQ:refinement-condition} with respect to the given initial condition $\bsm(0)$ and threshold $\theta_{\mathrm{r}}$ until $\rho \le\tol_{\mathrm{s}}$.
 
 \item Compute $\bsv_h^0$ and solve \eqref{EQ:LLG-weak-lambda-bdf}  once with the implicit Euler method, i.e.\ BDF(1), on the mesh $\Tri_h^0$ with the time step $\tilde{\tau}_1 = \tol_t/\norm{\bsv_h^0}{L^2(\Omega)}$ to obtain $\tilde{\bsv}_h^1$.
 We then define $\tau_1$ by \eqref{EQ:def_tau-1}, where we employ $\norm{\ptl_t^2\bsm}{L^2(\Omega)} \approx\norm{\tilde{\bsv}_h^1-\bsv_h^0}{L^2(\Omega)}/\tilde{\tau}_1$.
 
 \item Apply $k-1$ time steps to solve \eqref{EQ:LLG-weak-lambda-bdf} with the second order singly-diagonal implicit Runge--Kutta (SDIRK) method as proposed in \cite{Nishikawa:2019} on the mesh $\Tri_h^0$.
 Thus, we end up with the states $\bsm_h^j,\bsv_h^j$ for $j =0,\dots,k-1$.
\end{itemize}

\goodbreak
\underline{Time-stepping}
\begin{itemize}\tightlist
 \item For $n \ge k$: Interpolate the approximations $\bsm_h^{n-2},\dots,\bsm_h^{n-k}$ at times $t_{n-2},\dots,t_{n-k}$ on the grid $\Tri_h^{n-1}$ (the approximation $\bsm_h^{n-1}$ at time $t_{n-1}$ on grid $\Tri_h^{n-1}$ was already calculated). 

 \item Continue by normalizing the magnetization $\tilde \bsm_h^{n-1} = \bsm_h^{n-1}/|\bsm_h^{n-1}|$ prior to evaluating the spatial error estimator. This improves the effectiveness of coarsening. Note that the normalized magnetization is used only for the error computation  and is not propagated to the next step. Then use the spatial indicator $\rho$ from \eqref{EQ:local-error} for $\tilde \bsm_h^{n-1}$ to first coarsen the mesh using the threshold $\theta_{\mathrm{c}}$. Then refine the mesh using the threshold $\theta_{\mathrm{r}}$ until the tolerance $\tol_{\mathrm{s}}$ is achieved. 
 Both thresholds are described in Section \ref{SSEC:refinement}.
 This results in the new mesh $\Tri_h^{n}$.
 
 \item Compute the next time step size $\tau_{n}$ for BDF($k$) as described in Section~\ref{SSEC:stepselection} with $\tol_{\mathrm{t}}$, \eqref{EQ:def_tau-2} or \eqref{EQ:def_tau-3}, and $\norm{\ptl_t^{k+1}\bsm(t_{n-1})}{L^2(\Omega)} \approx\norm{\ptl_\tau^{k+1}\bsm_h(t_{n-1})}{L^2(\Omega)} \approx\norm{\ptl_\tau^{k}\bsv_h(t_{n-1})}{L^2(\Omega)}$ using \eqref{EQ:def_dh2y} or \eqref{EQ:def_dh3y}, respectively.

 \item Extrapolate $\hat{\bsm}_h^{n}$ and solve \eqref{EQ:LLG-weak-lambda-bdf} for $\bsm_h^{n}$ at time $t_{n} =t_{n-1}+\tau_n$, increase $n$ as long as $t_{n} <T$ and close with a final step to reach $T$.
\end{itemize}
\end{Algorithm}

%%%%%%%%%%%%%%%%%%%%%%%%%%%%%%%%%%%%%%%%%%%%%%%%%%%%%%%%%%%%%%%%%%%%%%%
\section{Numerical Experiments}\label{SEC:numerical-experiments}

In the following, we apply our algorithm to a number of examples in order to study its performance.
In order to measure the error, we choose (for $T=t_N$)
\begin{align}\label{EQ:error-norm}
   \textrm{err}_T = \max\limits_{0\le n\le N}
      \norm{\grad(\bsm(t_n)-\bsm_h^n)}{L^2(\Omega)}
\end{align}
which corresponds to the theoretical result in Section~\ref{SSEC:uniform-results} and to the chosen error indicator \eqref{EQ:local-error}.
The algorithm was implemented with the finite element library \texttt{deal.II} \cite{BangerthET:2023} implemented in \texttt{C++}.
In this work, valuable features of the \texttt{deal.II}  software include the ability to perform $h$-refinement, which allows us to refine and coarsen the mesh based on a fixed fraction of the estimated error based on \eqref{EQ:refinement-condition} and \eqref{EQ:coarsening-condition}.
In general, we only coarsen the spatial mesh every tenth time step.
We solved in each time step the saddle point problem \eqref{EQ:LLG-weak-lambda-bdf}. The resulting linear systems were solved using the UMFPACK package \cite{Davis:2004} for the LU-decomposition.

%%%%%%%%%%%%%%%%%%%%%%%%%%%%%%%%%%
\subsection{Example 1}\label{SSEC:example1}
The following example is taken from \cite[Sect.~9.2, (9.2)]{AFKL:2021}.
In this case, we let $\Omega =(0,1)\times(0,1)$ and provide a source term $\bsH_{\mathrm{ext}}$ such that the exact solution of \eqref{EQ:LLG-general-pde}--\eqref{EQ:LLG-general-iv} is
\begin{align*}
  \bsm(t,\bsx) =
    \Bmat{-(x_1^3-3/2\,x_1^2+1/4)\sin(3\pi t/T)\\
          \sqrt{1-(x_1^3-3/2\,x_1^2+1/4)^2}\\
          -(x_1^3-3/2\,x_1^2+1/4)\cos(3\pi t/T)}.
\end{align*}
We choose $T =0.1$, $C_{\text{e}} =1$, and $\alpha = 0.2$. 
In this example, the temporal error dominates the total error, so we fix $p =3$ for the polynomial degree in space.

\begin{figure}[ht]
    \center 
    \resizebox{0.48\textwidth}{6.2cm}{%
        % This file was created with tikzplotlib v0.10.1.
\begin{tikzpicture}

\definecolor{darkgray176}{RGB}{176,176,176}
\definecolor{limegreen}{RGB}{50,205,50}
\definecolor{purple}{RGB}{128,0,128}
\definecolor{royalblue}{RGB}{65,105,225}

\begin{axis}[
legend cell align={left},
legend style={
  fill opacity=0.8,
  draw opacity=1,
  text opacity=1,
  at={(0.97,0.03)},
  anchor=south east,
  draw=none
},
log basis x={10},
log basis y={10},
tick pos=both,
% x grid style={darkgray176},
xlabel={$\tol_{\mathrm{t}}$},
xmajorgrids,
% xmin=3.22158474947005e-08, xmax=0.000303142731278652,
xmode=log,
% xtick style={color=black},
% xtick={1e-09,1e-08,1e-07,1e-06,1e-05,0.0001,0.001,0.01},
% xticklabels={
%   \(\displaystyle {10^{-9}}\),
%   \(\displaystyle {10^{-8}}\),
%   \(\displaystyle {10^{-7}}\),
%   \(\displaystyle {10^{-6}}\),
%   \(\displaystyle {10^{-5}}\),
%   \(\displaystyle {10^{-4}}\),
%   \(\displaystyle {10^{-3}}\),
%   \(\displaystyle {10^{-2}}\)
% },
% y grid style={darkgray176},
ylabel={\(\displaystyle \textrm{err}_T \)},
ymajorgrids,
% ymin=5.64629319658561e-06, ymax=0.0640195499590323,
ymode=log,
% ytick style={color=black},
% ytick={1e-07,1e-06,1e-05,0.0001,0.001,0.01,0.1,1},
% yticklabels={
%   \(\displaystyle {10^{-7}}\),
%   \(\displaystyle {10^{-6}}\),
%   \(\displaystyle {10^{-5}}\),
%   \(\displaystyle {10^{-4}}\),
%   \(\displaystyle {10^{-3}}\),
%   \(\displaystyle {10^{-2}}\),
%   \(\displaystyle {10^{-1}}\),
%   \(\displaystyle {10^{0}}\)
% }
]
\addplot [thick, red, dashed, mark=*, mark size=1.5, mark options={solid}]
table {%
0.0002 0.04133
0.0001 0.02933
5e-05 0.02078
2.5e-05 0.01473
1.25e-05 0.01043
6.25e-06 0.007383
3.125e-06 0.005226
1.563e-06 0.003697
7.813e-07 0.002615
3.906e-07 0.00185
1.953e-07 0.001308
9.766e-08 0.0009252
4.883e-08 0.0006544
};
\addlegendentry{$k=1$}
\addplot [thick, red, dotted, forget plot]
table {%
4.883e-08 0.00052352
0.0002 0.0335046367171505
};
% \addlegendentry{eoc=0.5}
\addplot [thick, royalblue, dashed, mark=*, mark size=1.5, mark options={solid}]
table {%
0.0002 0.009641
0.0001 0.006035
5e-05 0.003781
2.5e-05 0.002367
1.25e-05 0.001486
6.25e-06 0.000933
3.125e-06 0.0005862
1.563e-06 0.0003686
7.813e-07 0.000232
3.906e-07 0.0001462
1.953e-07 9.239e-05
9.766e-08 5.876e-05
4.883e-08 3.794e-05
};
\addlegendentry{$k=2$}
\addplot [thick, royalblue, dotted, forget plot]
table {%
4.883e-08 3.0352e-05
0.0002 0.00735078757637011
};
% \addlegendentry{eoc=0.66}
\addplot [thick, limegreen, dashed, mark=*, mark size=1.5, mark options={solid}]
table {%
0.0002 0.01635
0.0001 0.01006
5e-05 0.00618
2.5e-05 0.003745
1.25e-05 0.002265
6.25e-06 0.001363
3.125e-06 0.0008188
1.563e-06 0.0004908
7.813e-07 0.0002937
3.906e-07 0.0001756
1.953e-07 0.0001052
9.766e-08 6.34e-05
4.883e-08 3.874e-05
};
\addlegendentry{$k=3$}
\addplot [thick, limegreen, dotted, forget plot]
table {%
4.883e-08 3.0992e-05
0.0002 0.0158674470197193
};
% \addlegendentry{eoc=0.75}
\addplot [thick, purple, dashed, mark=*, mark size=1.5, mark options={solid}]
table {%
0.0002 0.01086
0.0001 0.005982
5e-05 0.003162
2.5e-05 0.001632
1.25e-05 0.0008563
6.25e-06 0.0004523
3.125e-06 0.0002402
1.563e-06 0.0001248
7.813e-07 6.748e-05
3.906e-07 3.698e-05
1.953e-07 2.174e-05
9.766e-08 1.508e-05
4.883e-08 1.233e-05
};
\addlegendentry{$k=4$}
\addplot [thick, purple, dotted, forget plot]
table {%
4.883e-08 9.864e-06
0.0002 0.00765469129336043
};
% \addlegendentry{eoc=0.8}
\end{axis}

\end{tikzpicture}
    }
    \quad
    \resizebox{0.48\textwidth}{6.2cm}{%
        % This file was created with tikzplotlib v0.10.1.
\begin{tikzpicture}

\definecolor{darkgray176}{RGB}{176,176,176}
\definecolor{limegreen}{RGB}{50,205,50}
\definecolor{purple}{RGB}{128,0,128}
\definecolor{royalblue}{RGB}{65,105,225}

\begin{axis}[
legend cell align={left},
legend style={
  fill opacity=0.8,
  draw opacity=1,
  text opacity=1,
  at={(0.97,0.03)},
  anchor=south east,
  draw=none
},
log basis x={10},
log basis y={10},
tick pos=both,
% x grid style={darkgray176},
xlabel={$\tol_{\mathrm{t}}$},
xmajorgrids,
% xmin=3.22158474947005e-08, xmax=0.000303142731278652,
xmode=log,
% xtick style={color=black},
% xtick={1e-09,1e-08,1e-07,1e-06,1e-05,0.0001,0.001,0.01},
% xticklabels={
%   \(\displaystyle {10^{-9}}\),
%   \(\displaystyle {10^{-8}}\),
%   \(\displaystyle {10^{-7}}\),
%   \(\displaystyle {10^{-6}}\),
%   \(\displaystyle {10^{-5}}\),
%   \(\displaystyle {10^{-4}}\),
%   \(\displaystyle {10^{-3}}\),
%   \(\displaystyle {10^{-2}}\)
% },
% y grid style={darkgray176},
ylabel={\(\displaystyle \max_n \| |\bsm_h^n| - 1\|_{L^\infty(\Omega)}\)},
ymajorgrids,
% ymin=3.2993369902719e-07, ymax=0.0234260188237188,
ymode=log,
% ytick style={color=black},
% ytick={1e-08,1e-07,1e-06,1e-05,0.0001,0.001,0.01,0.1,1},
% yticklabels={
%   \(\displaystyle {10^{-8}}\),
%   \(\displaystyle {10^{-7}}\),
%   \(\displaystyle {10^{-6}}\),
%   \(\displaystyle {10^{-5}}\),
%   \(\displaystyle {10^{-4}}\),
%   \(\displaystyle {10^{-3}}\),
%   \(\displaystyle {10^{-2}}\),
%   \(\displaystyle {10^{-1}}\),
%   \(\displaystyle {10^{0}}\)
% }
]
\addplot [thick, red, dashed, mark=*, mark size=1.5, mark options={solid}]
table {%
0.0002 0.01393
0.0001 0.009884
5e-05 0.007004
2.5e-05 0.004964
1.25e-05 0.003515
6.25e-06 0.002488
3.125e-06 0.001761
1.563e-06 0.001246
7.813e-07 0.000881
3.906e-07 0.0006231
1.953e-07 0.0004406
9.766e-08 0.0003116
4.883e-08 0.0002203
};
\addlegendentry{$k=1$}
\addplot [thick, red, dotted, forget plot]
table {%
4.883e-08 0.00017624
0.0002 0.0112791434425248
};
% \addlegendentry{eoc=0.5}
\addplot [thick, royalblue, dashed, mark=*, mark size=1.5, mark options={solid}]
table {%
0.0002 0.001452
0.0001 0.0007371
5e-05 0.000372
2.5e-05 0.0001865
1.25e-05 9.372e-05
6.25e-06 4.699e-05
3.125e-06 2.348e-05
1.563e-06 1.187e-05
7.813e-07 6.031e-06
3.906e-07 3.108e-06
1.953e-07 1.646e-06
9.766e-08 9.146e-07
4.883e-08 5.482e-07
};
\addlegendentry{$k=2$}
\addplot [thick, royalblue, dotted, forget plot]
table {%
4.883e-08 2.741e-07
0.0002 0.0011226704894532
};
% \addlegendentry{eoc=1}
\addplot [thick, limegreen, dashed, mark=*, mark size=1.5, mark options={solid}]
table {%
0.0002 0.006643
0.0001 0.004093
5e-05 0.002512
2.5e-05 0.001523
1.25e-05 0.0009203
6.25e-06 0.0005536
3.125e-06 0.0003325
1.563e-06 0.0001992
7.813e-07 0.0001191
3.906e-07 7.117e-05
1.953e-07 4.252e-05
9.766e-08 2.541e-05
4.883e-08 1.517e-05
};
\addlegendentry{$k=3$}
\addplot [thick, limegreen, dotted, forget plot]
table {%
4.883e-08 1.0619e-05
0.0002 0.00543677142173462
};
% \addlegendentry{eoc=0.75}
\addplot [thick, purple, dashed, mark=*, mark size=1.5, mark options={solid}]
table {%
0.0002 0.004569
0.0001 0.002462
5e-05 0.001281
2.5e-05 0.0006512
1.25e-05 0.0003333
6.25e-06 0.0001704
3.125e-06 8.682e-05
1.563e-06 4.334e-05
7.813e-07 2.198e-05
3.906e-07 1.106e-05
1.953e-07 5.574e-06
9.766e-08 2.839e-06
4.883e-08 1.449e-06
};
\addlegendentry{$k=4$}
\addplot [thick, purple, dotted, forget plot]
table {%
4.883e-08 1.0143e-06
0.0002 0.00415441327053041
};
% \addlegendentry{eoc=1.0}
\end{axis}

\end{tikzpicture}
    }
    \caption{\label{FIG:Ex1-adapt-ts}
    Example~\ref{SSEC:example1}: Left: Convergence of $\text{err}_T$ with respect to $\tol_t$ for adaptive time-stepping for BDF(1) (red), BDF(2) (blue), BDF(3) (green) and BDF(4) (purple) with fixed $h = 1/16$ and polynomial degree $p=3$.
    The corresponding experimental rate of convergence is shown in the same color by the dotted line, here for $k = 1, \dots, 4$ the rate is $k / (k+1)$.
    Right: Convergence of $\max\limits_{0\le n\le N} \norm{|\bsm_h^n| - 1}{L^\infty(\Omega)}$ with respect to $\tol_t$ for adaptive time-stepping for BDF($k$), $k=1, \dots 4$, and the rate of $1/2$ for $k=1$, $3/4$ for $k=3$ and $1$ for $k=2, 4$ in dotted lines.
    }
\end{figure}

In Figure~\ref{FIG:Ex1-adapt-ts} (left) we show the error \eqref{EQ:error-norm} of the adaptive time-stepping method for BDF(1) to BDF(4) depending on the prescribed temporal tolerance $\tol_{\mathrm{t}}$ for a fixed mesh.
For BDF($k$) we have the convergence order $k$, but a local truncation error of order $k+1$, for uniform steps.
Since we equilibrate the local truncation error, we have $\tau^{k+1} \sim \tol_{\mathrm{t}}$ and therefore we expect $\text{err}_T \sim\tol_{\mathrm{t}}^{k/(k+1)}$ and this is indeed observed.
Since the discrete solution is not explicitly normalized, Figure~\ref{FIG:Ex1-adapt-ts} (right) displays the convergence rate of the maximal normalization error $\max\limits_{0\le n\le N} \norm{|\bsm_h^n| - 1}{L^\infty(\Omega)}$.
We obtain convergence of at least order $k / (k+1)$ for $k=1, \dots, 4$.

%%%%%%%%%%%%%%%%%%%%%%%%%%%%%%%%%%%
\subsection{Example 2}\label{SSEC:example2}

The following example is taken from \cite[Sect.~9.2, (9.1)]{AFKL:2021}.
In this case, we let $\Omega =(0,1)\times(0,1)$ and provide a source term $\bsH_{\mathrm{ext}}$ such that the exact solution of \eqref{EQ:LLG-general-pde}--\eqref{EQ:LLG-general-iv} is
\begin{align*}
  \bsm(t,\bsx) &=
  \begin{cases}
    \bmat{C_0(x_1-1/2)\,\ez^{-\frac{g(t)}{1/4-d(\bsx)}}\\
          C_0(x_2-1/2)\,\ez^{-\frac{g(t)}{1/4-d(\bsx)}}\\
          \sqrt{1-C_0^2d(\bsx)\,\ez^{-\frac{2g(t)}{1/4-d(\bsx)}}}}\quad
          & \mbox{if } d(\bsx) < 1/4\,,\vspace{4pt}\\
    \bmat{0\\ 0\\ 1} & \mbox{otherwise}\,,
  \end{cases}
\end{align*}
with $\bsm^0 =\bsm(0)$, $g(t) =(T_0+0.1)/(T_0+0.1-t)$, $d(\bsx) =|\bsx-[1,1]/2|^2$, $C_0 =400$.
We choose $C_{\mathrm{e}} =1$, $\alpha =0.2$, $T_0=0.06$ and $T=0.05$ as the final time.
In this example, the spatial error dominates the total error.
We use $\tau=T/100$, as well as $\theta_{\mathrm{r}}=0.85$ and $\theta_{\mathrm{c}} = 0.9$.

We perform numerical experiments for linear, quadratic and cubic finite elements on uniform and adaptive spatial meshes.
Time integration is done with uniform BDF(2).
We execute Algorithm~\ref{ALG:Algo1} for several mesh sizes $h$ and for the adaptation of space for several tolerances $\tol_{\mathrm{s}}$.
In Figure~\ref{FIG:Ex2-tol_s-vs-err_T} (left) we notice convergence in the error norm $\text{err}_T$ given by \eqref{EQ:error-norm} of order one in $\tol_{\mathrm{s}}$ for all polynomial degrees until saturation occurs due to the fixed time step.
However, in Figure~\ref{FIG:Ex2-tol_s-vs-err_T} (right), we notice smaller errors in the error norm $\text{err}_T$ for higher polynomial degrees if we relate these errors to the maximal degrees of freedom we obtained at a discrete time step during the simulation.

In Figure~\ref{FIG:Ex2-tol_s_vs_norm} we show the convergence of the maximal normalization error $\max\limits_{0\le n\le N} \norm{|\bsm_h^n| - 1}{L^\infty(\Omega)}$ with respect to the spatial tolerance $\tol_\textrm{s}$. 
For linear finite elements, we observe a convergence rate that is one order higher with respect to the spatial tolerance $\tol_\textrm{s}$ than for the higher order elements.

\begin{figure}[ht]
    \center 
    \resizebox{0.48\textwidth}{!}{%
        % This file was created with tikzplotlib v0.10.1.
\begin{tikzpicture}

\definecolor{darkgray176}{RGB}{176,176,176}
\definecolor{green}{RGB}{0,128,0}

\begin{axis}[
legend cell align={left},
legend style={
  fill opacity=0.8,
  draw opacity=1,
  text opacity=1,
  at={(0.97,0.03)},
  anchor=south east,
  draw=none
},
log basis x={10},
log basis y={10},
tick pos=both,
% x grid style={darkgray176},
xlabel={$\tol_\textrm{s}$},
xmajorgrids,
xmode=log,
% xtick style={color=black},
% xtick={1e-05,0.0001,0.001,0.01,0.1,1,10},
% xticklabels={
%   \(\displaystyle {10^{-5}}\),
%   \(\displaystyle {10^{-4}}\),
%   \(\displaystyle {10^{-3}}\),
%   \(\displaystyle {10^{-2}}\),
%   \(\displaystyle {10^{-1}}\),
%   \(\displaystyle {10^{0}}\),
%   \(\displaystyle {10^{1}}\)
% },
% y grid style={darkgray176},
ylabel={\(\displaystyle \textrm{err}_T \)},
ymajorgrids,
% ymin=4.20212097619745e-05, ymax=0.113945184518053,
ymode=log,
% ytick style={color=black},
% ytick={1e-06,1e-05,0.0001,0.001,0.01,0.1,1,10},
% yticklabels={
%   \(\displaystyle {10^{-6}}\),
%   \(\displaystyle {10^{-5}}\),
%   \(\displaystyle {10^{-4}}\),
%   \(\displaystyle {10^{-3}}\),
%   \(\displaystyle {10^{-2}}\),
%   \(\displaystyle {10^{-1}}\),
%   \(\displaystyle {10^{0}}\),
%   \(\displaystyle {10^{1}}\)
% }
]
\addplot [thick, red, dashed, mark=*, mark size=1.5, mark options={solid}]
table {%
0.2 0.02614
0.1333 0.01844
0.08889 0.01158
0.05926 0.00789
0.03951 0.005076
0.02634 0.003417
0.01756 0.002078
};
\addlegendentry{$p=1$}
\addplot [thick, blue, dashed, mark=*, mark size=1.5, mark options={solid}]
table {%
0.2 0.07955
0.1333 0.07955
0.08889 0.02276
0.05926 0.02276
0.03951 0.01464
0.02634 0.008418
0.01756 0.005365
0.01171 0.005027
0.007804 0.005184
0.005202 0.002467
0.003468 0.001775
0.002312 0.001033
0.001541 0.00102
0.001028 0.0006441
0.0006851 0.0003686
0.0004567 0.0002982
0.0003045 0.0002828
0.000203 0.0002135
};
\addlegendentry{$p=2$}
\addplot [thick, green, dashed, mark=*, mark size=1.5, mark options={solid}]
table {%
0.2 0.02089
0.1333 0.02089
0.08889 0.02089
0.05926 0.02089
0.03951 0.01157
0.02634 0.01157
0.01756 0.004571
0.01171 0.004571
0.007804 0.003547
0.005202 0.001581
0.003468 0.001531
0.002312 0.000633
0.001541 0.0006313
0.001028 0.000529
0.0006851 0.0002929
0.0004567 0.0001888
0.0003045 0.0001888
0.000203 0.0001749
};
\addlegendentry{$p=3$}
\addplot [thick, black, dotted]
table {%
0.000203 4.27e-05
0.2 0.0420689655172414
};
\addlegendentry{$\cO(\tol_\textrm{s})$}
\end{axis}

\end{tikzpicture}
    }
    \quad
    \resizebox{0.48\textwidth}{!}{%
        % This file was created with tikzplotlib v0.10.1.
\begin{tikzpicture}

\definecolor{darkgray176}{RGB}{176,176,176}
\definecolor{green}{RGB}{0,128,0}

\begin{axis}[
legend cell align={left},
legend style={
  fill opacity=0.8,
  draw opacity=1,
  text opacity=1,
  at={(0.03,0.03)},
  anchor=south west,
  draw=none
},
log basis x={10},
log basis y={10},
tick pos=both,
% x grid style={darkgray176},
xlabel={\(\displaystyle N_\textrm{s}^{\max}\)},
xmajorgrids,
xmode=log,
% xtick style={color=black},
% xtick={100,1000,10000,100000,1000000,10000000,100000000},
% xticklabels={
%   \(\displaystyle {10^{2}}\),
%   \(\displaystyle {10^{3}}\),
%   \(\displaystyle {10^{4}}\),
%   \(\displaystyle {10^{5}}\),
%   \(\displaystyle {10^{6}}\),
%   \(\displaystyle {10^{7}}\),
%   \(\displaystyle {10^{8}}\)
% },
% y grid style={darkgray176},
ylabel={\(\displaystyle \textrm{err}_T\)},
ymajorgrids,
ymode=log,
% ytick style={color=black},
% ytick={1e-05,0.0001,0.001,0.01,0.1,1,10},
% yticklabels={
%   \(\displaystyle {10^{-5}}\),
%   \(\displaystyle {10^{-4}}\),
%   \(\displaystyle {10^{-3}}\),
%   \(\displaystyle {10^{-2}}\),
%   \(\displaystyle {10^{-1}}\),
%   \(\displaystyle {10^{0}}\),
%   \(\displaystyle {10^{1}}\)
% }
]
\addplot [semithick, red, dashed, mark=*, mark size=1.5, mark options={solid}]
table {%
11080 0.1831
20640 0.1316
45960 0.08585
102200 0.05599
219700 0.0379
452100 0.02586
1131000 0.01646
};
\addlegendentry{$p=1$}
\addplot [semithick, blue, dashed, mark=*, mark size=1.5, mark options={solid}]
table {%
2052 0.1928
2052 0.1928
3652 0.1033
3652 0.1033
8964 0.0619
10370 0.04696
10880 0.03327
12480 0.02668
23300 0.02303
41090 0.01204
42240 0.009804
48900 0.007187
53440 0.006485
133700 0.003561
170600 0.00253
208100 0.002088
230100 0.0019
569300 0.001457
};
\addlegendentry{$p=2$}
\addplot [semithick, green, dashed, mark=*, mark size=1.5, mark options={solid}]
table {%
2500 0.07822
2500 0.07822
2500 0.07822
2500 0.07822
4452 0.05967
4452 0.05967
8004 0.01473
8004 0.01473
12130 0.01176
18920 0.006574
22820 0.005394
27650 0.002629
30980 0.002124
40740 0.001687
72520 0.00137
91750 0.001232
91750 0.001232
103300 0.001198
};
\addlegendentry{$p=3$}
\end{axis}

\end{tikzpicture}
    }
    \caption{\label{FIG:Ex2-tol_s-vs-err_T}
    Example~\ref{SSEC:example2}: Left: Spatial tolerance $\tol_{\mathrm{s}}$ vs.\ $\textrm{err}_T$. Right: Maximal spatial degrees of freedom $N_{\mathrm{s}}^\text{max} := \max\{N_{\mathrm{s}}(t_n)\fdg 0\le n\le 100\}$ vs.\ $\textrm{err}_T$, where $N_{\mathrm{s}}$ denotes the number of degrees of freedom.
    Various polynomial degrees of the finite element space are presented, i.e., linear (red), quadratic (blue) and cubic (green).}
\end{figure}

\begin{figure}[ht]
    \center 
    \resizebox{0.48\textwidth}{!}{%
        % This file was created with tikzplotlib v0.10.1.
\begin{tikzpicture}

\definecolor{darkgray176}{RGB}{176,176,176}
\definecolor{green}{RGB}{0,128,0}

\begin{axis}[
legend cell align={left},
legend style={
  fill opacity=0.8,
  draw opacity=1,
  text opacity=1,
  at={(0.03,0.97)},
  anchor=north west,
  draw=none
},
log basis x={10},
log basis y={10},
tick pos=both,
% x grid style={darkgray176},
xlabel={$\tol_\textrm{s}$},
xmajorgrids,
xmode=log,
% xtick style={color=black},
% xtick={1e-05,0.0001,0.001,0.01,0.1,1,10},
% xticklabels={
%   \(\displaystyle {10^{-5}}\),
%   \(\displaystyle {10^{-4}}\),
%   \(\displaystyle {10^{-3}}\),
%   \(\displaystyle {10^{-2}}\),
%   \(\displaystyle {10^{-1}}\),
%   \(\displaystyle {10^{0}}\),
%   \(\displaystyle {10^{1}}\)
% },
% y grid style={darkgray176},
ylabel={\(\displaystyle \max_n \| |\bsm_h^n| - 1\|_{L^\infty(\Omega)}\)},
ymajorgrids,
% ymin=4.76320565423018e-06, ymax=0.00934946341140485,
ymode=log,
% ytick style={color=black},
% ytick={1e-07,1e-06,1e-05,0.0001,0.001,0.01,0.1},
% yticklabels={
%   \(\displaystyle {10^{-7}}\),
%   \(\displaystyle {10^{-6}}\),
%   \(\displaystyle {10^{-5}}\),
%   \(\displaystyle {10^{-4}}\),
%   \(\displaystyle {10^{-3}}\),
%   \(\displaystyle {10^{-2}}\),
%   \(\displaystyle {10^{-1}}\)
% }
]
\addplot [red, dashed, mark=*, mark size=1.5, mark options={solid}]
table {%
0.2 0.003226
0.1333 0.001388
0.08889 0.0008108
0.05926 0.0002658
0.03951 0.000203
0.02634 5.793e-05
0.01756 2.684e-05
};
\addlegendentry{$p=1$}
\addplot [blue, dashed, mark=*, mark size=1.5, mark options={solid}]
table {%
0.2 0.006153
0.1333 0.006153
0.08889 0.001985
0.05926 0.001985
0.03951 0.0007071
0.02634 0.0005947
0.01756 0.0002384
0.01171 0.0002384
0.007804 0.0002295
0.005202 0.000105
0.003468 7.863e-05
0.002312 2.915e-05
0.001541 2.915e-05
0.001028 1.396e-05
0.0006851 8.428e-06
0.0004567 6.539e-06
0.0003045 5.094e-06
0.000203 3.883e-06
};
\addlegendentry{$p=2$}
\addplot [green, dashed, mark=*, mark size=1.5, mark options={solid}]
table {%
0.2 0.003631
0.1333 0.003631
0.08889 0.003631
0.05926 0.003631
0.03951 0.0006518
0.02634 0.0006518
0.01756 0.0002738
0.01171 0.0002738
0.007804 0.0002089
0.005202 0.0001106
0.003468 0.0001105
0.002312 1.982e-05
0.001541 1.982e-05
0.001028 1.393e-05
0.0006851 8.797e-06
0.0004567 6.362e-06
0.0003045 6.362e-06
0.000203 3.996e-06
};
\addlegendentry{$p=3$}
\addplot [thick, black, dotted]
table {%
0.000203 1.5532e-06
0.2 0.00153024630541872
};
\addlegendentry{$\cO(\tol_\textrm{s})$}
\addplot [thick, red, dotted]
table {%
0.01756 8.5888e-06
0.2 0.00111414946995916
};
\addlegendentry{$\cO(\tol_\textrm{s}^2)$}
\end{axis}

\end{tikzpicture}
    } 
    \caption{\label{FIG:Ex2-tol_s_vs_norm}
    Example~\ref{SSEC:example2}: Convergence of $\max\limits_{0\le n\le N} \norm{|\bsm_h^n| - 1}{L^\infty(\Omega)}$ with respect to the spatial tolerance $\tol_{\mathrm{s}}$.
    Various polynomial degrees of the finite element space are presented, i.e., linear (red), quadratic (blue) and cubic (green).}
\end{figure}

%%%%%%%%%%%%%%%%%%%%%%%%%%%%%%%%%%%
\subsection{Example 3}\label{SSEC:example3}\ %
Here we consider an example that appears to show singular behavior and has been considered in~\cite{BartelsProhl:2006,BKP:2008,BBP:2008,ChengShen:2023}.
However, we will demonstrate that computations on finer meshes and particularly adaptive computations produce a stable solution that does not appear to collapse into a singularity.
We let $\Omega =(-1/2,1/2)^2$ and choose constants $C_{\mathrm{e}} =1$ and $\alpha =1$.
We solve \eqref{EQ:LLG-general-pde}--\eqref{EQ:LLG-general-iv} with the initial data given by
\begin{align*}
  \bsm^0(\bsx) &=
  \begin{cases}
    \frac{1}{A^2+|\bsx|^2}\bmat{2A\bsx\\ A^2-|\bsx|^2}\quad
       & \mbox{if } |\bsx|\le 1/2\,,\vspace{4pt}\\
    \bmat{0\\ 0\\ -1} & \mbox{otherwise},
  \end{cases}
\end{align*}
where $A =(1-2|\bsx|)^4/s$, $s =16$.
No external field is applied, i.e.,~$\bsH_{\mathrm{ext}} =\bs 0$.
An exact solution is not known, and we have to rely on the qualitative behavior of the solution in order to judge the validity of the results.
We perform numerical experiments for the space- and time-adaptive Algorithm \ref{ALG:Algo1} and choose the mesh refinement parameter $\theta_{\mathrm{r}}=0.7$ and the mesh coarsening parameter $\theta_{\mathrm{c}} = 0.95$.
Time integration was done with BDF(2).

We compare the adaptive algorithm with uniform discretizations for several mesh sizes in Figure~\ref{FIG:BlowUp-Lagr}.
We observe a drop in $\norm{\nabla\bsm}{L^\infty(\Omega)}$ and energy \eqref{EQ:def_energy} at increasing time instances for decreasing uniform spatial step sizes $h\in\{1/8,1/10,1/12,1/14\}$.
For the finest uniform case and the adaptive case with temporal tolerance $\tol_{\mathrm{t}} =10^{-5}$ and spatial tolerance $\tol_{\mathrm{s}} = 0.4$ we do not observe such a drop over a much longer time period.

The takeaway of these experiments is that the adaptive algorithm produces a stable solution that does not appear to collapse into a singularity and thus requires far less computational effort than the uniform approach. 
This indicates that the error estimator is also useful in non-smooth situations, despite the theory only working under smoothness assumptions on the solution.

\begin{figure}[ht]
\center
\includegraphics[width=4.8cm]{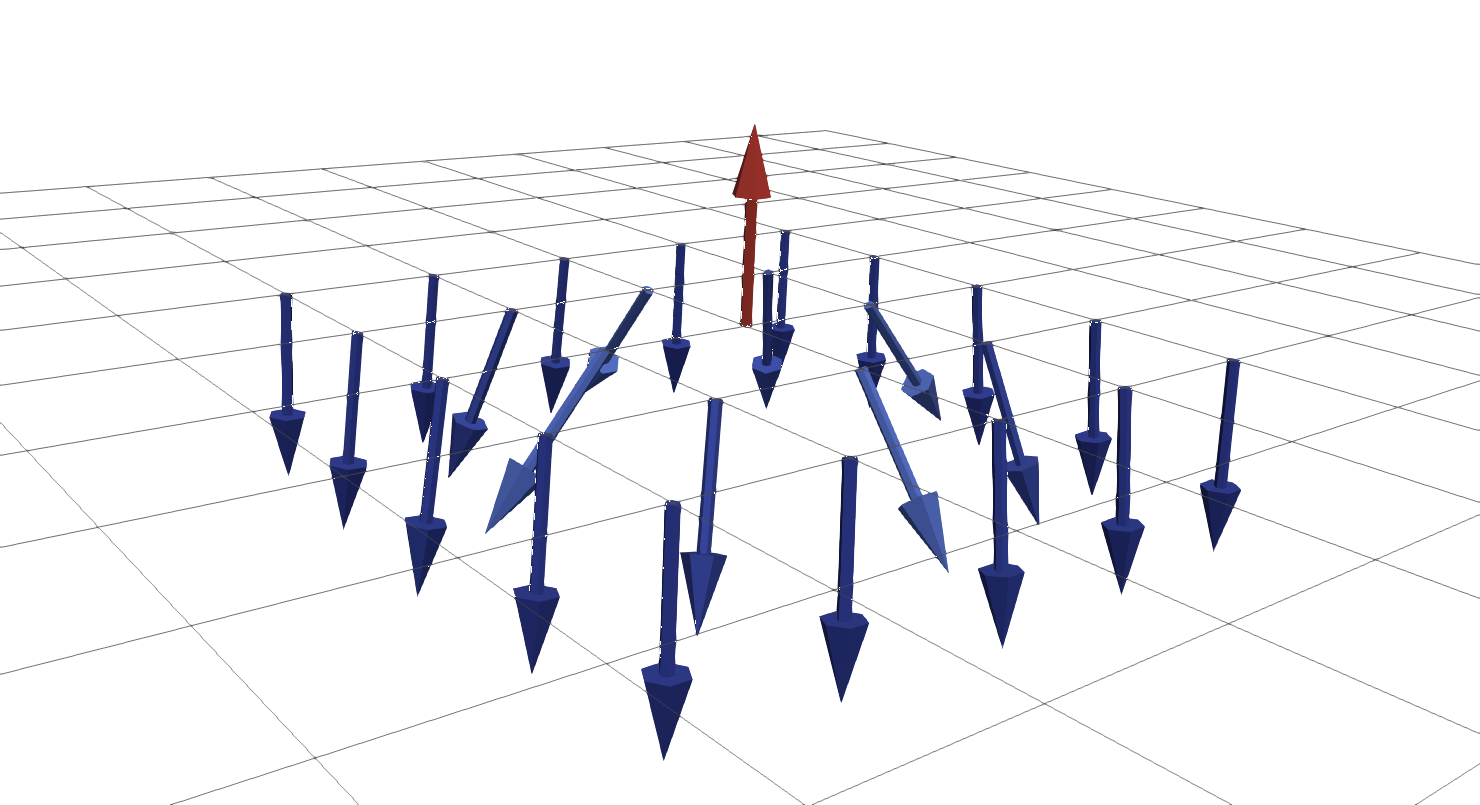} \quad
\includegraphics[width=4.8cm]{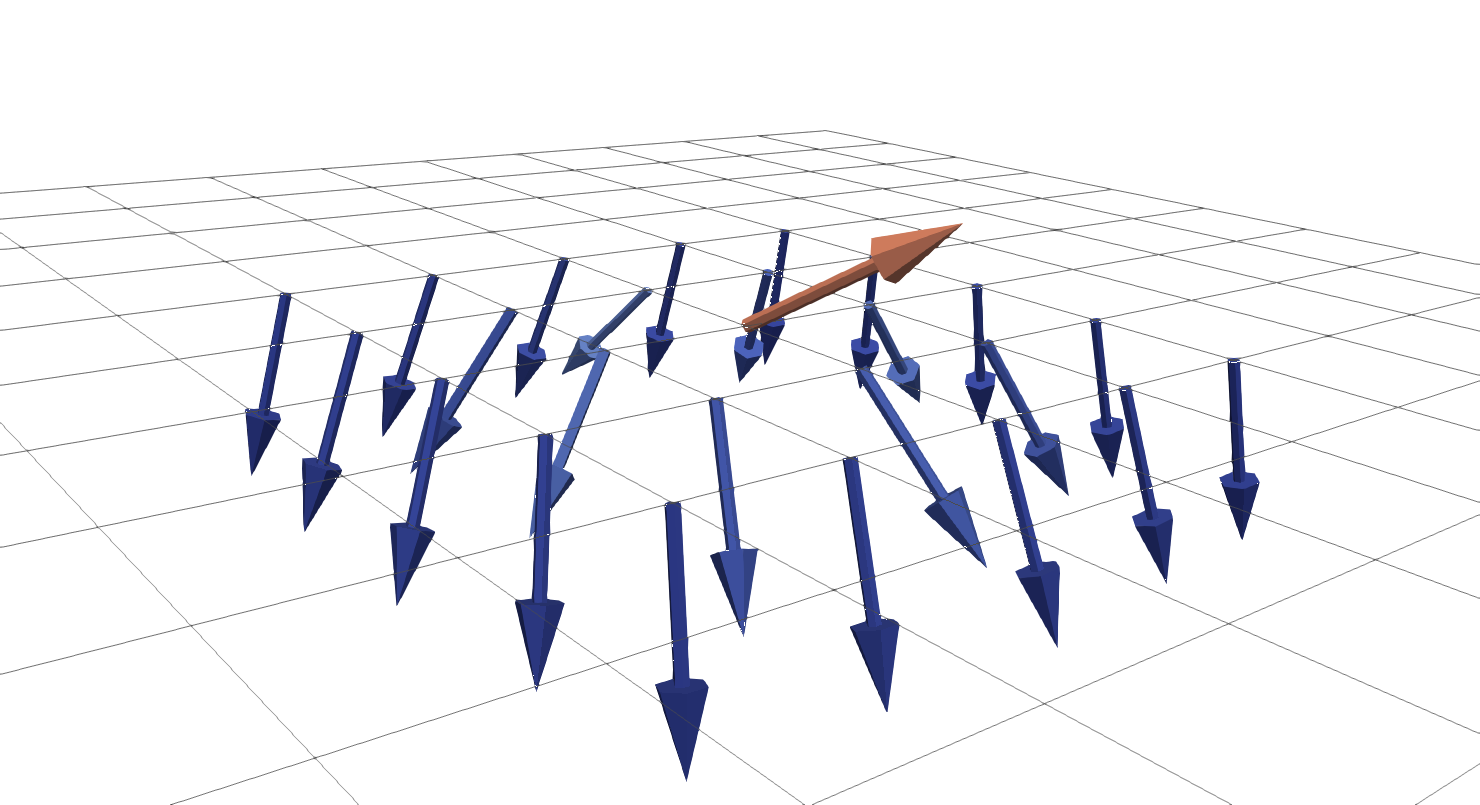} \quad
\includegraphics[width=4.8cm]{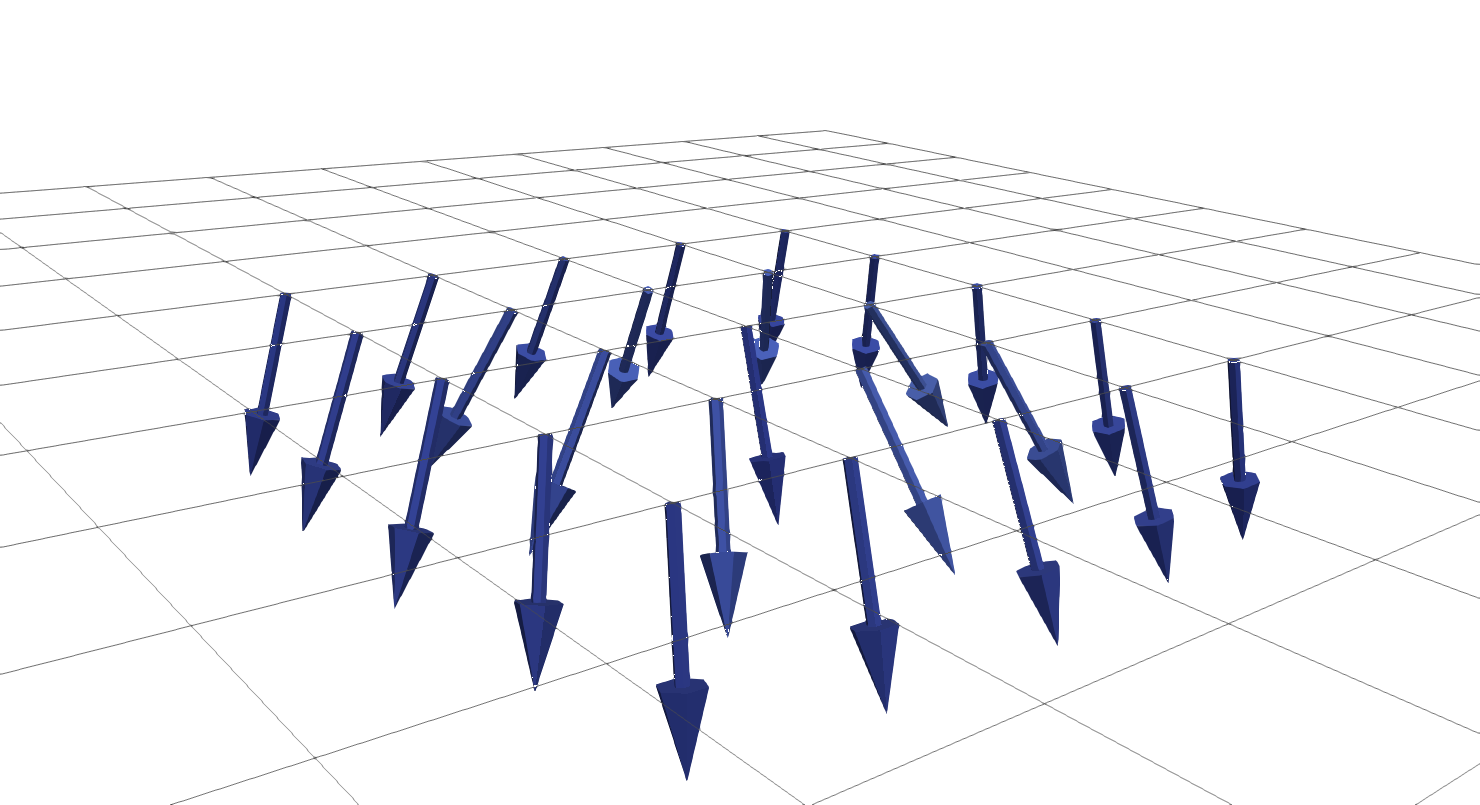}\\
\resizebox{0.48\textwidth}{6.2cm}{%
    % This file was created with tikzplotlib v0.10.1.
\begin{tikzpicture}

\definecolor{darkgray176}{RGB}{176,176,176}
\definecolor{darkorange2551490}{RGB}{255,149,0}
\definecolor{limegreen018569}{RGB}{0,185,69}
\definecolor{orangered255440}{RGB}{255,44,0}
\definecolor{slategray13291151}{RGB}{132,91,151}
\definecolor{teal1293165}{RGB}{12,93,165}

\begin{axis}[
legend cell align={left},
legend style={
  fill opacity=0.8,
  draw opacity=1,
  text opacity=1,
  at={(0.03,0.97)},
  anchor=north west,
  draw=none
},
log basis x={10},
log basis y={10},
tick pos=both,
x grid style={darkgray176},
xmajorgrids,
xlabel={t},
xmin=6.01726408703519e-08, xmax=1.06779600059997,
xmode=log,
xtick style={color=black},
% xtick={1e-10,1e-08,1e-06,0.0001,0.01,1,100,10000},
% xticklabels={
%   $\mathdefault{10^{-10}}$,
%   $\mathdefault{10^{-8}}$,
%   $\mathdefault{10^{-6}}$,
%   $\mathdefault{10^{-4}}$,
%   $\mathdefault{10^{-2}}$,
%   $\mathdefault{10^{0}}$,
%   $\mathdefault{10^{2}}$,
%   $\mathdefault{10^{4}}$
% },
y grid style={darkgray176},
ylabel={$\tau_n$},
ymin=7.03907014285313e-08, ymax=0.0396321523068279,
ymode=log,
ymajorgrids,
ytick style={color=black},
% ytick={1e-09,1e-08,1e-07,1e-06,1e-05,0.0001,0.001,0.01,0.1,1},
% yticklabels={
%   $\mathdefault{10^{-9}}$,
%   $\mathdefault{10^{-8}}$,
%   $\mathdefault{10^{-7}}$,
%   $\mathdefault{10^{-6}}$,
%   $\mathdefault{10^{-5}}$,
%   $\mathdefault{10^{-4}}$,
%   $\mathdefault{10^{-3}}$,
%   $\mathdefault{10^{-2}}$,
%   $\mathdefault{10^{-1}}$,
%   $\mathdefault{10^{0}}$
% }
]
\addplot [slategray13291151]
table {%
1.285e-07 1.285e-07
3.855e-07 2.57e-07
8.994e-07 5.139e-07
1.9274e-06 1.028e-06
3.9834e-06 2.056e-06
8.0944e-06 4.111e-06
1.63174e-05 8.223e-06
3.27674e-05 1.645e-05
6.03874e-05 2.762e-05
8.51274e-05 2.474e-05
0.0001126574 2.753e-05
0.0001450574 3.24e-05
0.0001811674 3.611e-05
0.0002219974 4.083e-05
0.0002679374 4.594e-05
0.0003192074 5.127e-05
0.0003758674 5.666e-05
0.0004376774 6.181e-05
0.0005045874 6.691e-05
0.0005766674 7.208e-05
0.0006540574 7.739e-05
0.0007369474 8.289e-05
0.0008255974 8.865e-05
0.0009203574 9.476e-05
0.0010216574 0.0001013
0.0010469774 2.532e-05
0.0010547734 7.796e-06
0.0010703634 1.559e-05
0.0010871434 1.678e-05
0.0011043134 1.717e-05
0.0011248134 2.05e-05
0.0011486834 2.387e-05
0.0011760134 2.733e-05
0.0012075934 3.158e-05
0.0012441734 3.658e-05
0.0012866034 4.243e-05
0.0013360234 4.942e-05
0.0013938534 5.783e-05
0.0014618834 6.803e-05
0.0015423734 8.049e-05
0.0016379834 9.561e-05
0.0017512834 0.0001133
0.0018835834 0.0001323
0.0020344834 0.0001509
0.0022019834 0.0001675
0.0023840834 0.0001821
0.0025793834 0.0001953
0.0027873834 0.000208
0.0030082834 0.0002209
0.0032424834 0.0002342
0.0034903834 0.0002479
0.0037525834 0.0002622
0.0040295834 0.000277
0.0043219834 0.0002924
0.0046304834 0.0003085
0.0049557834 0.0003253
0.0052986834 0.0003429
0.0056598834 0.0003612
0.0060402834 0.0003804
0.0064406834 0.0004004
0.0068619834 0.0004213
0.0073050834 0.0004431
0.0077709834 0.0004659
0.0082605834 0.0004896
0.0087749834 0.0005144
0.0093153834 0.0005404
0.0098827834 0.0005674
0.0104784834 0.0005957
0.0111036834 0.0006252
0.0117596834 0.000656
0.0124477834 0.0006881
0.0131693834 0.0007216
0.0139258834 0.0007565
0.0147187834 0.0007929
0.0155496834 0.0008309
0.0164202834 0.0008706
0.0173322834 0.000912
0.0182873834 0.0009551
0.0192873834 0.001
0.0203343834 0.001047
0.0214303834 0.001096
0.0225773834 0.001147
0.0237773834 0.0012
0.0250323834 0.001255
0.0263443834 0.001312
0.0277163834 0.001372
0.0291513834 0.001435
0.0306523834 0.001501
0.0322233834 0.001571
0.0338683834 0.001645
0.0355923834 0.001724
0.0374003834 0.001808
0.0392973834 0.001897
0.0412883834 0.001991
0.0433803834 0.002092
0.0455783834 0.002198
0.0478873834 0.002309
0.0503143834 0.002427
0.0528663834 0.002552
0.0555493834 0.002683
0.0583723834 0.002823
0.0613453834 0.002973
0.0644793834 0.003134
0.0677853834 0.003306
0.0712763834 0.003491
0.0749653834 0.003689
0.0788673834 0.003902
0.0829963834 0.004129
0.0873663834 0.00437
0.0919923834 0.004626
0.0968883834 0.004896
0.1020663834 0.005178
0.1075353834 0.005469
0.1133043834 0.005769
0.1193813834 0.006077
0.1257753834 0.006394
0.1324983834 0.006723
0.1395633834 0.007065
0.1469873834 0.007424
0.1547873834 0.0078
0.1629783834 0.008191
0.1715743834 0.008596
0.1805833834 0.009009
0.1900043834 0.009421
0.1998173834 0.009813
0.2099673834 0.01015
0.2203273834 0.01036
0.2307073834 0.01038
0.2408573834 0.01015
0.2506013834 0.009744
0.2600063834 0.009405
0.2698433834 0.009837
0.2796763834 0.009833
0.2909963834 0.01132
0.2991763834 0.00818
0.3053813834 0.006205
0.3108143834 0.005433
0.3188323834 0.008018
0.3229573834 0.004125
0.3277563834 0.004799
0.3331373834 0.005381
0.3393653834 0.006228
0.3476713834 0.008306
0.3567163834 0.009045
0.3671863834 0.01047
0.3794163834 0.01223
0.3930763834 0.01366
0.4084363834 0.01536
0.4249763834 0.01654
0.4432363834 0.01826
0.4621763834 0.01894
0.4838863834 0.02171
0.5000163834 0.01613
};
\end{axis}

\end{tikzpicture}
}
\caption{\label{FIG:BlowUp-solution-Lagr}
Example~\ref{SSEC:example3}: Top row: Solution $\bsm_h$ for the uniform mesh $h=1/12$ and uniform time step $\tau=10^{-3}$ at $t=0$ (left), $t=0.068$ (middle) and $t=0.071$ (right) according to \eqref{EQ:LLG-weak-lambda-bdf}. The singularity appears in the center of the vector field, where $\bsm$ points upward in a very small neighborhood only.
Bottom: Evolution of the time step size in the adaptive computation. 
}
\end{figure}

Three snapshots in time for the case $h=1/12$ illustrate the flipping of the central vector in Figure~\ref{FIG:BlowUp-solution-Lagr} (top row).
The evolution of the time step sizes for the adaptive method corresponding to Figure~\ref{FIG:BlowUp-Lagr} is shown in Figure~\ref{FIG:BlowUp-solution-Lagr} (bottom).

\begin{figure}[ht] \center
    \resizebox{!}{8cm}{%
        \input{new_plots/BlowUp_energy}
    }
 \caption{\label{FIG:BlowUp-Lagr}
 Example~\ref{SSEC:example3}: We compare simulations with uniform meshes and time steps with a fully adaptive simulation. The time step size for the uniform simulations is $\tau=10^{-3}$ and the mesh size ranges from $h=1/8$ to $h=1/14$. The adaptive simulation starts with a mesh size of $h=1/2$.
 The singularity, i.e., the rotation of the entire magnetization towards $[0,0,-1]$, is observed via a sharp decay of the energy \eqref{EQ:def_energy} as well as a spike in the $\norm{\nabla\bsm}{L^\infty}$-norm.
 We observe that only the finest uniform simulation avoids the singularity (that requires a total number of 1\,682\,000 degrees of freedom).
 In contrast, the adaptive simulation achieves the same with just 753\,540 degrees of freedom. 
 } 
\end{figure}

\subsection{Example 4}\label{SSEC:example5}\ 
The following prescribed solution, taken from \cite[Sect.~3.7]{Banas:2005}, simulates a moving local structure.
We take $\Omega =(0,1)\times (0,1)$, $T =0.4$, $C_{\mathrm{e}} =1$, $\alpha =1$ and provide the exact solution
\begin{align*}
    \bsm(t,\bsx) =
	\frac{\sqrt{2}}{2}
	\begin{bmatrix}
		0 \\
		\sqrt{2 - \ez^{-2 D r(t, \bsx)} } \\
		\ez^{-D r(t, \bsx)}
	\end{bmatrix},
\end{align*}
with $r(t,\bsx) = (x_1 - t - t_0)^2 + (x_2 - t - t_0)^2$ and $D =200$, $t_0 =0.3$.
Initial conditions and an external field $\bsH_{\mathrm{ext}}$ have been chosen accordingly.

We perform numerical experiments for adaptive cubic finite elements in space and adaptive time stepping with BDF($4$), where we apply the fully adaptive Algorithm~\ref{ALG:Algo1} with mesh refinement parameter $\theta_{\mathrm{r}} = 0.85$, $\theta_{\mathrm{c}} = 0.9$ and tolerances $\tol_{\mathrm{s}} =10^{-4}$ and $\tol_{\mathrm{t}} =10^{-6}$.

\begin{center}
\begin{figure}[t!]
	\centering
	\begin{subfigure}{0.32\textwidth}
		\centering
		\includegraphics[width=\linewidth]{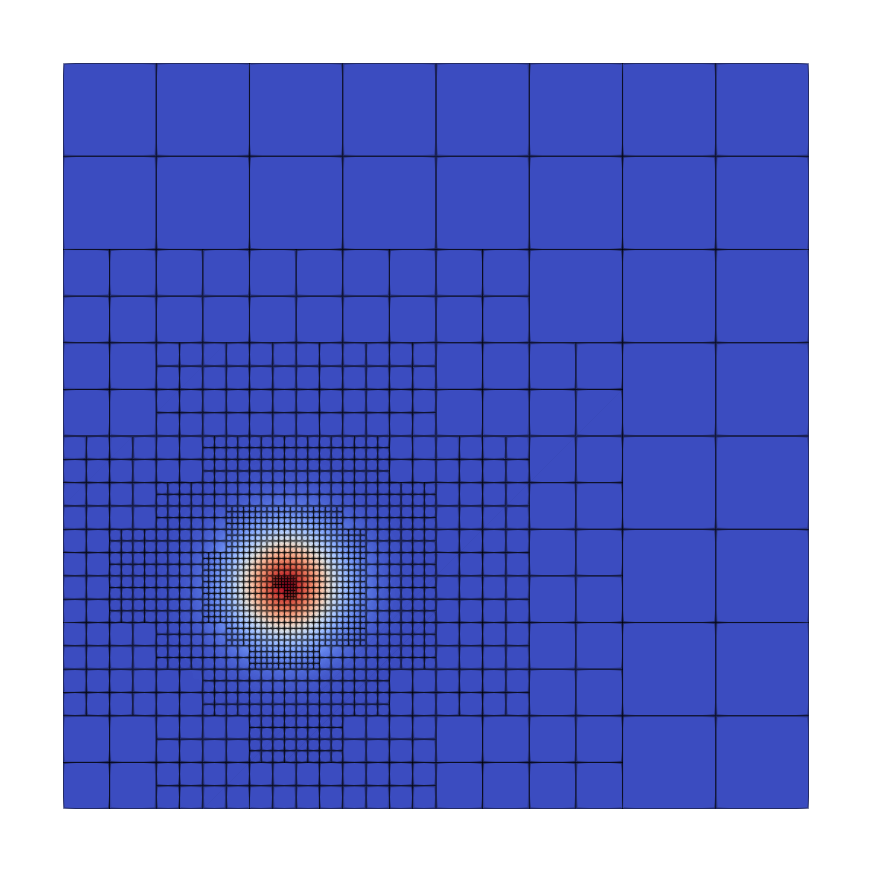}
		\caption{$t = 0$}
	\end{subfigure}
	\begin{subfigure}{0.32\textwidth}
		\centering
		\includegraphics[width=\linewidth]{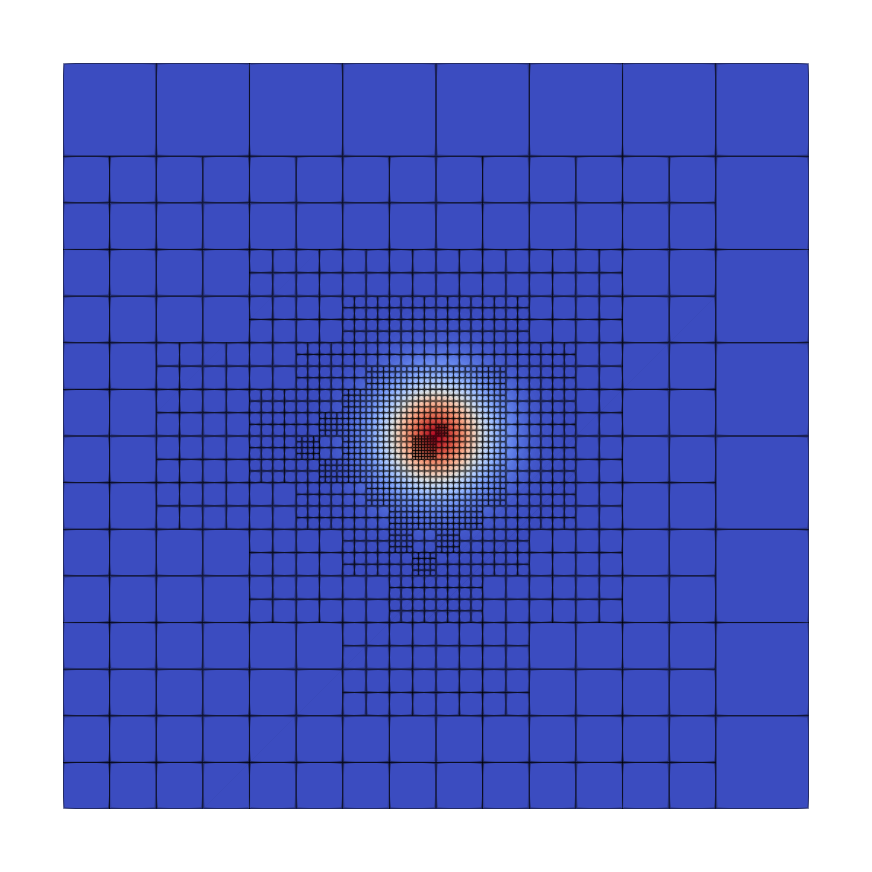}
		\caption{$t = T/2$}
	\end{subfigure}
	\begin{subfigure}{0.32\textwidth}
		\centering
		\includegraphics[width=\linewidth]{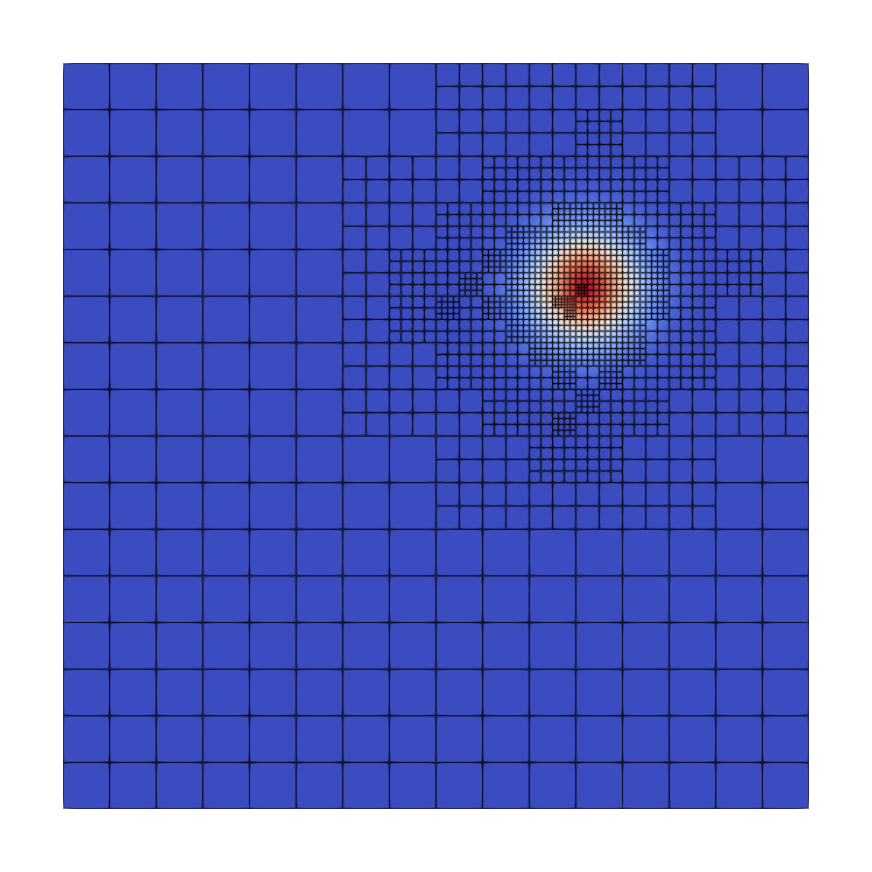}
		\caption{$t = T$}
	\end{subfigure}
	\caption{
        \label{FIG:bump-meshes}
		Example~\ref{SSEC:example5}: Evolution of the spatial mesh using the space and time adaptive Algorithm~\ref{ALG:Algo1}. The color represents the $z$-component of the magnetization.
	}
\end{figure}
\end{center}

In Figure~\ref{FIG:bump-meshes} we show the mesh at different time instances.
The adapted mesh follows the bump's position and adapts the spatial grid accordingly.
A comparison between a uniform and an adaptive computation is presented in Figure~\ref{FIG:bump-uniform-vs-adaptive}.
The adaptive algorithm significantly reduces the error compared to the uniform simulation under comparable computational effort.

\begin{figure}[htb!]
    \centering
    \begin{subfigure}{0.48\textwidth}
        \centering
        \resizebox{\textwidth}{!}{%
            % This file was created with tikzplotlib v0.10.1.
\begin{tikzpicture}

\definecolor{darkgray176}{RGB}{176,176,176}
\definecolor{darkorange2551490}{RGB}{255,149,0}
\definecolor{limegreen018569}{RGB}{0,185,69}
\definecolor{teal1293165}{RGB}{12,93,165}

\begin{axis}[
legend cell align={left},
height=7cm,
width=10cm,
legend style={fill opacity=0.8, draw opacity=1, text opacity=1, at={(0.91,0.5)}, anchor=east, draw=none},
log basis y={10},
tick pos=both,
xlabel={$t$},
xmajorgrids,
ylabel={$\norm{\grad(\bsm(t_n)-\bsm_h^n)}{L^2(\Omega)}$},
ymajorgrids,
ymode=log
]
\addplot [teal1293165, dashed, mark=*, mark size=2, mark options={solid}]
table {%
0 0.012315
0.0025641 0.0070866
0.0051282 0.0039277
0.0076923 0.0031778
0.0102564 0.0030879
0.0128205 0.0030717
0.0153846 0.0030881
0.0179487 0.003108
0.0205128 0.0031017
0.0230769 0.0030797
0.025641 0.0030487
0.0282051 0.0030247
0.0307692 0.0030241
0.0333333 0.0030547
0.0358974 0.0031075
0.0384615 0.0031619
0.0410256 0.0031953
0.0435897 0.003193
0.0461538 0.0031535
0.0487179 0.0030894
0.051282 0.0030234
0.0538461 0.0029794
0.0564102 0.0029717
0.0589743 0.0029976
0.0615384 0.0030394
0.0641025 0.003073
0.0666666 0.0030792
0.0692307 0.0030512
0.0717948 0.0029972
0.0743589 0.0029376
0.076923 0.0028968
0.0794871 0.0028922
0.0820512 0.002925
0.0846153 0.0029795
0.0871794 0.0030312
0.0897435 0.0030583
0.0923076 0.0030498
0.0948717 0.00301
0.0974358 0.0029566
0.0999999 0.0029137
0.102564 0.0029018
0.1051281 0.0029273
0.1076922 0.0029791
0.1102563 0.0030342
0.1128204 0.0030693
0.1153845 0.0030703
0.1179486 0.0030369
0.1205127 0.0029836
0.1230768 0.0029337
0.1256409 0.0029099
0.128205 0.0029235
0.1307691 0.0029682
0.1333332 0.0030235
0.1358973 0.0030653
0.1384614 0.0030758
0.1410255 0.0030507
0.1435896 0.0030004
0.1461537 0.0029465
0.1487178 0.0029131
0.1512819 0.0029157
0.153846 0.002953
0.1564101 0.0030078
0.1589742 0.0030557
0.1615383 0.0030762
0.1641024 0.0030609
0.1666665 0.003016
0.1692306 0.0029606
0.1717947 0.0029193
0.1743588 0.0029109
0.1769229 0.0029395
0.179487 0.0029917
0.1820511 0.0030441
0.1846152 0.0030739
0.1871793 0.0030688
0.1897434 0.0030309
0.1923075 0.0029761
0.1948716 0.0029283
0.1974357 0.0029094
0.1999998 0.0029281
0.2025639 0.0029757
0.205128 0.0030305
0.2076921 0.0030687
0.2102562 0.003074
0.2128203 0.0030444
0.2153844 0.0029921
0.2179485 0.0029398
0.220512600000001 0.002911
0.223076700000001 0.0029191
0.225640800000001 0.0029603
0.228204900000001 0.0030156
0.230769000000001 0.0030607
0.233333100000001 0.0030763
0.235897200000001 0.003056
0.238461300000001 0.0030082
0.241025400000001 0.0029533
0.243589500000001 0.0029158
0.246153600000001 0.0029129
0.248717700000001 0.002946
0.251281800000001 0.0029998
0.253845900000001 0.0030502
0.256410000000001 0.0030755
0.258974100000001 0.0030652
0.261538200000001 0.0030236
0.264102300000001 0.0029683
0.266666400000001 0.0029235
0.269230500000001 0.0029098
0.271794600000001 0.0029336
0.274358700000001 0.0029837
0.276922800000001 0.0030375
0.279486900000001 0.0030717
0.282051000000001 0.0030718
0.284615100000001 0.0030379
0.287179200000001 0.0029841
0.289743300000001 0.0029339
0.292307400000001 0.0029099
0.294871500000001 0.0029233
0.297435600000001 0.0029679
0.299999700000001 0.0030233
0.302563800000001 0.003065
0.305127900000001 0.0030756
0.307692000000001 0.0030505
0.310256100000001 0.0030003
0.312820200000001 0.0029464
0.315384300000001 0.0029131
0.317948400000001 0.0029157
0.320512500000001 0.002953
0.323076600000001 0.0030078
0.325640700000001 0.0030558
0.328204800000001 0.0030763
0.330768900000001 0.003061
0.333333000000001 0.0030161
0.335897100000001 0.0029607
0.338461200000001 0.0029194
0.341025300000001 0.002911
0.343589400000001 0.0029396
0.346153500000001 0.0029918
0.348717600000001 0.0030441
0.351281700000001 0.003074
0.353845800000001 0.0030689
0.356409900000001 0.003031
0.358974000000001 0.0029762
0.361538100000001 0.0029284
0.364102200000001 0.0029095
0.366666300000001 0.0029282
0.369230400000001 0.0029758
0.371794500000001 0.0030306
0.374358600000001 0.0030688
0.376922700000001 0.0030741
0.379486800000001 0.0030445
0.382050900000001 0.0029922
0.384615000000001 0.0029399
0.387179100000001 0.0029111
0.389743200000001 0.0029192
0.392307300000001 0.0029604
0.394871400000001 0.0030157
0.397435500000001 0.0030608
};
\addlegendentry{uni}
\addplot [limegreen018569, dashed, mark=square*, mark size=2, mark options={solid}]
table {%
0 0.00012553
6.529e-06 0.0001253
1.9587e-05 0.0001253
4.5703e-05 0.00012545
8.2637e-05 0.00012531
0.000126958 0.00012519
0.000180143 0.00012513
0.000243965 0.00012508
0.000320551 0.00012502
0.000412454 0.00012496
0.000522734 0.00012489
0.000655074 0.00012483
0.000813884 0.00012476
0.001004454 0.00012469
0.001233134 0.00012463
0.001507554 0.0001246
0.001836864 0.00012461
0.002232034 0.0001247
0.002706234 0.00012491
0.003275274 0.0001253
0.003958124 0.00012599
0.004777544 0.0001271
0.005760844 0.00012883
0.006940844 0.00012487
0.008356844 0.00012808
0.010055944 0.00013295
0.012094944 0.00014052
0.014541744 0.00014817
0.017477844 0.0001653
0.021001144 0.00022969
0.025088244 0.0003068
0.028457144 0.00035047
0.031096744 0.0003672
0.033545744 0.00036576
0.036063744 0.00036171
0.038606344 0.00035599
0.041259944 0.00034605
0.043935744 0.00033616
0.046707944 0.00032
0.049888544 0.00030233
0.053286344 0.00029841
0.056544844 0.00030325
0.059610744 0.00031299
0.062512644 0.00031762
0.065462044 0.00032536
0.068568544 0.00033189
0.071663944 0.00033483
0.074592744 0.00033694
0.077445744 0.00033789
0.080488144 0.0003384
0.083721544 0.00034132
0.086920744 0.0003455
0.090026544 0.00034847
0.092978144 0.00034897
0.095754944 0.00034487
0.098467044 0.00034364
0.101142744 0.00034137
0.103829244 0.00033588
0.106564744 0.00033281
0.109444344 0.00032809
0.112710144 0.0003279
0.116104244 0.00033958
0.119293144 0.00035113
0.122288244 0.0003616
0.125187144 0.00036794
0.128188044 0.00036938
0.131334344 0.00037296
0.134462544 0.00037341
0.137428344 0.00037275
0.140264844 0.00037099
0.143206444 0.00036754
0.146360744 0.00036414
0.149570544 0.0003626
0.152734844 0.00036253
0.155740644 0.00036217
0.158527644 0.00035826
0.161220044 0.00035757
0.163878444 0.0003561
0.166553444 0.00035138
0.169278544 0.00034905
0.172148344 0.00034546
0.175406444 0.00034633
0.178794644 0.0003566
0.181981444 0.00036845
0.184979544 0.00037902
0.187880844 0.00038558
0.190880744 0.00038737
0.194023344 0.00039085
0.197146044 0.00039139
0.200106244 0.00039079
0.202938844 0.00038905
0.205879444 0.00038566
0.209035844 0.00038169
0.212246944 0.0003801
0.215409344 0.00037995
0.218409444 0.00037947
0.221186644 0.00037803
0.223868944 0.00037734
0.226519644 0.0003735
0.229201544 0.00036786
0.231980444 0.00036548
0.234923544 0.00036246
0.238113044 0.00036328
0.241302644 0.00036712
0.244397744 0.0003593
0.247448744 0.0003707
0.250449144 0.00037784
0.253474944 0.00038163
0.256607244 0.00038697
0.259764744 0.00038949
0.262738244 0.00039032
0.265574644 0.00038952
0.268514244 0.0003866
0.271653344 0.00038247
0.274858144 0.00038018
0.278024144 0.00037916
0.281029544 0.00037799
0.283812644 0.00037608
0.286494244 0.0003752
0.289140244 0.00037118
0.291816544 0.00036479
0.294587644 0.00036225
0.297524144 0.00035902
0.300710144 0.00035966
0.303895144 0.00036336
0.306987244 0.00037212
0.310039044 0.00038329
0.313040844 0.00039044
0.316066744 0.00039442
0.319197844 0.00039885
0.322353044 0.00040156
0.325323844 0.00040253
0.328158244 0.00040185
0.331097344 0.00039908
0.334237344 0.00039507
0.337442844 0.00039284
0.340607944 0.00039181
0.343609844 0.00039062
0.346385944 0.00038869
0.349060444 0.00038777
0.351701444 0.00038384
0.354373444 0.00037757
0.357139244 0.0003751
0.360069744 0.00037191
0.363253544 0.00037246
0.366442444 0.00037605
0.369539944 0.00037747
0.372593544 0.00038868
0.375592544 0.00039586
0.378614244 0.00039984
0.381743344 0.00040523
0.384899444 0.00040788
0.387870944 0.00040883
0.390704644 0.00040816
0.393642144 0.00040569
0.396781144 0.00040043
};
\addlegendentry{adapt}
\end{axis}
\end{tikzpicture}
        }
        \caption{$H^1(\Omega)$-error}
    \end{subfigure}
    \hfill
    \begin{subfigure}{0.48\textwidth}
        \centering
        \resizebox{\textwidth}{!}{%
            % This file was created with tikzplotlib v0.10.1.
\begin{tikzpicture}

\definecolor{darkgray176}{RGB}{176,176,176}
\definecolor{darkorange2551490}{RGB}{255,149,0}
\definecolor{limegreen018569}{RGB}{0,185,69}
\definecolor{teal1293165}{RGB}{12,93,165}

\begin{axis}[
legend cell align={left},
height=7cm,
width=10cm,
legend style={fill opacity=0.8, draw opacity=1, text opacity=1, at={(0.91,0.5)}, anchor=east, draw=none},
log basis y={10},
tick pos=both,
xlabel={$t$},
xmajorgrids,
ylabel={$\norm{|\bsm_h^n| - 1}{L^\infty(\Omega)}$},
ymajorgrids,
ymode=log
]
\addplot [teal1293165, dashed, mark=*, mark size=2, mark options={solid}]
table {%
0 0.00013767
0.0025641 0.00014686
0.0051282 0.00012811
0.0076923 0.00011433
0.0102564 8.0066e-05
0.0128205 0.00017084
0.0153846 0.00022396
0.0179487 0.00022913
0.0205128 0.00022259
0.0230769 0.0002208
0.025641 0.00021991
0.0282051 0.00021432
0.0307692 0.00020246
0.0333333 0.00019139
0.0358974 0.00019972
0.0384615 0.00020959
0.0410256 0.00021238
0.0435897 0.00020695
0.0461538 0.00019402
0.0487179 0.00017585
0.051282 0.00017149
0.0538461 0.00017107
0.0564102 0.00017073
0.0589743 0.00017045
0.0615384 0.00017019
0.0641025 0.00016995
0.0666666 0.00016974
0.0692307 0.00016957
0.0717948 0.00016944
0.0743589 0.00016932
0.076923 0.00016923
0.0794871 0.00016915
0.0820512 0.0001691
0.0846153 0.00016905
0.0871794 0.00016902
0.0897435 0.00016899
0.0923076 0.00016897
0.0948717 0.00016896
0.0974358 0.00016895
0.0999999 0.00016894
0.102564 0.00016894
0.1051281 0.00016893
0.1076922 0.00016893
0.1102563 0.00016893
0.1128204 0.00016893
0.1153845 0.00016893
0.1179486 0.00016893
0.1205127 0.00016893
0.1230768 0.00016892
0.1256409 0.00016892
0.128205 0.00016892
0.1307691 0.00016892
0.1333332 0.00016892
0.1358973 0.00016892
0.1384614 0.00016892
0.1410255 0.00016892
0.1435896 0.00016892
0.1461537 0.00016892
0.1487178 0.00016892
0.1512819 0.00016892
0.153846 0.00016892
0.1564101 0.00016892
0.1589742 0.00016892
0.1615383 0.00016892
0.1641024 0.00016892
0.1666665 0.00016892
0.1692306 0.00016892
0.1717947 0.00016892
0.1743588 0.00016892
0.1769229 0.00016892
0.179487 0.00016892
0.1820511 0.00016892
0.1846152 0.00016892
0.1871793 0.00016892
0.1897434 0.00016892
0.1923075 0.00016892
0.1948716 0.00016892
0.1974357 0.00016892
0.1999998 0.00016892
0.2025639 0.00016892
0.205128 0.00016892
0.2076921 0.00016892
0.2102562 0.00016892
0.2128203 0.00016892
0.2153844 0.00016892
0.2179485 0.00016892
0.220512600000001 0.00016892
0.223076700000001 0.00016892
0.225640800000001 0.00016892
0.228204900000001 0.00016892
0.230769000000001 0.00016892
0.233333100000001 0.00016892
0.235897200000001 0.00016892
0.238461300000001 0.00016892
0.241025400000001 0.00016892
0.243589500000001 0.00016892
0.246153600000001 0.00016892
0.248717700000001 0.00016892
0.251281800000001 0.00016892
0.253845900000001 0.00016892
0.256410000000001 0.00016892
0.258974100000001 0.00016892
0.261538200000001 0.00016892
0.264102300000001 0.00016892
0.266666400000001 0.00016892
0.269230500000001 0.00016892
0.271794600000001 0.00016892
0.274358700000001 0.00016892
0.276922800000001 0.00016892
0.279486900000001 0.00016892
0.282051000000001 0.00016892
0.284615100000001 0.00016892
0.287179200000001 0.00016892
0.289743300000001 0.00016892
0.292307400000001 0.00016892
0.294871500000001 0.00016892
0.297435600000001 0.00016892
0.299999700000001 0.00016892
0.302563800000001 0.00016892
0.305127900000001 0.00016892
0.307692000000001 0.00016892
0.310256100000001 0.00016892
0.312820200000001 0.00016892
0.315384300000001 0.00016892
0.317948400000001 0.00016892
0.320512500000001 0.00016892
0.323076600000001 0.00016892
0.325640700000001 0.00016892
0.328204800000001 0.00016892
0.330768900000001 0.00016892
0.333333000000001 0.00016892
0.335897100000001 0.00016892
0.338461200000001 0.00016892
0.341025300000001 0.00016892
0.343589400000001 0.00016892
0.346153500000001 0.00016892
0.348717600000001 0.00016892
0.351281700000001 0.00016892
0.353845800000001 0.00016892
0.356409900000001 0.00016892
0.358974000000001 0.00016892
0.361538100000001 0.00016892
0.364102200000001 0.00016892
0.366666300000001 0.00016892
0.369230400000001 0.00016892
0.371794500000001 0.00016892
0.374358600000001 0.00016892
0.376922700000001 0.00016892
0.379486800000001 0.00016892
0.382050900000001 0.00016892
0.384615000000001 0.00016892
0.387179100000001 0.00016892
0.389743200000001 0.00016892
0.392307300000001 0.00016892
0.394871400000001 0.00016892
0.397435500000001 0.00016892
};
\addlegendentry{uni}
\addplot [limegreen018569, dashed, mark=square*, mark size=2, mark options={solid}]
table {%
0 1.3966e-06
6.529e-06 1.3966e-06
1.9587e-05 1.3975e-06
4.5703e-05 1.3988e-06
8.2637e-05 1.3996e-06
0.000126958 1.4013e-06
0.000180143 1.4033e-06
0.000243965 1.4057e-06
0.000320551 1.4088e-06
0.000412454 1.4124e-06
0.000522734 1.4165e-06
0.000655074 1.4215e-06
0.000813884 1.4276e-06
0.001004454 1.4349e-06
0.001233134 1.4436e-06
0.001507554 1.454e-06
0.001836864 1.4663e-06
0.002232034 1.4811e-06
0.002706234 1.4986e-06
0.003275274 1.5192e-06
0.003958124 1.5434e-06
0.004777544 1.5714e-06
0.005760844 1.6034e-06
0.006940844 1.64e-06
0.008356844 1.6815e-06
0.010055944 1.8309e-06
0.012094944 3.6014e-06
0.014541744 8.3758e-06
0.017477844 2.0824e-05
0.021001144 4.9014e-05
0.025088244 7.7171e-05
0.028457144 9.2556e-05
0.031096744 9.9185e-05
0.033545744 0.00010066
0.036063744 0.00010044
0.038606344 9.9853e-05
0.041259944 9.9057e-05
0.043935744 9.7629e-05
0.046707944 9.4778e-05
0.049888544 9.0092e-05
0.053286344 8.4631e-05
0.056544844 7.9693e-05
0.059610744 7.5955e-05
0.062512644 7.3131e-05
0.065462044 7.0787e-05
0.068568544 6.8814e-05
0.071663944 6.7361e-05
0.074592744 6.6457e-05
0.077445744 6.6049e-05
0.080488144 6.5956e-05
0.083721544 6.6127e-05
0.086920744 6.6664e-05
0.090026544 6.7083e-05
0.092978144 6.748e-05
0.095754944 6.7817e-05
0.098467044 6.8091e-05
0.101142744 6.8322e-05
0.103829244 6.8527e-05
0.106564744 6.8725e-05
0.109444344 6.8955e-05
0.112710144 6.9203e-05
0.116104244 6.941e-05
0.119293144 6.7307e-05
0.122288244 6.7392e-05
0.125187144 7.2765e-05
0.128188044 7.5492e-05
0.131334344 7.6702e-05
0.134462544 7.6947e-05
0.137428344 7.5547e-05
0.140264844 7.3026e-05
0.143206444 6.951e-05
0.146360744 6.7542e-05
0.149570544 6.7544e-05
0.152734844 6.7545e-05
0.155740644 6.7546e-05
0.158527644 6.7546e-05
0.161220044 6.7546e-05
0.163878444 6.7546e-05
0.166553444 6.7546e-05
0.169278544 6.7546e-05
0.172148344 6.7546e-05
0.175406444 6.7546e-05
0.178794644 6.7546e-05
0.181981444 6.7546e-05
0.184979544 6.7546e-05
0.187880844 7.1257e-05
0.190880744 7.4157e-05
0.194023344 7.6066e-05
0.197146044 7.6171e-05
0.200106244 7.4694e-05
0.202938844 7.2407e-05
0.205879444 6.9046e-05
0.209035844 6.7546e-05
0.212246944 6.7546e-05
0.215409344 6.7546e-05
0.218409444 6.7546e-05
0.221186644 6.7546e-05
0.223868944 6.7546e-05
0.226519644 6.7546e-05
0.229201544 6.7546e-05
0.231980444 6.7546e-05
0.234923544 6.7546e-05
0.238113044 6.7546e-05
0.241302644 6.7546e-05
0.244397744 5.864e-05
0.247448744 5.9886e-05
0.250449144 6.5967e-05
0.253474944 6.9967e-05
0.256607244 7.2643e-05
0.259764744 7.2964e-05
0.262738244 7.2087e-05
0.265574644 7.0586e-05
0.268514244 6.7214e-05
0.271653344 6.303e-05
0.274858144 5.864e-05
0.278024144 5.864e-05
0.281029544 5.864e-05
0.283812644 5.864e-05
0.286494244 5.864e-05
0.289140244 5.864e-05
0.291816544 5.864e-05
0.294587644 5.864e-05
0.297524144 5.864e-05
0.300710144 5.864e-05
0.303895144 5.864e-05
0.306987244 5.864e-05
0.310039044 5.9622e-05
0.313040844 6.5718e-05
0.316066744 6.9979e-05
0.319197844 7.2683e-05
0.322353044 7.317e-05
0.325323844 7.2432e-05
0.328158244 7.112e-05
0.331097344 6.7518e-05
0.334237344 6.3324e-05
0.337442844 5.8716e-05
0.340607944 5.864e-05
0.343609844 5.864e-05
0.346385944 5.864e-05
0.349060444 5.864e-05
0.351701444 5.864e-05
0.354373444 5.864e-05
0.357139244 5.864e-05
0.360069744 5.864e-05
0.363253544 5.864e-05
0.366442444 5.864e-05
0.369539944 5.864e-05
0.372593544 5.9789e-05
0.375592544 6.5908e-05
0.378614244 7.0264e-05
0.381743344 7.2895e-05
0.384899444 7.3398e-05
0.387870944 7.2687e-05
0.390704644 7.1368e-05
0.393642144 6.7743e-05
0.396781144 6.351e-05
};
\addlegendentry{adapt}
\end{axis}
\end{tikzpicture}
        }
        \caption{Normalization error}
    \end{subfigure}
    \vspace{0.5cm}
    \begin{subfigure}{0.48\textwidth}
        \centering
        \resizebox{\textwidth}{!}{%
            % This file was created with tikzplotlib v0.10.1.
\begin{tikzpicture}

\definecolor{darkgray176}{RGB}{176,176,176}
\definecolor{darkorange2551490}{RGB}{255,149,0}
\definecolor{limegreen018569}{RGB}{0,185,69}
\definecolor{teal1293165}{RGB}{12,93,165}

\begin{axis}[
legend cell align={left},
height=7cm,
width=10cm,
legend style={fill opacity=0.8, draw opacity=1, text opacity=1, at={(0.91,0.5)}, anchor=east, draw=none},
log basis y={10},
tick pos=both,
xlabel={$t$},
xmajorgrids,
ylabel={$N_x$},
ymajorgrids,
ymode=log
]
\addplot [teal1293165, dashed, mark=*, mark size=2, mark options={solid}]
table {%
0.0025641 64516
0.0051282 64516
0.0076923 64516
0.0102564 64516
0.0128205 64516
0.0153846 64516
0.0179487 64516
0.0205128 64516
0.0230769 64516
0.025641 64516
0.0282051 64516
0.0307692 64516
0.0333333 64516
0.0358974 64516
0.0384615 64516
0.0410256 64516
0.0435897 64516
0.0461538 64516
0.0487179 64516
0.051282 64516
0.0538461 64516
0.0564102 64516
0.0589743 64516
0.0615384 64516
0.0641025 64516
0.0666666 64516
0.0692307 64516
0.0717948 64516
0.0743589 64516
0.076923 64516
0.0794871 64516
0.0820512 64516
0.0846153 64516
0.0871794 64516
0.0897435 64516
0.0923076 64516
0.0948717 64516
0.0974358 64516
0.0999999 64516
0.102564 64516
0.1051281 64516
0.1076922 64516
0.1102563 64516
0.1128204 64516
0.1153845 64516
0.1179486 64516
0.1205127 64516
0.1230768 64516
0.1256409 64516
0.128205 64516
0.1307691 64516
0.1333332 64516
0.1358973 64516
0.1384614 64516
0.1410255 64516
0.1435896 64516
0.1461537 64516
0.1487178 64516
0.1512819 64516
0.153846 64516
0.1564101 64516
0.1589742 64516
0.1615383 64516
0.1641024 64516
0.1666665 64516
0.1692306 64516
0.1717947 64516
0.1743588 64516
0.1769229 64516
0.179487 64516
0.1820511 64516
0.1846152 64516
0.1871793 64516
0.1897434 64516
0.1923075 64516
0.1948716 64516
0.1974357 64516
0.1999998 64516
0.2025639 64516
0.205128 64516
0.2076921 64516
0.2102562 64516
0.2128203 64516
0.2153844 64516
0.2179485 64516
0.220512600000001 64516
0.223076700000001 64516
0.225640800000001 64516
0.228204900000001 64516
0.230769000000001 64516
0.233333100000001 64516
0.235897200000001 64516
0.238461300000001 64516
0.241025400000001 64516
0.243589500000001 64516
0.246153600000001 64516
0.248717700000001 64516
0.251281800000001 64516
0.253845900000001 64516
0.256410000000001 64516
0.258974100000001 64516
0.261538200000001 64516
0.264102300000001 64516
0.266666400000001 64516
0.269230500000001 64516
0.271794600000001 64516
0.274358700000001 64516
0.276922800000001 64516
0.279486900000001 64516
0.282051000000001 64516
0.284615100000001 64516
0.287179200000001 64516
0.289743300000001 64516
0.292307400000001 64516
0.294871500000001 64516
0.297435600000001 64516
0.299999700000001 64516
0.302563800000001 64516
0.305127900000001 64516
0.307692000000001 64516
0.310256100000001 64516
0.312820200000001 64516
0.315384300000001 64516
0.317948400000001 64516
0.320512500000001 64516
0.323076600000001 64516
0.325640700000001 64516
0.328204800000001 64516
0.330768900000001 64516
0.333333000000001 64516
0.335897100000001 64516
0.338461200000001 64516
0.341025300000001 64516
0.343589400000001 64516
0.346153500000001 64516
0.348717600000001 64516
0.351281700000001 64516
0.353845800000001 64516
0.356409900000001 64516
0.358974000000001 64516
0.361538100000001 64516
0.364102200000001 64516
0.366666300000001 64516
0.369230400000001 64516
0.371794500000001 64516
0.374358600000001 64516
0.376922700000001 64516
0.379486800000001 64516
0.382050900000001 64516
0.384615000000001 64516
0.387179100000001 64516
0.389743200000001 64516
0.392307300000001 64516
0.394871400000001 64516
0.397435500000001 64516
0.399999600000001 64516
};
\addlegendentry{uni}
\addplot [limegreen018569, dashed, mark=square*, mark size=2, mark options={solid}]
table {%
6.529e-06 52408
1.9587e-05 52408
4.5703e-05 52408
8.2637e-05 52408
0.000126958 52408
0.000180143 52408
0.000243965 52408
0.000320551 52408
0.000412454 52408
0.000522734 52408
0.000655074 52408
0.000813884 52408
0.001004454 52408
0.001233134 52408
0.001507554 52408
0.001836864 52408
0.002232034 52408
0.002706234 52408
0.003275274 52408
0.003958124 52408
0.004777544 52408
0.005760844 52408
0.006940844 52408
0.008356844 53816
0.010055944 53816
0.012094944 53816
0.014541744 53816
0.017477844 54792
0.021001144 56088
0.025088244 55448
0.028457144 55448
0.031096744 55336
0.033545744 56744
0.036063744 58472
0.038606344 58472
0.041259944 57608
0.043935744 58692
0.046707944 57604
0.049888544 57716
0.053286344 58580
0.056544844 56420
0.059610744 56420
0.062512644 57828
0.065462044 59556
0.068568544 59556
0.071663944 59556
0.074592744 59444
0.077445744 59556
0.080488144 59988
0.083721544 60852
0.086920744 62964
0.090026544 56996
0.092978144 57860
0.095754944 59268
0.098467044 62076
0.101142744 62076
0.103829244 60556
0.106564744 59260
0.109444344 59260
0.112710144 60236
0.116104244 60124
0.119293144 59052
0.122288244 57644
0.125187144 57964
0.128188044 58828
0.131334344 59692
0.134462544 58828
0.137428344 58828
0.140264844 59260
0.143206444 60236
0.146360744 60668
0.149570544 59820
0.152734844 60684
0.155740644 61548
0.158527644 60028
0.161220044 61756
0.163878444 61756
0.166553444 60780
0.169278544 59912
0.172148344 59368
0.175406444 60344
0.178794644 61208
0.181981444 58296
0.184979544 58184
0.187880844 58504
0.190880744 59368
0.194023344 60232
0.197146044 58504
0.200106244 58504
0.202938844 58936
0.205879444 59912
0.209035844 60344
0.212246944 62868
0.215409344 63732
0.218409444 64596
0.221186644 61236
0.223868944 62964
0.226519644 62532
0.229201544 63396
0.231980444 61988
0.234923544 60900
0.238113044 61876
0.241302644 62100
0.244397744 60692
0.247448744 60692
0.250449144 60692
0.253474944 60468
0.256607244 61332
0.259764744 60468
0.262738244 60468
0.265574644 60900
0.268514244 61876
0.271653344 62308
0.274858144 61460
0.278024144 62324
0.281029544 63188
0.283812644 61668
0.286494244 62532
0.289140244 62100
0.291816544 62964
0.294587644 61988
0.297524144 60900
0.300710144 61876
0.303895144 62100
0.306987244 60692
0.310039044 61124
0.313040844 61124
0.316066744 60900
0.319197844 61764
0.322353044 60036
0.325323844 60036
0.328158244 60468
0.331097344 61444
0.334237344 61876
0.337442844 62804
0.340607944 63668
0.343609844 64532
0.346385944 61172
0.349060444 62868
0.351701444 62436
0.354373444 63300
0.357139244 61892
0.360069744 60804
0.363253544 61780
0.366442444 62004
0.369539944 60596
0.372593544 60596
0.375592544 60596
0.378614244 60372
0.381743344 61236
0.384899444 60372
0.387870944 60372
0.390704644 60804
0.393642144 61780
0.396781144 61572
0.399999944 60276
};
\addlegendentry{adapt}
\end{axis}
\end{tikzpicture}
        }
        \caption{Degrees of freedom}
    \end{subfigure}
    \hfill
    \begin{subfigure}{0.48\textwidth}
        \centering
        \resizebox{\textwidth}{!}{%
            % This file was created with tikzplotlib v0.10.1.
\begin{tikzpicture}

\definecolor{darkgray176}{RGB}{176,176,176}
\definecolor{darkorange2551490}{RGB}{255,149,0}
\definecolor{limegreen018569}{RGB}{0,185,69}
\definecolor{teal1293165}{RGB}{12,93,165}

\begin{axis}[
legend cell align={left},
height=7cm,
width=10cm,
legend style={fill opacity=0.8, draw opacity=1, text opacity=1, at={(0.91,0.5)}, anchor=east, draw=none},
log basis y={10},
tick pos=both,
xlabel={$t$},
xmajorgrids,
ylabel={$\tau_n$},
ymajorgrids,
ymode=log
]
\addplot [teal1293165, dashed, mark=*, mark size=2, mark options={solid}]
table {%
0.0025641 0.0025641
0.0051282 0.0025641
0.0076923 0.0025641
0.0102564 0.0025641
0.0128205 0.0025641
0.0153846 0.0025641
0.0179487 0.0025641
0.0205128 0.0025641
0.0230769 0.0025641
0.025641 0.0025641
0.0282051 0.0025641
0.0307692 0.0025641
0.0333333 0.0025641
0.0358974 0.0025641
0.0384615 0.0025641
0.0410256 0.0025641
0.0435897 0.0025641
0.0461538 0.0025641
0.0487179 0.0025641
0.051282 0.0025641
0.0538461 0.0025641
0.0564102 0.0025641
0.0589743 0.0025641
0.0615384 0.0025641
0.0641025 0.0025641
0.0666666 0.0025641
0.0692307 0.0025641
0.0717948 0.0025641
0.0743589 0.0025641
0.076923 0.0025641
0.0794871 0.0025641
0.0820512 0.0025641
0.0846153 0.0025641
0.0871794 0.0025641
0.0897435 0.0025641
0.0923076 0.0025641
0.0948717 0.0025641
0.0974358 0.0025641
0.0999999 0.0025641
0.102564 0.0025641
0.1051281 0.0025641
0.1076922 0.0025641
0.1102563 0.0025641
0.1128204 0.0025641
0.1153845 0.0025641
0.1179486 0.0025641
0.1205127 0.0025641
0.1230768 0.0025641
0.1256409 0.0025641
0.128205 0.0025641
0.1307691 0.0025641
0.1333332 0.0025641
0.1358973 0.0025641
0.1384614 0.0025641
0.1410255 0.0025641
0.1435896 0.0025641
0.1461537 0.0025641
0.1487178 0.0025641
0.1512819 0.0025641
0.153846 0.0025641
0.1564101 0.0025641
0.1589742 0.0025641
0.1615383 0.0025641
0.1641024 0.0025641
0.1666665 0.0025641
0.1692306 0.0025641
0.1717947 0.0025641
0.1743588 0.0025641
0.1769229 0.0025641
0.179487 0.0025641
0.1820511 0.0025641
0.1846152 0.0025641
0.1871793 0.0025641
0.1897434 0.0025641
0.1923075 0.0025641
0.1948716 0.0025641
0.1974357 0.0025641
0.1999998 0.0025641
0.2025639 0.0025641
0.205128 0.0025641
0.2076921 0.0025641
0.2102562 0.0025641
0.2128203 0.0025641
0.2153844 0.0025641
0.2179485 0.0025641
0.220512600000001 0.0025641
0.223076700000001 0.0025641
0.225640800000001 0.0025641
0.228204900000001 0.0025641
0.230769000000001 0.0025641
0.233333100000001 0.0025641
0.235897200000001 0.0025641
0.238461300000001 0.0025641
0.241025400000001 0.0025641
0.243589500000001 0.0025641
0.246153600000001 0.0025641
0.248717700000001 0.0025641
0.251281800000001 0.0025641
0.253845900000001 0.0025641
0.256410000000001 0.0025641
0.258974100000001 0.0025641
0.261538200000001 0.0025641
0.264102300000001 0.0025641
0.266666400000001 0.0025641
0.269230500000001 0.0025641
0.271794600000001 0.0025641
0.274358700000001 0.0025641
0.276922800000001 0.0025641
0.279486900000001 0.0025641
0.282051000000001 0.0025641
0.284615100000001 0.0025641
0.287179200000001 0.0025641
0.289743300000001 0.0025641
0.292307400000001 0.0025641
0.294871500000001 0.0025641
0.297435600000001 0.0025641
0.299999700000001 0.0025641
0.302563800000001 0.0025641
0.305127900000001 0.0025641
0.307692000000001 0.0025641
0.310256100000001 0.0025641
0.312820200000001 0.0025641
0.315384300000001 0.0025641
0.317948400000001 0.0025641
0.320512500000001 0.0025641
0.323076600000001 0.0025641
0.325640700000001 0.0025641
0.328204800000001 0.0025641
0.330768900000001 0.0025641
0.333333000000001 0.0025641
0.335897100000001 0.0025641
0.338461200000001 0.0025641
0.341025300000001 0.0025641
0.343589400000001 0.0025641
0.346153500000001 0.0025641
0.348717600000001 0.0025641
0.351281700000001 0.0025641
0.353845800000001 0.0025641
0.356409900000001 0.0025641
0.358974000000001 0.0025641
0.361538100000001 0.0025641
0.364102200000001 0.0025641
0.366666300000001 0.0025641
0.369230400000001 0.0025641
0.371794500000001 0.0025641
0.374358600000001 0.0025641
0.376922700000001 0.0025641
0.379486800000001 0.0025641
0.382050900000001 0.0025641
0.384615000000001 0.0025641
0.387179100000001 0.0025641
0.389743200000001 0.0025641
0.392307300000001 0.0025641
0.394871400000001 0.0025641
0.397435500000001 0.0025641
0.399999600000001 0.0025641
};
\addlegendentry{uni}
\addplot [limegreen018569, dashed, mark=square*, mark size=2, mark options={solid}]
table {%
6.529e-06 6.529e-06
1.9587e-05 1.3058e-05
4.5703e-05 2.6116e-05
8.2637e-05 3.6934e-05
0.000126958 4.4321e-05
0.000180143 5.3185e-05
0.000243965 6.3822e-05
0.000320551 7.6586e-05
0.000412454 9.1903e-05
0.000522734 0.00011028
0.000655074 0.00013234
0.000813884 0.00015881
0.001004454 0.00019057
0.001233134 0.00022868
0.001507554 0.00027442
0.001836864 0.00032931
0.002232034 0.00039517
0.002706234 0.0004742
0.003275274 0.00056904
0.003958124 0.00068285
0.004777544 0.00081942
0.005760844 0.0009833
0.006940844 0.00118
0.008356844 0.001416
0.010055944 0.0016991
0.012094944 0.002039
0.014541744 0.0024468
0.017477844 0.0029361
0.021001144 0.0035233
0.025088244 0.0040871
0.028457144 0.0033689
0.031096744 0.0026396
0.033545744 0.002449
0.036063744 0.002518
0.038606344 0.0025426
0.041259944 0.0026536
0.043935744 0.0026758
0.046707944 0.0027722
0.049888544 0.0031806
0.053286344 0.0033978
0.056544844 0.0032585
0.059610744 0.0030659
0.062512644 0.0029019
0.065462044 0.0029494
0.068568544 0.0031065
0.071663944 0.0030954
0.074592744 0.0029288
0.077445744 0.002853
0.080488144 0.0030424
0.083721544 0.0032334
0.086920744 0.0031992
0.090026544 0.0031058
0.092978144 0.0029516
0.095754944 0.0027768
0.098467044 0.0027121
0.101142744 0.0026757
0.103829244 0.0026865
0.106564744 0.0027355
0.109444344 0.0028796
0.112710144 0.0032658
0.116104244 0.0033941
0.119293144 0.0031889
0.122288244 0.0029951
0.125187144 0.0028989
0.128188044 0.0030009
0.131334344 0.0031463
0.134462544 0.0031282
0.137428344 0.0029658
0.140264844 0.0028365
0.143206444 0.0029416
0.146360744 0.0031543
0.149570544 0.0032098
0.152734844 0.0031643
0.155740644 0.0030058
0.158527644 0.002787
0.161220044 0.0026924
0.163878444 0.0026584
0.166553444 0.002675
0.169278544 0.0027251
0.172148344 0.0028698
0.175406444 0.0032581
0.178794644 0.0033882
0.181981444 0.0031868
0.184979544 0.0029981
0.187880844 0.0029013
0.190880744 0.0029999
0.194023344 0.0031426
0.197146044 0.0031227
0.200106244 0.0029602
0.202938844 0.0028326
0.205879444 0.0029406
0.209035844 0.0031564
0.212246944 0.0032111
0.215409344 0.0031624
0.218409444 0.0030001
0.221186644 0.0027772
0.223868944 0.0026823
0.226519644 0.0026507
0.229201544 0.0026819
0.231980444 0.0027789
0.234923544 0.0029431
0.238113044 0.0031895
0.241302644 0.0031896
0.244397744 0.0030951
0.247448744 0.003051
0.250449144 0.0030004
0.253474944 0.0030258
0.256607244 0.0031323
0.259764744 0.0031575
0.262738244 0.0029735
0.265574644 0.0028364
0.268514244 0.0029396
0.271653344 0.0031391
0.274858144 0.0032048
0.278024144 0.003166
0.281029544 0.0030054
0.283812644 0.0027831
0.286494244 0.0026816
0.289140244 0.002646
0.291816544 0.0026763
0.294587644 0.0027711
0.297524144 0.0029365
0.300710144 0.003186
0.303895144 0.003185
0.306987244 0.0030921
0.310039044 0.0030518
0.313040844 0.0030018
0.316066744 0.0030259
0.319197844 0.0031311
0.322353044 0.0031552
0.325323844 0.0029708
0.328158244 0.0028344
0.331097344 0.0029391
0.334237344 0.00314
0.337442844 0.0032055
0.340607944 0.0031651
0.343609844 0.0030019
0.346385944 0.0027761
0.349060444 0.0026745
0.351701444 0.002641
0.354373444 0.002672
0.357139244 0.0027658
0.360069744 0.0029305
0.363253544 0.0031838
0.366442444 0.0031889
0.369539944 0.0030975
0.372593544 0.0030536
0.375592544 0.002999
0.378614244 0.0030217
0.381743344 0.0031291
0.384899444 0.0031561
0.387870944 0.0029715
0.390704644 0.0028337
0.393642144 0.0029375
0.396781144 0.003139
0.399999944 0.0032188
};
\addlegendentry{adapt}
\end{axis}
\end{tikzpicture}
        }
        \caption{Time-step size}
    \end{subfigure}
    \caption{
        \label{FIG:bump-uniform-vs-adaptive}
		Example~\ref{SSEC:example5}:
        Comparison of the uniform and the adaptive algorithm.
        The adaptive method (Algorithm~\ref{ALG:Algo1}) has been run using cubic finite elements and BDF(4) with tolerances $\tol_{\mathrm{s}} = 10^{-4}$, $\tol_{\mathrm{t}} = 10^{-6}$, resulting in $N_t = 156$ time steps.
        For the uniform method we used $\tau = T/N_t$, $h = 1/42$.
    }
\end{figure}

%%%%%%%%%%%%%%%%%%%%%%%%%%%%%%%%%%%
\subsection{Example 5}\label{SSEC:example6}\ 
In this example we present numerical results for an extension of the model presented in Section~\ref{SSEC:LLG} in a domain $\Omega \subset\IR^3$.
Here the effective field is
\begin{align}\label{EQ:LLG-Heff-DMI}
   \bsH_{\text{eff}}(\bsm)
    =  C_{\mathrm{e}} \Delta\bsm
     + \bsH_{\mathrm{ext}}
     - C_{\mathrm{DMI}} \grad\times\bsm
     + \bsH_{\mathrm{str}}(\bsm),
\end{align}
where the third term is the Dzyaloshinskii--Moriya interaction (DMI) and $\bsH_{\mathrm{str}}(\bsm)$ is the magnetic stray field \cite{HPPRSS:2019}.
The DMI interaction will change the boundary condition to $C_{\mathrm{e}}\ptl_{\bsn}\bsm =-1/2\,C_{\mathrm{DMI}}\,\bsn\times\bsm$, while to get $\bsH_{\mathrm{str}}(\bsm)$ one has to solve a potential equation in $\IR^3$ with $\bsm|_{\Omega}$ as a source term.
This problem is solved using a boundary integral approach, for details we refer to \cite[Alg.~12]{PRS:2018}.

We consider the evolution of a magnetic Skyrmion within a disk of diameter $160 \unit{nm}$ and height $10 \unit{nm}$, an example adopted from \cite{DDPR:2022}.
We take the reference length $L_{\mathrm{ref}} = 1 \unit{nm}$ and perform the simulation over a total time of $T = 1 \unit{ns}$.
With realistic parameters adopted from \cite{DDPR:2022} we find the dimensionless parameter
$\alpha = 0.28$, $C_{\mathrm{e}} = 94.77$, $C_{\mathrm{DMI}} = 17.05$.

For the numerical discretization, we apply the fully adaptive Algorithm~\ref{ALG:Algo1} with quadratic finite elements and BDF($2$). 
We take the mesh refinement parameters $\theta_{\mathrm{r}} = 0.85$, $\theta_{\mathrm{c}} = 0.9$ and tolerances $\tol_{\mathrm{s}} = 1.5$ and $\tol_{\mathrm{t}} =10^{-1}$.

%\commsk{bzw. relative Toleranzen $\tilde \tol_{\mathrm{s}} = 3 \cdot 10^{-3}$ and $\tilde \tol_{\mathrm{t}} = 2 \cdot 10^{-4}$.}

\begin{center}
\begin{figure}[t!]
	\centering
	\begin{subfigure}{0.4\textwidth}
		\centering
		\includegraphics[width=\linewidth]{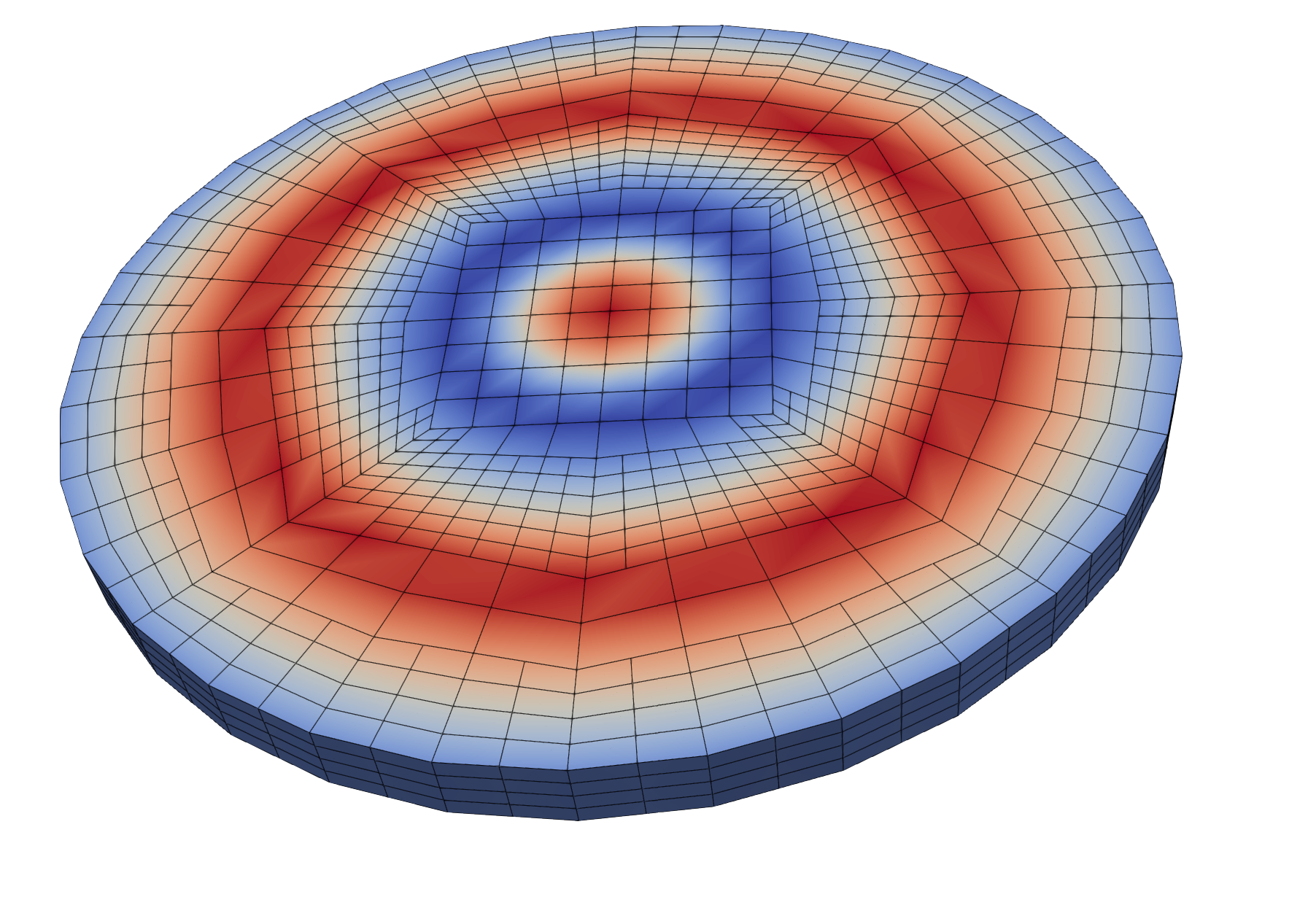}
		% \caption{...}
	\end{subfigure}
	\begin{subfigure}{0.4\textwidth}
		\centering
		\includegraphics[width=\linewidth]{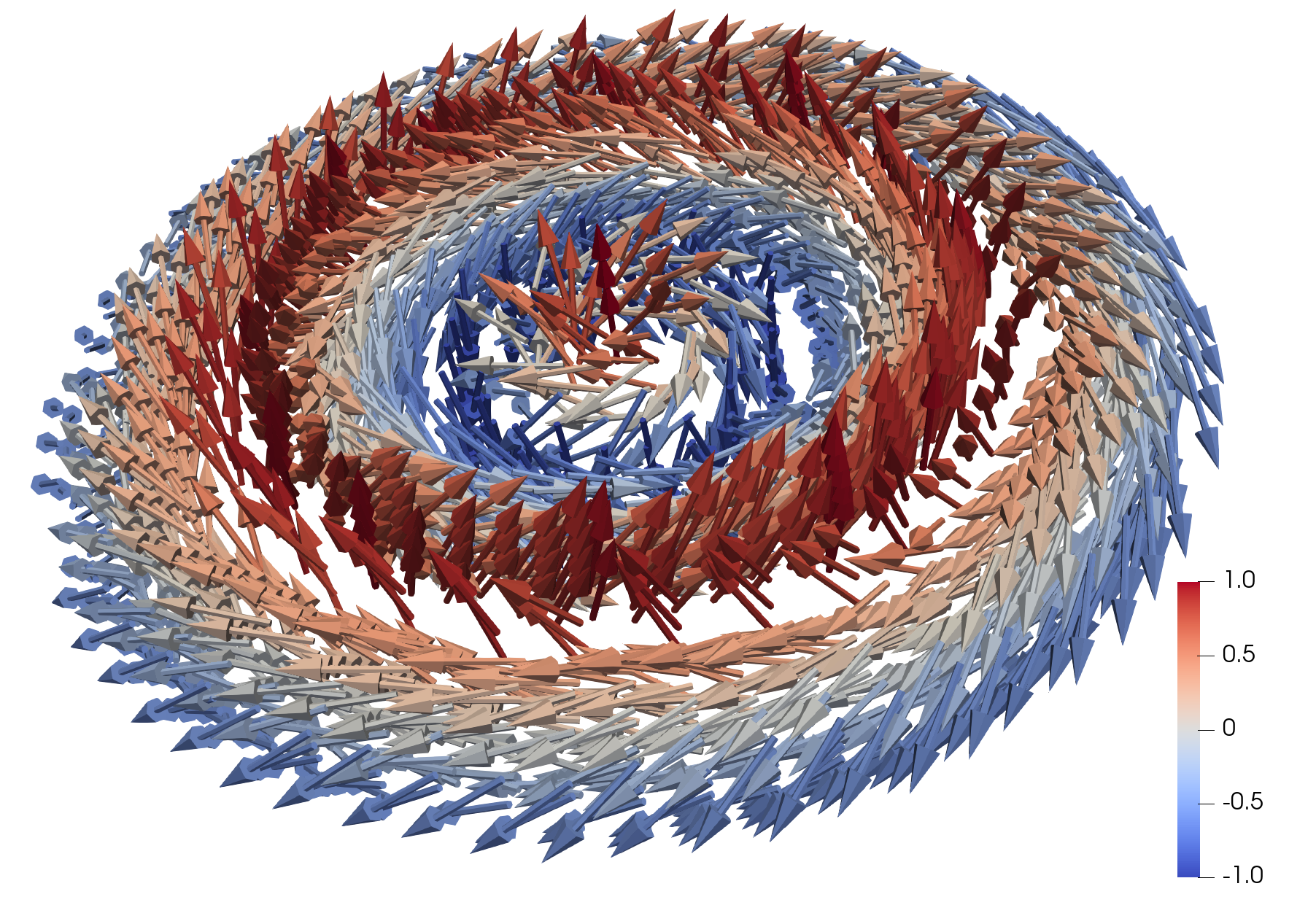}
		% \caption{...}
	\end{subfigure}
	\caption{
        \label{FIG:skyrmion-mesh}
        Example~\ref{SSEC:example6}:
		Skyrmion at $t=T$.
		The background color represents the $z$-component of the magnetization.
        % $h_{\min} \approx 5.5$, $h_{\max}\approx 16$, \commwd{???}
	}
\end{figure}
\end{center}

\begin{figure}[htb!]
    \centering
    \begin{subfigure}{0.48\textwidth}
        \centering
        \resizebox{\textwidth}{!}{%
            \input{new_plots/Skyrmion/avg_z}
        }
        \caption{Average $z$-component}
    \end{subfigure}
    \hfill
    \begin{subfigure}{0.48\textwidth}
        \centering
        \resizebox{\textwidth}{!}{%
            \input{new_plots/Skyrmion/norm_err}
        }
        \caption{Normalization error}
    \end{subfigure}
    \\
    \begin{subfigure}{0.48\textwidth}
        \centering
        \resizebox{\textwidth}{!}{%
            \input{new_plots/Skyrmion/spacedim}
        }
        \caption{Degrees of freedom}
    \end{subfigure}
    \hfill
    \begin{subfigure}{0.48\textwidth}
        \centering
        \resizebox{\textwidth}{!}{%
            \input{new_plots/Skyrmion/timesteps}
        }
        \caption{Time-step size}
    \end{subfigure}
    \caption{
        \label{FIG:Skyrmion-adaptive}
		Example~\ref{SSEC:example6}:
        The adaptive method (Algorithm~\ref{ALG:Algo1}) has been run using quadratic finite elements and BDF(2).
        %with tolerances $\tol_{\mathrm{s}} = XXX$, $\tol_{\mathrm{t}} = XXX$.
    }
\end{figure}

%%%%%%%%%%%%%%%%%%%%%%%%%%%%%%%%%%%%%%%%%%%%%%%%%%%%%%%%%%%%%%%%%%%%%%%
\section*{Conclusion}

We propose a time- and space-adaptive algorithm for the numerical approximation of the Landau--Lifshitz--Gilbert equation.
Under certain regularity assumptions on the exact solutions, we show that the numerical approximation satisfies an energy bound similarly to that of the exact solution.
Numerical experiments demonstrate the advantages of adaptive algorithms over uniform approaches, such as increased convergence speed (Example~\ref{SSEC:example1}) or increased stability (Example~\ref{SSEC:example3}).
Many interesting theoretical questions remain open for future research. Is the error estimator an upper bound for the error? Does the adaptive algorithm converge towards the exact solution and, if yes, with optimal rate (see optimal rates for adaptive time-stepping in~\cite{F:2022})? Does the energy bound (Theorem~\ref{THM:stability}) hold under reduced regularity assumptions?

%%%%%%%%%%%%%%%%%%%%%%%%%%%%%%%%%%%%%%%%%%%%%%%%%%%%%%%%%%%%%%%%%%%%%%%
\section{Appendix: $\bsG$-matrix for BDF($2$) with variable time-steps}

Recall the variable step-size BDF($2$) coefficients
\begin{align}\label{EQ:delta-def}
   \delta_0(\kappa) = \frac{1+2\kappa}{1+\kappa},\qquad
   \delta_1(\kappa) = -(1+\kappa),\qquad
   \delta_2(\kappa) = \frac{\kappa^2}{1+\kappa},
\end{align}
and 
\begin{align}\label{EQ:g-def}
    g(\kappa) := \frac{3\kappa^2-2\kappa-1}{(1+\kappa)^2},
\end{align}
where $g(\kappa)$ is the threshold such that $\Real\,\delta(\kappa;s)/(1-\eta s)>0$ for $|s|<1$ holds if $\eta>g(\kappa)$, with $\delta(\kappa; s) := \delta_0(\kappa) + \delta_1(\kappa)s + \delta_2(\kappa)s^2$,~cf.~Step~1 of the proof of Theorem~\ref{THM:stability}.
Note that $g(\kappa) < 1$ if and only if $\kappa < 1 + \sqrt{2}$, so the interval $(g(\kappa), 1)$ used below is nonempty for $0 < \kappa < 1 + \sqrt2$.
Dahlquist's theorem \cite[Lem.~A.1]{AFKL:2021} then guarantees a symmetric positive definite matrix $\bsG$ with an identity of the form \eqref{EQ:Gidentity} below, but it does not provide any information on how $\bsG$ depends on $\kappa$, which we provide in Lemma~\ref{LEM:Gconstruction} and Proposition~\ref{PROP:Lipschitz}.

\begin{Lemma}[Monotonicity of $g$]\label{LEM:gmono}
$g$ given by \eqref{EQ:g-def} is strictly increasing on $(0,\infty)$, with $g'(\kappa)=8\kappa/(1+\kappa)^3>0$.
\end{Lemma}
\begin{proof}
Direct differentiation gives
\begin{align*}
   g'(\kappa) 
   %= \frac{(6\kappa-2)(1+\kappa)^2-(3\kappa^2-2\kappa-1)\cdot2(1+\kappa)}{(1+\kappa)^4}
   = \frac{(6\kappa-2)(1+\kappa)-2(3\kappa^2-2\kappa-1)}{(1+\kappa)^3}
   = \frac{8\kappa}{(1+\kappa)^3},
\end{align*}
which is positive for $\kappa>0$.
\end{proof}

\begin{Lemma}[Positivity of $\sigma$]\label{LEM:sigmapos}
Let $0 < \kappa < 1 + \sqrt2$ and $\eta\in(g(\kappa),1)$. Set
\begin{align}
  A(\kappa,\eta)
  &:= 2\kappa^2+(1+\kappa)^2\big(\eta-g(\kappa)\big), \label{EQ:Adef}\\
  \sigma(\kappa,\eta)
  &:= \frac{A(\kappa,\eta)-(1+\kappa)^2
     \sqrt{(\eta-g(\kappa))(1+\eta)}}{4(1+\kappa)}. \label{EQ:sigmadef}
\end{align}
Then $\sigma(\kappa,\cdot)$ is strictly decreasing on $(g(\kappa),1]$, and $\sigma(\kappa,\eta)>0$.
\end{Lemma}
\begin{proof}
Fix $\kappa>0$ and abbreviate $a:=\eta-g(\kappa)$ and $b:=1+\eta$, so that the square root in \eqref{EQ:sigmadef} is $\sqrt{ab}$.
Since
\begin{align*}
   b-a = 1+g(\kappa) = \frac{4\kappa^2}{(1+\kappa)^2} > 0,
\end{align*}
we have $0<a<b$ for every $\eta\in(g(\kappa),1]$.
In particular $\sqrt{ab}>0$ and thus $\sigma(\kappa,\cdot)$ is well defined and differentiable on $(g(\kappa),1]$.

\emph{Step 1: $\sigma$ is strictly decreasing in $\eta$}.
We have $\partial_\eta A=(1+\kappa)^2$ and $\partial_\eta(ab)=b+a$, hence
$\partial_\eta\sqrt{ab}=(a+b)/(2\sqrt{ab})$.
Differentiating \eqref{EQ:sigmadef} therefore gives
\begin{align*}
   \frac{\partial\sigma}{\partial\eta}
   = \frac{1+\kappa}{4}\big(1-\frac{a+b}{2\sqrt{ab}}\big)
   = -\frac{1+\kappa}{8\sqrt{ab}}\big(\sqrt b-\sqrt a\big)^2 < 0,
\end{align*}
where we used $2\sqrt{ab}-a-b=-(\sqrt b-\sqrt a)^2$ and $a\neq b$ since $b > a$. 
%\commwd{da $b>a$?}.

\emph{Step 2: $\sigma(\kappa,1)>0$}.
For $\eta=1$ we have $(1+\kappa)^2\big(1-g(\kappa)\big)=2(1+2\kappa-\kappa^2)$, which is positive because $g(\kappa)<1$ since $\kappa < 1 + \sqrt2$.
Hence $A(\kappa,1)=2\kappa^2+2(1+2\kappa-\kappa^2)=2(2\kappa+1)$ and
$(1+\kappa)^2\sqrt{2(1-g(\kappa))}=2(1+\kappa)\sqrt{1+2\kappa-\kappa^2}$, so that
\begin{align*}
   \sigma(\kappa,1)
   = \frac{2\kappa+1-(1+\kappa)\sqrt{1+2\kappa-\kappa^2}}{2(1+\kappa)}.
\end{align*}
Both terms in the numerator are positive, and
\begin{align*}
   (1+\kappa)^2(1+2\kappa-\kappa^2)
   = (1+2\kappa+\kappa^2)(1+2\kappa-\kappa^2)
   = (1+2\kappa)^2-\kappa^4<(2\kappa+1)^2.
\end{align*}
Taking square roots gives $(1+\kappa)\sqrt{1+2\kappa-\kappa^2}<2\kappa+1$, i.e.\ $\sigma(\kappa,1)>0$.

In summary, for $\eta\in(g(\kappa),1)$, Steps~1 and~2 give $\sigma(\kappa,\eta)>\sigma(\kappa,1)>0$.
\end{proof}

\begin{Lemma}[A key identity for $\sigma$]\label{LEM:sigmaid}
Let $0 < \kappa < 1 + \sqrt2$ and $\eta\in(g(\kappa),1)$, and let $A(\kappa,\eta)$ and $\sigma(\kappa,\eta)$ be given by \eqref{EQ:Adef}--\eqref{EQ:sigmadef}.
Then
\begin{align}\label{EQ:sigmaid}
   2\Big(\sigma(\kappa,\eta)+\frac{\delta_2(\kappa)^2}{4\sigma(\kappa,\eta)}\Big)
   = \delta_0(\kappa)-\delta_2(\kappa)+(1+\kappa)\eta.
\end{align}
\end{Lemma}
\begin{proof}
For convenience, we abbreviate $A=A(\kappa,\eta)$, $\sigma=\sigma(\kappa,\eta)$ and $\delta_i = \delta_i(\kappa)$.
Note that $\sigma>0$ by Lemma~\ref{LEM:sigmapos}, so the left-hand side is well defined.
With $a=\eta-g(\kappa)$ and $b=1+\eta$ as in the proof of Lemma~\ref{LEM:sigmapos}, we have $A-2\kappa^2=(1+\kappa)^2a$ and, since $b-a=4\kappa^2/(1+\kappa)^2$, also $A+2\kappa^2=(1+\kappa)^2b$.
Multiplying both equations gives $A^2-4\kappa^4=(1+\kappa)^4ab$.
Now set $\hat\sigma:=\big(A+(1+\kappa)^2\sqrt{ab}\big)/(4(1+\kappa))$.
Then
\begin{align*}
   \sigma\hat\sigma = \frac{A^2-(1+\kappa)^4ab}{16(1+\kappa)^2}
   = \frac{A^2-(A+2\kappa^2)(A-2\kappa^2)}{16(1+\kappa)^2}
   = \frac{\kappa^4}{4(1+\kappa)^2}
   = \frac{\delta_2^2}{4},
\end{align*}
so $\delta_2^2/(4\sigma)=\hat\sigma$ and therefore $2\big(\sigma+\delta_2^2/(4\sigma)\big)=2(\sigma+\hat\sigma)=A/(1+\kappa)$.
Finally, $A=(1+\kappa)^2\eta+1+2\kappa-\kappa^2$ and $\delta_0-\delta_2=(1+2\kappa-\kappa^2)/(1+\kappa)$ give $A/(1+\kappa)=\delta_0-\delta_2+(1+\kappa)\eta$.
\end{proof}

Following the approach in \cite[(1.7)]{AkrivisKatsoprinakis:2016} and \cite[V.6,~Ex.~6.5]{HairerWanner:2010}, we determine explicitly the $\bsG$-matrix for variable step-size BDF($2$) in the following lemma, which is ensured to exist by Dahlquist's theorem \cite[Lem.~A.1]{AFKL:2021} if the step ratio satisfies $\kappa < 1 + \sqrt{2}$.

\begin{Lemma}[Explicit $\bsG$-matrix]\label{LEM:Gconstruction}
Let $V$ be a real Hilbert space with inner product $(\cdot,\cdot)=(\cdot,\cdot)_V$ and norm $|w|:=(w,w)^{1/2}$.
Let $0 < \kappa < 1 + \sqrt2$ and $\eta\in(g(\kappa),1)$, and let $A(\kappa,\eta)$ and $\sigma(\kappa,\eta)$ be given by \eqref{EQ:Adef}--\eqref{EQ:sigmadef}, and set
\begin{equation}\label{EQ:lammudef}
  \mu(\kappa,\eta) := -\frac{\delta_2(\kappa)}{2\sigma(\kappa,\eta)}, \qquad
  \lambda(\kappa,\eta) := 1-\mu(\kappa,\eta)
\end{equation}
which are well defined since $\sigma(\kappa,\eta)>0$ by Lemma~\ref{LEM:sigmapos}. Then the matrix
\begin{align*}
   \bsG(\kappa,\eta) := 
   \begin{pmatrix}
      g_{11}(\kappa,\eta) & g_{12}(\kappa,\eta)\\
      g_{12}(\kappa,\eta) & g_{22}(\kappa,\eta)
   \end{pmatrix}
\end{align*}
with entries
\begin{align*}
   g_{11}(\kappa,\eta) = \frac{\delta_2(\kappa)^2}{4\sigma(\kappa,\eta)},\quad
   g_{22}(\kappa,\eta) = \delta_0(\kappa)-\sigma(\kappa,\eta),\quad
   g_{12}(\kappa,\eta) = \frac{\delta_2(\kappa)}{2}(\eta-1)-g_{11}(\kappa,\eta)
\end{align*}
is symmetric positive definite and satisfies, for all $v_0,v_1,v_2\in V$,
\begin{multline}
   \big(\delta_0(\kappa)v_2+\delta_1(\kappa)v_1+\delta_2(\kappa)v_0\,,
      \,v_2-\eta v_1\big)
   = \Bigl[g_{11}(v_1,v_1)+2g_{12}(v_1,v_2)+g_{22}(v_2,v_2)\Bigr] \\
     -\Bigl[g_{11}(v_0,v_0)+2g_{12}(v_0,v_1)+g_{22}(v_1,v_1)\Bigr]
     +\sigma(\kappa,\eta)\,\bigl|v_2-\lambda(\kappa,\eta)v_1-\mu(\kappa,\eta)v_0\bigr|^2. \label{EQ:Gidentity}
\end{multline}
\end{Lemma}
\begin{proof}
In the following, we suppress the arguments $\kappa$ and $\eta$.
Expanding both sides of \eqref{EQ:Gidentity} and matching the coefficients yields the six equations
\begin{alignat*}{3}
  \delta_0-g_{22} &= \sigma, &\qquad
  \delta_1-\delta_0\eta-2g_{12} &= -2\sigma\lambda, &\qquad
  -\delta_1\eta-g_{11}+g_{22} &= \sigma\lambda^2,\\
  \delta_2 &= -2\sigma\mu, &\qquad
  -\delta_2\eta+2g_{12} &= 2\sigma\lambda\mu, &\qquad
  g_{11} &= \sigma\mu^2.
\end{alignat*}
Throughout, we use $\delta_0+\delta_2=-\delta_1=1+\kappa$, $\lambda=1-\mu$, and \eqref{EQ:sigmaid}.
The first equation holds by the definition of $g_{22}$. 
The fourth equation holds by the definition of $\mu$.
The fifth equation holds because, by the definition of $g_{12}$, its left-hand side equals $-\delta_2-2g_{11}$, and by the fourth and sixth equations its right-hand side equals $2\sigma\mu-2\sigma\mu^2=-\delta_2-2g_{11}$.
The sixth equation holds by the definitions of $g_{11}$ and $\mu$, since $\sigma\mu^2=\delta_2^2/(4\sigma)$.
The second equation is equivalent to \eqref{EQ:sigmaid}: by the definition of $g_{12}$, its left-hand side equals $\delta_1+\delta_2-(\delta_0+\delta_2)\eta+2g_{11}=-\delta_0-(1+\kappa)\eta+2g_{11}$, and by the fourth equation its right-hand side equals $-2\sigma+2\sigma\mu=-2\sigma-\delta_2$.
The third equation is also equivalent to \eqref{EQ:sigmaid}: by the fourth and sixth equations, $\sigma\lambda^2=\sigma-2\sigma\mu+\sigma\mu^2=\sigma+\delta_2+g_{11}$, so with $g_{22}=\delta_0-\sigma$ the equation reads $(1+\kappa)\eta-g_{11}+\delta_0-\sigma=\sigma+\delta_2+g_{11}$.
This concludes \eqref{EQ:Gidentity}.

We continue with the positive definiteness of $\bsG$. 
The matrix $\bsG$ is symmetric and $g_{11}=\delta_2^2/(4\sigma)>0$ by Lemma~\ref{LEM:sigmapos}, so it suffices to show $\det(\bsG)>0$. 
By \eqref{EQ:sigmaid} and \eqref{EQ:delta-def}, we have
\begin{align*}
   g_{22}-g_{11} = \delta_0-(\sigma+g_{11})
   = \tfrac12\big(\delta_0+\delta_2-(1+\kappa)\eta\big)
   = \tfrac12(1+\kappa)(1-\eta)
\end{align*}
and thus get
\begin{align*}
    \frac{\delta_2}{2}(\eta-1) = \frac{\delta_2}{1+\kappa}(g_{11}-g_{22}).
\end{align*}
Inserting this into the definition of $g_{12}$ and using $\delta_2-(1+\kappa)=-\delta_0$ gives
\begin{align*}
   g_{12} = -\frac{\delta_0g_{11}+\delta_2g_{22}}{1+\kappa}.
\end{align*}
With $(1+\kappa)^2=(\delta_0+\delta_2)^2$ we obtain
\begin{align*}
   (1+\kappa)^2\det(\bsG)
   = (\delta_0+\delta_2)^2g_{11}g_{22}-(\delta_0g_{11}+\delta_2g_{22})^2
   = (g_{22}-g_{11})\big(\delta_0^2g_{11}-\delta_2^2g_{22}\big).
\end{align*}
The first factor equals $\tfrac12(1+\kappa)(1-\eta)$, and inserting $g_{11}=\delta_2^2/(4\sigma)$ and $g_{22}=\delta_0-\sigma$ into the second factor gives
$\delta_0^2g_{11}-\delta_2^2g_{22}=\delta_2^2(\delta_0-2\sigma)^2/(4\sigma)$. Therefore
\begin{align*}
   \det(\bsG) = \frac{(1-\eta)\,\delta_2^2\,(\delta_0-2\sigma)^2}%
      {8(1+\kappa)\,\sigma},
\end{align*}
which is positive for $\eta<1$ since $2\sigma\neq\delta_0$ (precisely $2 \sigma < \delta_2 < \delta_0$). 
\end{proof}

\begin{Proposition}[Lipschitz continuity of $\bsG$]\label{PROP:Lipschitz}
Let $0<\kappa_{\min} \le \kappa_{\max} < 1 + \sqrt2$, set $K:=[\kappa_{\min},\kappa_{\max}]$, and fix $\eta_\star\in\big(g(\kappa_{\max}),1\big)$. 
Then, $\kappa\mapsto\bsG(\kappa,\eta_\star)$ is continuously differentiable on $K$ and it holds
\begin{align}\label{EQ:GLipschitz}
   \bigl\|\bsG(\kappa,\eta_\star)-\bsG(\kappa',\eta_\star)\bigr\|_2
   \le L\,|\kappa-\kappa'| \qquad\text{for all }\kappa,\kappa'\in K ,
\end{align}
where $L:=\max_{\kappa\in K}\bigl\|\partial_\kappa\bsG(\kappa,\eta_\star)\bigr\|_2<\infty$, and $\norm{\cdot}{2}$ denotes the spectral norm.
Further, the smallest and largest eigenvalues $\gamma_-(\kappa)$ and $\gamma_+(\kappa)$ of $\bsG(\kappa,\eta_\star)$ are bounded uniformly on $K$ by 
\begin{align*}
   \gamma_-^{\min} := \min_{\kappa\in K}\gamma_-(\kappa)>0,\qquad
   \gamma_+^{\max} := \max_{\kappa\in K}\gamma_+(\kappa)<\infty.
\end{align*}
\end{Proposition}
\begin{proof}
Since $g$ is continuous and strictly increasing by Lemma~\ref{LEM:gmono}, $J:=\{\kappa>0:\ g(\kappa)<\eta_\star\}$ is an open interval of the form $(0,\kappa_\star)$, and $K\subset J$ since $g(\kappa_{\max})<\eta_\star$.
Since $\kappa\mapsto\sigma(\kappa,\eta_\star)$ is continuously differentiable on $J$, and $\sigma(\kappa,\eta_\star)>0$ on $J$ by Lemma~\ref{LEM:sigmapos}, the entries $g_{11}$, $g_{22}$ and $g_{12}$ are continuously differentiable on $J$ as well.
Hence $\kappa\mapsto\partial_\kappa\bsG(\kappa,\eta_\star)$ is continuous on $K\subset J$ with $L=\max_{\kappa\in K}\bigl\|\partial_\kappa\bsG(\kappa,\eta_\star)\bigr\|_2<\infty$.

Now, let $\kappa,\kappa'\in K$. Then the fundamental theorem of calculus gives
\[
   \bsG(\kappa,\eta_\star)-\bsG(\kappa',\eta_\star)
   = \int_{\kappa'}^{\kappa}\partial_s\bsG(s,\eta_\star)\,\mathrm{d}s.
\]
Bounding the norm of the integral by $L$ yields \eqref{EQ:GLipschitz}.

Since $\kappa \mapsto \bsG(\kappa,\eta_\star)$ is continuous on $K$, the eigenvalues $\gamma_-(\kappa)$ and $\gamma_+(\kappa)$ are continuous on $K$ as well.
Since $K$ is compact, the eigenvalues attain their minimum and maximum on $K$.
Finally, since $\bsG(\kappa,\eta_\star)$ is positive definite by Lemma~\ref{LEM:Gconstruction}, we have $\gamma_-(\kappa)>0$ for every $\kappa\in K$, and hence $\infty > \gamma_+^{\max} \ge \gamma_-^{\min}>0$.
\end{proof}

\begin{Corollary}[Dependence of the $\bsG$-norm on $\kappa$]\label{COR:commutator}
Let $V$ be a real Hilbert space with inner product $(\cdot,\cdot)$ and norm $|\cdot|$.
For a symmetric matrix $\bsG \in \IR^{2\times2}$ and $W=(w_1,w_2)\in V^2$ set
\begin{align*}
   \norm{W}{\bsG}^2:=\sum_{i,\ell=1}^2 g_{i\ell}\,(w_i,w_\ell).
\end{align*}
Then, with $K$, $\eta_\star$ and $L$ as in Proposition~\ref{PROP:Lipschitz}, it holds
\begin{align}
   \Bigl|\norm{W}{\bsG(\kappa,\eta_\star)}^2
      -\norm{W}{\bsG(\kappa',\eta_\star)}^2\Bigr|
   \le L\,|\kappa-\kappa'|\,\big(|w_1|^2+|w_2|^2\big)
   \label{EQ:Gnormdiff}
\end{align}
for all $\kappa,\kappa'\in K,\ W\in V^2$.
\end{Corollary}

%%%%%%%%%%%%%%%%%%%%%%%%%%%%%%%%%%%%%%%%%%%%%%%%%%%%%%%%%%%%%%%%%%%%%%%
\subsection*{Acknowledgement}
Willy Dörfler, Michael Feischl, and Stefan Karch gratefully acknowledge the support of the Deutsche Forschungsgemeinschaft (DFG) within the SFB 1173 'Wave Phenomena' (Project ID 258734477).
Michael Feischl additionally acknowledges support of the FWF (Austrian Science Fund) within SFB 65 "Taming Complexity in Partial Differential Systems" as well as project I6667-N.
Funding was also received from the European Research Council (ERC) under the Horizon 2020 research and innovation program (Grant agreement No.~101125225).
We thank J. Bohn for his contribution on a previous version of the paper.

%%%%%%%%%%%%%%%%%%%%%%%%%%%%%%%%%%%%%%%%%%%%%%%%%%%%%%%%%%%%%%%%%%%%%%%
\bibliographystyle{plain}
\bibliography{project_B8}

@article{guo2025recovery,
  title={Recovery techniques for finite element methods},
  author={Guo, Hailong and Zhang, Zhimin},
  journal={Advances in Applied Mechanics},
  volume={60},
  pages={399--463},
  year={2025},
  publisher={Elsevier}
}

@article{bartels2002each,
  title={{Each averaging technique yields reliable a posteriori error control in FEM on unstructured grids. Part II: Higher order FEM}},
  author={Bartels, S{\"o}ren and Carstensen, Carsten},
  journal={Mathematics of computation},
  volume={71},
  number={239},
  pages={971--994},
  year={2002}
}

@article {zzest,
    AUTHOR = {Zienkiewicz, Olgierd C. and Zhu, J. Z.},
     TITLE = {A simple error estimator and adaptive procedure for practical
              engineering analysis},
   JOURNAL = {Internat. J. Numer. Methods Engrg.},
  FJOURNAL = {International Journal for Numerical Methods in Engineering},
    VOLUME = {24},
      YEAR = {1987},
    NUMBER = {2},
     PAGES = {337--357},
      ISSN = {0029-5981,1097-0207},
   MRCLASS = {73K25},
  MRNUMBER = {875306},
       DOI = {10.1002/nme.1620240206},
       URL = {https://doi.org/10.1002/nme.1620240206},
}

@article {NevanlinnaOdeh:1981,
    AUTHOR = {Nevanlinna, Olavi and Odeh, Farouk},
     TITLE = {Multiplier techniques for linear multistep methods},
   JOURNAL = {Numer. Funct. Anal. Optim.},
  FJOURNAL = {Numerical Functional Analysis and Optimization. An
              International Journal},
    VOLUME = {3},
      YEAR = {1981},
    NUMBER = {4},
     PAGES = {377--423},
      ISSN = {0163-0563},
   MRCLASS = {65L05},
  MRNUMBER = {636736},
MRREVIEWER = {Peter Alfeld},
       DOI = {10.1080/01630568108816097},
       URL = {https://doi.org/10.1080/01630568108816097},
}

@article {ChengShen:2023,
    AUTHOR = {Cheng, Qing and Shen, Jie},
     TITLE = {Length preserving numerical schemes for {L}andau-{L}ifshitz
              equation based on {L}agrange multiplier approaches},
   JOURNAL = {SIAM J. Sci. Comput.},
  FJOURNAL = {SIAM Journal on Scientific Computing},
    VOLUME = {45},
      YEAR = {2023},
    NUMBER = {2},
     PAGES = {A530--A553},
      ISSN = {1064-8275},
   MRCLASS = {65M70 (35K61 65N12 65N22)},
  MRNUMBER = {4579735},
MRREVIEWER = {Xiang Wang},
       DOI = {10.1137/22M1501143},
       URL = {https://doi.org/10.1137/22M1501143},
}

@article {LiaoZhang:2021,
    AUTHOR = {Liao, Hong-lin and Zhang, Zhimin},
     TITLE = {Analysis of adaptive {BDF}2 scheme for diffusion equations},
   JOURNAL = {Math. Comp.},
  FJOURNAL = {Mathematics of Computation},
    VOLUME = {90},
      YEAR = {2021},
    NUMBER = {329},
     PAGES = {1207--1226},
      ISSN = {0025-5718,1088-6842},
   MRCLASS = {65M06 (65M12 65M50)},
  MRNUMBER = {4232222},
MRREVIEWER = {De\ Kang\ Mao},
       DOI = {10.1090/mcom/3585},
       URL = {https://doi.org/10.1090/mcom/3585},
}

@article {AFKL:2021,
    AUTHOR = {Akrivis, Georgios and Feischl, Michael and Kov\'{a}cs, Bal\'{a}zs
              and Lubich, Christian},
     TITLE = {Higher-order linearly implicit full discretization of the
              {L}andau--{L}ifshitz--{G}ilbert equation},
   JOURNAL = {Math. Comp.},
  FJOURNAL = {Mathematics of Computation},
    VOLUME = {90},
      YEAR = {2021},
    NUMBER = {329},
     PAGES = {995--1038},
      ISSN = {0025-5718},
   MRCLASS = {65M60 (35Q60 65L06 65M12 65M15)},
  MRNUMBER = {4232216},
MRREVIEWER = {Hamdullah Y\"{u}cel},
       DOI = {10.1090/mcom/3597},
       URL = {https://doi.org/10.1090/mcom/3597},
}

@article {BKW:2024,
    AUTHOR = {Bartels, S\"{o}ren and Kov\'{a}cs, Bal\'{a}zs and Wang, Zhangxian},
     TITLE = {Error analysis for the numerical approximation of the harmonic
              map heat flow with nodal constraints},
   JOURNAL = {IMA J. Numer. Anal.},
  FJOURNAL = {IMA Journal of Numerical Analysis},
    VOLUME = {44},
      YEAR = {2024},
    NUMBER = {2},
     PAGES = {633--653},
      ISSN = {0272-4979},
   MRCLASS = {65M12 (65M60 80A19)},
  MRNUMBER = {4727106},
       DOI = {10.1093/imanum/drad037},
       URL = {https://doi.org/10.1093/imanum/drad037},
}

@article{fert,
Author = {Fert, Albert and Cros, Vincent and Sampaio, Joao},
Year = {2013},
Title = {Skyrmions on the track},
Journal = {Nat. Nanotechnol.},
Pages ={152--156},
Volume ={8}
}

@article {harmonicmap,
    AUTHOR = {Coron, Jean M.},
     TITLE = {Nonuniqueness for the heat flow of harmonic maps},
   JOURNAL = {Ann. Inst. H. Poincar\'{e} C Anal. Non Lin\'{e}aire},
  FJOURNAL = {Annales de l'Institut Henri Poincar\'{e} C. Analyse Non Lin\'{e}aire},
    VOLUME = {7},
      YEAR = {1990},
    NUMBER = {4},
     PAGES = {335--344},
      ISSN = {0294-1449},
   MRCLASS = {58E20 (58G11)},
  MRNUMBER = {1067779},
MRREVIEWER = {Giorgio Valli},
       DOI = {10.1016/S0294-1449(16)30295-5},
       URL = {https://doi.org/10.1016/S0294-1449(16)30295-5},
}

@article {weakstrong,
    AUTHOR = {Di Fratta, Giovanni and Innerberger, Michael and Praetorius,
              Dirk},
     TITLE = {Weak-strong uniqueness for the {L}andau-{L}ifshitz-{G}ilbert
              equation in micromagnetics},
   JOURNAL = {Nonlinear Anal. Real World Appl.},
  FJOURNAL = {Nonlinear Analysis. Real World Applications. An International
              Multidisciplinary Journal},
    VOLUME = {55},
      YEAR = {2020},
     PAGES = {103122, 13},
      ISSN = {1468-1218},
   MRCLASS = {35Q60 (35A01 35Q82 82D40)},
  MRNUMBER = {4077401},
MRREVIEWER = {Gaetano Siciliano},
       DOI = {10.1016/j.nonrwa.2020.103122},
       URL = {https://doi.org/10.1016/j.nonrwa.2020.103122},
}

@article{adaptiveLLG1,
title = {Adaptive geometric integration applied to a {3D} micromagnetic solver},
journal = {Journal of Magnetism and Magnetic Materials},
volume = {518},
pages = {167409},
year = {2021},
issn = {0304-8853},
doi = {https://doi.org/10.1016/j.jmmm.2020.167409},
url = {https://www.sciencedirect.com/science/article/pii/S0304885320323763},
author = {Riccardo Ferrero and Alessandra Manzin}
}

@article {AkrivisKatsoprinakis:2016,
    AUTHOR = {Akrivis, Georgios and Katsoprinakis, Emmanuil},
     TITLE = {Backward difference formulae: new multipliers and stability properties for parabolic equations},
   JOURNAL = {Math. Comp.},
  FJOURNAL = {Mathematics of Computation},
    VOLUME = {85},
      YEAR = {2016},
    NUMBER = {301},
     PAGES = {2195–2216},
      ISSN = {0025-5718},
   MRCLASS = {65J08 (65L06)},
  MRNUMBER = {MR3511279},
MRREVIEWER = {},
       DOI = {10.1090/mcom3055},
       URL = {https://doi.org/10.1090/mcom3055},
}

@article {AlougesSoyeur:1992,
    AUTHOR = {Alouges, Fran\c{c}ois and Soyeur, Alain},
     TITLE = {On global weak solutions for {L}andau-{L}ifshitz equations:
              existence and nonuniqueness},
   JOURNAL = {Nonlinear Anal.},
  FJOURNAL = {Nonlinear Analysis. Theory, Methods \& Applications. An
              International Multidisciplinary Journal},
    VOLUME = {18},
      YEAR = {1992},
    NUMBER = {11},
     PAGES = {1071--1084},
      ISSN = {0362-546X},
   MRCLASS = {35Q99 (35D05 58E20 82D40)},
  MRNUMBER = {1167422},
MRREVIEWER = {Kotik K. Lee},
       DOI = {10.1016/0362-546X(92)90196-L},
       URL = {https://doi.org/10.1016/0362-546X(92)90196-L},
}

@article {Alouges:2008,
    AUTHOR = {Alouges, Fran\c{c}ois},
     TITLE = {A new finite element scheme for {L}andau--{L}ifchitz equations},
   JOURNAL = {Discrete Contin. Dyn. Syst. Ser. S},
  FJOURNAL = {Discrete and Continuous Dynamical Systems. Series S},
    VOLUME = {1},
      YEAR = {2008},
    NUMBER = {2},
     PAGES = {187--196},
      ISSN = {1937-1632},
   MRCLASS = {65M60 (35K55 82B80 82D40)},
  MRNUMBER = {2379897},
MRREVIEWER = {Etienne Emmrich},
       DOI = {10.3934/dcdss.2008.1.187},
       URL = {https://doi.org/10.3934/dcdss.2008.1.187},
}

@article {An:2016,
    AUTHOR = {An, Rong},
     TITLE = {Optimal error estimates of linearized {C}rank--{N}icolson
              {G}alerkin method for {L}andau--{L}ifshitz equation},
   JOURNAL = {J. Sci. Comput.},
  FJOURNAL = {Journal of Scientific Computing},
    VOLUME = {69},
      YEAR = {2016},
    NUMBER = {1},
     PAGES = {1--27},
      ISSN = {0885-7474},
   MRCLASS = {65M60 (35K55 35Q60 65M15)},
  MRNUMBER = {3542789},
MRREVIEWER = {Dmitriy Leykekhman},
       DOI = {10.1007/s10915-016-0181-1},
       URL = {https://doi.org/10.1007/s10915-016-0181-1},
}

@phdthesis{Banas:2005,
	title        = {On dynamical {M}icromagnetism with {M}agnetostriction},
	author       = {Ba{\v{n}}as, {\v L}ubom\'{\i}r},
	year         = {2005},
	school       = {Ghent University},
	status       = {thesis}
}

@ARTICLE{Baensch:1991,
 AUTHOR={B{\"a}nsch, Eberhard},
 TITLE={Local mesh refinement in $2$ and $3$ dimensions},
 JOURNAL={Impact Comput. Sci. Engrg.},
 VOLUME=3,
 YEAR=1991,
 PAGES={181-191},
}

@article {BankYserentant:2015,
    AUTHOR = {Bank, Randolph E. and Yserentant, Harry},
     TITLE = {A note on interpolation, best approximation, and the
              saturation property},
   JOURNAL = {Numer. Math.},
  FJOURNAL = {Numerische Mathematik},
    VOLUME = {131},
      YEAR = {2015},
    NUMBER = {1},
     PAGES = {199--203},
      ISSN = {0029-599X},
   MRCLASS = {65N30 (41A40 41A50 65D05 65N15 65N50)},
  MRNUMBER = {3383332},
MRREVIEWER = {Andrea Bressan},
       DOI = {10.1007/s00211-014-0687-0},
       URL = {https://doi.org/10.1007/s00211-014-0687-0},
}

@article {BBP:2008,
    AUTHOR = {Ba\v{n}as, {\v L}ubom\'{\i}r and Bartels, S\"{o}ren and Prohl, Andreas},
     TITLE = {A convergent implicit finite element discretization of the
              {M}axwell--{L}andau--{L}ifshitz--{G}ilbert equation},
   JOURNAL = {SIAM J. Numer. Anal.},
  FJOURNAL = {SIAM Journal on Numerical Analysis},
    VOLUME = {46},
      YEAR = {2008},
    NUMBER = {3},
     PAGES = {1399--1422},
      ISSN = {0036-1429},
   MRCLASS = {65M60 (35Q60 94A08)},
  MRNUMBER = {2390999},
       DOI = {10.1137/070683064},
       URL = {https://doi.org/10.1137/070683064},
}

@article {BKP:2008,
    AUTHOR = {Bartels, S\"{o}ren and Ko, Joy and Prohl, Andreas},
     TITLE = {Numerical analysis of an explicit approximation scheme for the
              {L}andau-{L}ifshitz-{G}ilbert equation},
   JOURNAL = {Math. Comp.},
  FJOURNAL = {Mathematics of Computation},
    VOLUME = {77},
      YEAR = {2008},
    NUMBER = {262},
     PAGES = {773--788},
      ISSN = {0025-5718},
   MRCLASS = {82D40 (35K65 35Q60 65M12)},
  MRNUMBER = {2373179},
MRREVIEWER = {Yaniv Almog},
       DOI = {10.1090/S0025-5718-07-02079-0},
       URL = {https://doi.org/10.1090/S0025-5718-07-02079-0},
}

@article {BartelsProhl:2006,
    AUTHOR = {Bartels, S\"{o}ren and Prohl, Andreas},
     TITLE = {Convergence of an implicit finite element method for the
              {L}andau--{L}ifshitz--{G}ilbert equation},
   JOURNAL = {SIAM J. Numer. Anal.},
  FJOURNAL = {SIAM Journal on Numerical Analysis},
    VOLUME = {44},
      YEAR = {2006},
    NUMBER = {4},
     PAGES = {1405--1419},
      ISSN = {0036-1429},
   MRCLASS = {65M60 (82D40)},
  MRNUMBER = {2257110},
MRREVIEWER = {Anne Nouri},
       DOI = {10.1137/050631070},
       URL = {https://doi.org/10.1137/050631070},
}

@article {DGLL:2021,
    AUTHOR = {DeCaria, Victor and Guzel, Ahmet and Layton, William and Li, Yi},
     TITLE = {A variable stepsize, variable order family of low complexity},
   JOURNAL = {SIAM J. Sci. Comput.},
  FJOURNAL = {SIAM Journal on Scientific Computing},
    VOLUME = {43},
      YEAR = {2021},
    NUMBER = {3},
     PAGES = {A2130--A2160},
      ISSN = {1064-8275},
   MRCLASS = {65L06 (65L04 65M20)},
  MRNUMBER = {4271476},
       DOI = {10.1137/19M1258153},
       URL = {https://doi.org/10.1137/19M1258153},
}

@article {DPPRS:2020,
    AUTHOR = {Di Fratta, Giovanni and Pfeiler, Carl-Martin and Praetorius,
              Dirk and Ruggeri, Michele and Stiftner, Bernhard},
     TITLE = {Linear second-order {IMEX}-type integrator for the (eddy
              current) {L}andau-{L}ifshitz-{G}ilbert equation},
   JOURNAL = {IMA J. Numer. Anal.},
  FJOURNAL = {IMA Journal of Numerical Analysis},
    VOLUME = {40},
      YEAR = {2020},
    NUMBER = {4},
     PAGES = {2802--2838},
      ISSN = {0272-4979},
   MRCLASS = {65M60 (65M12 78M10)},
  MRNUMBER = {4167063},
       DOI = {10.1093/imanum/drz046},
       URL = {https://doi.org/10.1093/imanum/drz046},
}

@article {FeischlTran:2017_fembem,
    AUTHOR = {Feischl, Michael and Tran, Thanh},
     TITLE = {The eddy current--{LLG} equations: {FEM}-{BEM} coupling and a
              priori error estimates},
   JOURNAL = {SIAM J. Numer. Anal.},
  FJOURNAL = {SIAM Journal on Numerical Analysis},
    VOLUME = {55},
      YEAR = {2017},
    NUMBER = {4},
     PAGES = {1786--1819},
      ISSN = {0036-1429},
   MRCLASS = {65M60 (35Q41 35Q60 35R60 60H15 65M12 65M15 82D45)},
  MRNUMBER = {3679317},
       DOI = {10.1137/16M1065161},
       URL = {https://doi.org/10.1137/16M1065161},
}

@article {FeischlTran:2017_existence,
    AUTHOR = {Feischl, Michael and Tran, Thanh},
     TITLE = {Existence of regular solutions of the
              {L}andau--{L}ifshitz--{G}ilbert equation in 3{D} with natural
              boundary conditions},
   JOURNAL = {SIAM J. Math. Anal.},
  FJOURNAL = {SIAM Journal on Mathematical Analysis},
    VOLUME = {49},
      YEAR = {2017},
    NUMBER = {6},
     PAGES = {4470--4490},
      ISSN = {0036-1410},
   MRCLASS = {35Q60 (35K55 35R60 60H15 65M12 82D40)},
  MRNUMBER = {3723324},
MRREVIEWER = {Peter Bernard Weichman},
       DOI = {10.1137/16M1103427},
       URL = {https://doi.org/10.1137/16M1103427},
}

@article {F:2022,
    AUTHOR = {Feischl, Michael},
     TITLE = {Inf-sup stability implies quasi-orthogonality},
   JOURNAL = {Math. Comp.},
  FJOURNAL = {Mathematics of Computation},
    VOLUME = {91},
      YEAR = {2022},
    NUMBER = {337},
     PAGES = {2059--2094},
      ISSN = {0025-5718},
   MRCLASS = {65N30 (15A23 65N50)},
  MRNUMBER = {4451456},
       DOI = {10.1090/mcom/3748},
       URL = {https://doi.org/10.1090/mcom/3748},
}

@article {Grigorieff:1983,
    AUTHOR = {Grigorieff, Rolf D.},
     TITLE = {Stability of multistep-methods on variable grids},
   JOURNAL = {Numer. Math.},
  FJOURNAL = {Numerische Mathematik},
    VOLUME = {42},
      YEAR = {1983},
    NUMBER = {3},
     PAGES = {359--377},
      ISSN = {0029-599X},
   MRCLASS = {65L20 (65L05)},
  MRNUMBER = {723632},
MRREVIEWER = {F. J. Murray},
       DOI = {10.1007/BF01389580},
       URL = {https://doi.org/10.1007/BF01389580},
}

@book {HairerWanner:2010,
    AUTHOR = {Hairer, Ernst and Wanner, Gerhard},
     TITLE = {Solving ordinary differential equations. {II}},
    SERIES = {Springer Series in Computational Mathematics},
    VOLUME = {14},
      NOTE = {Stiff and differential-algebraic problems,
              Second revised edition, paperback},
 PUBLISHER = {Springer-Verlag, Berlin},
      YEAR = {2010},
     PAGES = {xvi+614},
      ISBN = {978-3-642-05220-0},
   MRCLASS = {65-02 (65Lxx)},
  MRNUMBER = {2657217},
       DOI = {10.1007/978-3-642-05221-7},
       URL = {https://doi.org/10.1007/978-3-642-05221-7},
}

@article {HEGP:2015,
    AUTHOR = {Hay, Alexander and Etienne, Stephane and Garon, Andr{\'e} and Pelletier, Dominique},
     TITLE = {Time-integration for {ALE} simulations of fluid-structure
              interaction problems: stepsize and order selection based on
              the {BDF}},
   JOURNAL = {Comput. Methods Appl. Mech. Engrg.},
  FJOURNAL = {Computer Methods in Applied Mechanics and Engineering},
    VOLUME = {295},
      YEAR = {2015},
     PAGES = {172--195},
      ISSN = {0045-7825,1879-2138},
   MRCLASS = {65M60 (74F10 76D05 76M10)},
  MRNUMBER = {3388830},
MRREVIEWER = {Shangerganesh\ Lingeshwaran},
       DOI = {10.1016/j.cma.2015.06.006},
       URL = {https://doi.org/10.1016/j.cma.2015.06.006},
}

@article {HPPRSS:2019,
    AUTHOR = {Hrkac, Gino and Pfeiler, Carl-Martin and Praetorius, Dirk and
              Ruggeri, Michele and Segatti, Antonio and Stiftner, Bernhard},
     TITLE = {Convergent tangent plane integrators for the simulation of
              chiral magnetic skyrmion dynamics},
   JOURNAL = {Adv. Comput. Math.},
  FJOURNAL = {Advances in Computational Mathematics},
    VOLUME = {45},
      YEAR = {2019},
    NUMBER = {3},
     PAGES = {1329--1368},
      ISSN = {1019-7168},
   MRCLASS = {65M60 (35K55 35Q60 65M12 78M10)},
  MRNUMBER = {3955721},
       DOI = {10.1007/s10444-019-09667-z},
       URL = {https://doi.org/10.1007/s10444-019-09667-z},
}

@article {KVBAT:2014,
    AUTHOR = {Kritsikis, Evaggelos and Vaysset, A. and Buda-Prejbeanu, Liliana D. and
              Alouges, Fran\c{c}ois and Toussaint, Jean-C.},
     TITLE = {Beyond first-order finite element schemes in micromagnetics},
   JOURNAL = {J. Comput. Phys.},
  FJOURNAL = {Journal of Computational Physics},
    VOLUME = {256},
      YEAR = {2014},
     PAGES = {357--366},
      ISSN = {0021-9991},
   MRCLASS = {65M60 (78A30 78A48 78M10 82D40)},
  MRNUMBER = {3117413},
MRREVIEWER = {Weimin Han},
       DOI = {10.1016/j.jcp.2013.08.035},
       URL = {https://doi.org/10.1016/j.jcp.2013.08.035},
}

@article {GearTu:1974,
    AUTHOR = {Gear, C. William and Tu, K. W.},
     TITLE = {The effect of variable mesh size on the stability of multistep
              methods},
   JOURNAL = {SIAM J. Numer. Anal.},
  FJOURNAL = {SIAM Journal on Numerical Analysis},
    VOLUME = {11},
      YEAR = {1974},
     PAGES = {1025--1043},
      ISSN = {0036-1429},
   MRCLASS = {65L05},
  MRNUMBER = {368436},
MRREVIEWER = {D.\ J.\ Evans},
       DOI = {10.1137/0711079},
       URL = {https://doi.org/10.1137/0711079},
}

@article {GearWatanabe:1974,
    AUTHOR = {Gear, C. William and Watanabe, Daniel S.},
     TITLE = {Stability and convergence of variable order multistep methods},
   JOURNAL = {SIAM J. Numer. Anal.},
  FJOURNAL = {SIAM Journal on Numerical Analysis},
    VOLUME = {11},
      YEAR = {1974},
     PAGES = {1044--1058},
      ISSN = {0036-1429},
   MRCLASS = {65L05},
  MRNUMBER = {368437},
MRREVIEWER = {D.\ J.\ Evans},
       DOI = {10.1137/0711080},
       URL = {https://doi.org/10.1137/0711080},
}

@article {PRS:2018,
    AUTHOR = {Praetorius, Dirk and Ruggeri, Michele and Stiftner, Bernhard},
     TITLE = {Convergence of an implicit-explicit midpoint scheme for
              computational micromagnetics},
   JOURNAL = {Comput. Math. Appl.},
  FJOURNAL = {Computers \& Mathematics with Applications. An International
              Journal},
    VOLUME = {75},
      YEAR = {2018},
    NUMBER = {5},
     PAGES = {1719--1738},
      ISSN = {0898-1221},
   MRCLASS = {65M60 (82D40)},
  MRNUMBER = {3766546},
       DOI = {10.1016/j.camwa.2017.11.028},
       URL = {https://doi.org/10.1016/j.camwa.2017.11.028},
}

@book {Prohl:2001,
    AUTHOR = {Prohl, Andreas},
     TITLE = {Computational micromagnetism},
    SERIES = {Advances in Numerical Mathematics},
 PUBLISHER = {B. G. Teubner, Stuttgart},
      YEAR = {2001},
     PAGES = {xviii+304},
      ISBN = {3-519-00358-9},
   MRCLASS = {82D40 (65N30 74F15 82-02 82C80 82D30)},
  MRNUMBER = {1885923},
MRREVIEWER = {Thomas P. Svobodny},
       DOI = {10.1007/978-3-663-09498-2},
       URL = {https://doi.org/10.1007/978-3-663-09498-2},
}

@article {ShampineReichelt:1997,
    AUTHOR = {Shampine, Lawrence F. and Reichelt, Mark W.},
     TITLE = {The {MATLAB} {ODE} suite},
      NOTE = {Dedicated to C. William Gear on the occasion of his 60th
              birthday},
   JOURNAL = {SIAM J. Sci. Comput.},
  FJOURNAL = {SIAM Journal on Scientific Computing},
    VOLUME = {18},
      YEAR = {1997},
    NUMBER = {1},
     PAGES = {1--22},
      ISSN = {1064-8275},
   MRCLASS = {65Y15 (34-04 34A50 65L05 65L06)},
  MRNUMBER = {1433374},
       DOI = {10.1137/S1064827594276424},
       URL = {https://doi.org/10.1137/S1064827594276424},
}

@article{CAC:2014,
  author   = {E. Alberdi Celaya and J. J. Anza Aguirrezabala and Panagiotis Chatzipantelidis},
  title    = {Implementation of an Adaptive {BDF2} Formula and Comparison with the {MATLAB} {Ode15s}},
  journal  = {Procedia Computer Science},
  year     = {2014},
  volume   = {29},
  pages    = {1014-1026},
  issn     = {1877-0509},
  note     = {2014 International Conference on Computational Science},
  doi      = {https://doi.org/10.1016/j.procs.2014.05.091},
  url      = {https://www.sciencedirect.com/science/article/pii/S1877050914002683},
}

@misc{AFP:2025,
      title={{BDF2-type integrator for Landau-Lifshitz-Gilbert equation in micromagnetics. Part I: Unconditional weak convergence to weak solutions}}, 
      author={Michele Aldé and Michael Feischl and Dirk Praetorius},
      year={2025},
      eprint={2511.22000},
      archivePrefix={arXiv},
      primaryClass={math.NA},
      url={https://arxiv.org/abs/2511.22000}, 
}

@article {Gao:2014,
    AUTHOR = {Gao, Huadong},
     TITLE = {Optimal error estimates of a linearized backward {E}uler {FEM}
              for the {L}andau--{L}ifshitz equation},
   JOURNAL = {SIAM J. Numer. Anal.},
  FJOURNAL = {SIAM Journal on Numerical Analysis},
    VOLUME = {52},
      YEAR = {2014},
    NUMBER = {5},
     PAGES = {2574--2593},
      ISSN = {0036-1429},
   MRCLASS = {65M60 (35K51 35Q60 65M12 78M12)},
  MRNUMBER = {3273326},
MRREVIEWER = {Sarangam Majumdar},
       DOI = {10.1137/130936476},
       URL = {https://doi.org/10.1137/130936476},
}

@article{Gilbert:1955,
 AUTHOR = {Gilbert, Thomas L.},
 TITLE = {A {L}agrangian formulation of the gyromagnetic equation of the
   magnetic field},
 FJOURNAL = {Physical Reviews},
 JOURNAL = {Phys. Rev.},
 VOLUME = {100},
 YEAR = {1955},
 PAGES = {1243--1255}
}

@book {Ainsworth:2000,
    AUTHOR = {Ainsworth, Mark and Oden, J. Tinsley},
     TITLE = {A posteriori error estimation in finite element analysis},
    SERIES = {Pure and Applied Mathematics (New York)},
 PUBLISHER = {Wiley-Interscience [John Wiley \& Sons], New York},
      YEAR = {2000},
     PAGES = {xx+240},
      ISBN = {0-471-29411-X},
   MRCLASS = {65-02 (65N15)},
  MRNUMBER = {1885308},
MRREVIEWER = {Ricardo G. Dur\'{a}n},
       DOI = {10.1002/9781118032824},
       URL = {https://doi.org/10.1002/9781118032824},
}

@article{LandauLifshitz:1935,
 AUTHOR = {Landau, Lew D. and Lifshitz, Jewgeni M.},
 TITLE = {On the theory of the dispersion of magnetic permeability in
   ferromagnetic bodies},
 JOURNAL = {Phys. Z. Sowjetunion},
 VOLUME = {8},
 YEAR = {1935},
 PAGES = {153--168}
}

@article {LTZ:2020,
    AUTHOR = {Liao, Hong-lin and Tang, Tao and Zhou, Tao},
     TITLE = {On energy stable, maximum-principle preserving, second-order
              {BDF} scheme with variable steps for the {Allen}--{Cahn}
              equation},
   JOURNAL = {SIAM J. Numer. Anal.},
  FJOURNAL = {SIAM Journal on Numerical Analysis},
    VOLUME = {58},
      YEAR = {2020},
    NUMBER = {4},
     PAGES = {2294--2314},
      ISSN = {0036-1429},
   MRCLASS = {65M60 (35B25 35J20 35K58 65M12)},
  MRNUMBER = {4134034},
MRREVIEWER = {Mohammad Asadzadeh},
       DOI = {10.1137/19M1289157},
       URL = {https://doi.org/10.1137/19M1289157},
}

@article {Schrefl:1999,
 Author  = {Schrefl, Thomas},
 Title   = {Finite elements in numerical micromagnetics. {Part I}: {Granular} hard magnets},
 Journal = {Journal of Magnetism and Magnetic Materials},
 Year    = {1999},
 Number  = {1},
 Pages   = {45--65},
 Volume  = {207},
 ISSN    = {0304-8853},
 Doi     = {10.1016/S0304-8853(99)00532-6},
 Url     = {https://doi.org/10.1016/S0304-8853(99)00532-6},
}

@article {DDPR:2022,
    AUTHOR = {Davoli, Elisa and Di Fratta, Giovanni and Praetorius, Dirk and
      Ruggeri, Michele},
     TITLE = {Micromagnetics of thin films in the presence of
      {Dzyaloshinskii--Moriya} interaction},
   JOURNAL = {Math. Models Methods Appl. Sci.},
  FJOURNAL = {Mathematical Models and Methods in Applied Sciences},
    VOLUME = {32},
      YEAR = {2022},
    NUMBER = {5},
     PAGES = {911-939},
      ISSN = {---},
   MRCLASS = {---},
  MRNUMBER = {---},
MRREVIEWER = {---},
       DOI = {https://doi.org/10.1142/S0218202522500208},
       URL = {https://doi.org/10.1142/S0218202522500208},
}

@article{Davis:2004,
author = {Davis, Timothy A.},
title = {Algorithm 832: {UMFPACK} {V}4.3---an unsymmetric-pattern multifrontal method},
journal = {ACM Transactions on Mathematical Software},
year = {2004},
publisher = {Association for Computing Machinery},
address = {New York, NY, USA},
volume = {30},
number = {2},
issn = {0098-3500},
url = {https://doi.org/10.1145/992200.992206},
doi = {10.1145/992200.992206},
pages = {196–199},
numpages = {4}
}

@Article{BangerthET:2023,
  title   = {The \texttt{deal.II} Library, Version 9.5},
  author  = {Daniel Arndt and Wolfgang Bangerth and Maximilian Bergbauer and
             Marco Feder and Marc Fehling and Johannes Heinz and
             Timo Heister and Luca Heltai and Martin Kronbichler and
             Matthias Maier and Peter Munch and Jean-Paul Pelteret and
             Bruno Turcksin and David Wells and Stefano Zampini},
  journal = {Journal of Numerical Mathematics},
  year    = {2023},
  doi     = {10.1515/jnma-2023-0089},
  pages   = {231--246},
  volume  = {31},
  number  = {3},
  url     = {https://dealii.org/deal95-preprint.pdf}
}

@article {Nishikawa:2019,
    AUTHOR = {Nishikawa, Hiroaki},
     TITLE = {On large start-up error of {BDF}2},
   JOURNAL = {J. Comput. Phys.},
  FJOURNAL = {Journal of Computational Physics},
    VOLUME = {392},
      YEAR = {2019},
     PAGES = {456--461},
      ISSN = {0021-9991,1090-2716},
   MRCLASS = {65L06 (65L20)},
  MRNUMBER = {3948731},
MRREVIEWER = {Emanuele\ Galligani},
       DOI = {10.1016/j.jcp.2019.04.070},
       URL = {https://doi.org/10.1016/j.jcp.2019.04.070},
}
%%%%%%%%%%%%%%%%%%%%%%%%%%%%%%%%%%%%%%%%%%%%%%%%%%%%%%%%%%%%%%%%%%%%%%%
\end{document}